\documentclass{amsart}

\usepackage{amssymb,amscd,amsthm,latexsym}
\usepackage[all]{xy}

\usepackage[pagebackref, colorlinks, citecolor=red, linkcolor=red]{hyperref}

\def\Aff{ \mathrm{Aff} }
\def\Ab{ \mathrm{Ab} }

\def\Pt{ \mathrm{Pt} }
\def\Nil{ \mathrm{Nil} }

\def\Fld{ \mathrm{Fld} }
\def\Gr{ \mathrm{Grp} }

\def\TG{ \mathrm{TriGrp}_\star }

\def\cleft{\hbox{[\kern-.16em\hbox{[}}}
\def\cright{\hbox{]\kern-.16em\hbox{]}}}

\def\CC{ \mathbb{C} }
\def\DD{ \mathbb{D} }
\def\EE{ \mathbb{E} }

\def\ZZ{ \mathbb{Z} }

\def\inc{ \hookrightarrow }
\def\into{ \rightarrowtail }
\def\onto{ \twoheadrightarrow }
\def\dto{ \rightrightarrows }
\def\lto{ \longrightarrow }
\def\otl{ \longleftarrow }
\def\lrto{ \leftrightarrows }

\def\eps{ \epsilon }
\def\Sg{ \mathfrak{S} }
\def\sg{ \sigma }
\def\OO{ \mathbb{O} }
\def\Ff{ \mathcal{F} }

\def\cleft{\hbox{[\kern-.16em\hbox{[}}}
\def\cright{\hbox{]\kern-.16em\hbox{]}}}

\theoremstyle{plain}
\newtheorem{dfn}{Definition}[section]
\newtheorem{thm}[dfn]{Theorem}
\newtheorem{lma}[dfn]{Lemma}
\newtheorem{prp}[dfn]{Proposition}
\newtheorem{cor}[dfn]{Corollary}

\theoremstyle{remark}
\newtheorem{rmk}[dfn]{Remark}

\begin{document}

\title{Central reflections and nilpotency in exact Mal'tsev categories}

\begin{abstract}We study nilpotency in the context of exact Mal'tsev categories taking central extensions as the primitive notion. This yields a nilpotency tower which is analysed from the perspective of Goodwillie's functor calculus.

We show in particular that the reflection into the subcategory of $n$-nilpotent objects is the universal endofunctor of degree $n$ if and only if every $n$-nilpotent object is $n$-folded. In the special context of a semi-abelian category, an object is $n$-folded precisely when its Higgins commutator of length $n+1$ vanishes.\end{abstract}

\author{Clemens Berger}
\address{Universit\'{e} de Nice, Lab. J. A. Dieudonn\'{e}, Parc Valrose, 06108 Nice Cedex, France}
\email{cberger@math.unice.fr}
\author{Dominique Bourn}
\address{Lab. J. Liouville, CNRS Fed. Rech. 2956, Calais, France}
\email{bourn@univ-littoral.fr}
\date{November 25, 2016}
\subjclass{17B30, 18D35, 18E10, 18G50, 20F18}
\keywords{Nilpotency, Mal'tsev category, central extension, Goodwillie functor calculus}

\maketitle

\tableofcontents

\section*{Introduction}

This text investigates \emph{nilpotency} in the context of \emph{exact Mal'tsev categories}. Our purpose is twofold: basic phenomena of nilpotency are treated through universal properties rather than through commutator calculus, emphasising the fundamental role played by central extensions; nilpotency is then linked to an algebraic form of \emph{Goodwillie's functor calculus} \cite{G2}. This leads to a global understanding of nilpotency in terms of functors with bounded degree.

A \emph{Mal'tsev category} is a finitely complete category in which reflexive relations are equivalence relations \cite{CKP, CLP}. Important examples of exact Mal'tsev categories are \emph{Mal'tsev varieties} \cite{Ma} and \emph{semi-abelian categories} \cite{JMT}. The simplicial objects of an exact Mal'tsev category are ``internal'' Kan complexes (cf. \cite{CKP,Tim}). %In a varietal context, this property implies every simplicial object has an underlying simplicial set fulfilling Kan's extension condition \cite{CKP}.

Nilpotency is classically understood via the vanishing of iterated commutators: in a Mal'tsev variety by means of so-called Smith commutators \cite{Sm, FM}, in a semi-abelian category by means of so-called Huq commutators \cite{Huq,EV}. The first aim of this text is to promote another point of view which seems more intrinsic to us and is based on the notion of \emph{central extension}, by which we mean a regular epimorphism with central kernel relation. The \emph{$n$-nilpotent objects} are defined to be those which can be linked to a terminal object by a chain of $n$ consecutive central extensions. This notion of nilpotency is equivalent to the two aforementioned notions in their respective contexts (cf. Proposition \ref{iteratedSmith}). In particular, we get the usual notions of $n$-nilpotent group, $n$-nilpotent Lie algebra and $n$-nilpotent loop \cite{Br}. A category is called $n$-nilpotent if all its objects are $n$-nilpotent.

For any exact Mal'tsev category with binary sums, the full subcategory spanned by the $n$-nilpotent objects is a \emph{reflective Birkhoff subcategory} (cf. Theorem \ref{noBirkhoff}). This generalises the analogous known results for Mal'tsev varieties \cite{Sm,FM} and semi-abelian categories \cite{EV}. We denote the reflection into the subcategory of $n$-nilpotent objects by $I^n$ and the unit of the adjunction at an object $X$ by $\eta^n_X:X\onto I^n(X)$.

Since an $n$-nilpotent object is a fortiori $(n+1)$-nilpotent, the different reflections assemble into the following \emph{nilpotency tower}

$$ \xymatrix@=15pt{
&  & X \ar@{->>}[lldd] \ar@{->>}[ldd]^>>>>>>>>{\!\!\eta_X^1} \ar@{->>}[dd]^>>>>>{\eta_X^2}  \ar@{->>}[ddr]^>>>>{\eta_X^n} \ar@{->>}[ddrr]^>>>>{\eta_X^{n+1}}&& \\
&&&&\\
{\star\;}\ar@{<<-}[r]  & {I^1(X)\;}\ar@{<<-}[r]  & {I^2(X)\;}\ar@{<<.}[r]  & {I^n(X)\;}\ar@{<<-}[r]  & {I^{n+1}(X)\;} \ar@{<<.}[r]  &
                  }$$
in which the successive quotient maps $I^{n+1}(X)\onto I^n(X)$ are central extensions.

Pointed categories with binary sums will be ubiquitous throughout the text; we call them \emph{$\sg$-pointed} for short. Among $\sg$-pointed exact Mal'tsev categories we characterise the $n$-nilpotent ones as those for which the comparison maps $$\theta_{X,Y}:X+Y\onto X\times Y$$ are $(n-1)$-fold central extensions (cf. Theorem \ref{ndiscrepancy}). The nilpotency class of a $\sg$-pointed exact Mal'tsev category measures thus the discrepancy between binary sum and binary product.  If $n=1$, binary sum and binary product coincide, and all objects are abelian group objects. A $\sg$-pointed exact Mal'tsev category is $1$-nilpotent if and only if it is an abelian category (cf. Corollary \ref{abelian}). The unit $\eta_X^1:X\onto I^1(X)$ of the first Birkhoff reflection is \emph{abelianisation} (cf. Proposition \ref{abelianization}). Moreover, the successive kernels of the nilpotency tower are abelian group objects as well. This situation is reminiscent of what happens in Goodwillie's functor calculus \cite{G2} where ``infinite loop spaces'' play the role of abelian group objects. The second aim of our study of nilpotency was to get a deeper understanding of this analogy.

Goodwillie's notions \cite{G2} of \emph{cross-effect} and \emph{degree} of a functor translate well into our algebraic setting: for each $(n+1)$-tuple $(X_1,\dots,X_{n+1})$ of objects of a $\sg$-pointed category and each based endofunctor $F$ we define an $(n+1)$-cube $\Xi^F_{X_1,\dots,X_{n+1}}$ consisting of the images $F(X_{i_1}+\cdots+X_{i_k})$ for all subsequences of $(X_1,\dots,X_{n+1})$ together with the obvious contraction maps. We say that a functor $F$ is \emph{of degree $\leq n$} if these $(n+1)$-cubes are limit-cubes for all choices of $(X_1,\dots,X_{n+1})$.

We denote by $\theta^F_{X_1,\dots,X_{n+1}}:F(X_1+\cdots+X_{n+1})\to P^F_{X_1,\dots,X_{n+1}}$ the comparison map towards the limit of the punctured $(n+1)$-cube so that $F$ is of degree $\leq n$ if and only if $\theta_{X_1,\dots,X_{n+1}}^F$ is invertible for each choice of $(n+1)$-tuple. The kernel of $\theta^F_{X_1,\dots,X_{n+1}}$ is a \emph{$(n+1)$-st cross-effect} of $F$, denoted $cr^F_{n+1}(X_1,\dots,X_{n+1})$.

A based endofunctor $F$ is \emph{linear}, i.e. of degree $\leq 1$, if and only if $F$ takes binary sums to binary products. In a semi-abelian category, the second cross-effects $cr_2^F(X,Y)$ measure thus the failure of linearity of $F$. If $F$ is the identity functor, we drop $F$ from the notation so that $cr_2(X,Y)$ denotes the kernel of the comparison map $\theta_{X,Y}:X+Y\to X\times Y$. This kernel is often denoted $X\diamond Y$ and called the \emph{co-smash product} of $X$ and $Y$ (cf. \cite{CJ} and Remarks \ref{diamond} and \ref{history}).

%Goodwillie \cite{G2} defined for endofunctors of the category of based spaces a tower of approximating endofunctors such that the successive fibres of this tower take values in infinite loop spaces. Infinite loop spaces play a similar role in the homotopy theory of based spaces, as abelian group objects in pointed Mal'tsev categories. The Goodwillie tower of the identity functor of the category of based spaces is thus a homotopical analog of the nilpotency tower of a $\sg$-pointed exact Mal'tsev category.

An endofunctor of a semi-abelian (or homological \cite{BB}) category is of degree $\leq n$ if and only if all its cross-effects of order $n+1$ vanish. For functors taking values in abelian categories, our cross-effects agree with the original cross-effects of Eilenberg-Mac Lane \cite{EM} (cf. Remark \ref{history}). For functors taking values in $\sg$-pointed categories with pullbacks, our cross-effects agree with those of Hartl-Loiseau \cite{HL} and Hartl-Van der Linden \cite{HV}, defined as kernel intersections (cf. Definition \ref{ncube}).

A Goodwillie type characterisation of the nilpotency tower amounts to the property that for each $n$, the reflection $I^n$ into the Birkhoff subcategory of $n$-nilpotent objects is the \emph{universal endofunctor of degree $\leq n$}. In fact, every endofunctor of degree $\leq n$ of a $\sg$-pointed exact Mal'tsev category takes values in $n$-nilpotent objects (cf. Proposition \ref{nadditive}). The reflection $I^n$ is of degree $\leq n$ if and only if the identity functor of the Birkhoff subcategory of $n$-nilpotent objects itself is of degree $\leq n$. In the present article we have mainly been investigating this last property.

The property holds for $n=1$ because the identity functor of an abelian category is linear. However, already for $n=2$, there are examples of $2$-nilpotent semi-abelian categories which are not \emph{quadratic}, i.e. do not have an identity functor of degree $\leq 2$ (cf. Section \ref{Moufang}). We show that a $\sg$-pointed exact Mal'tsev category is quadratic if and only if the category is $2$-nilpotent \emph{and} algebraically distributive, i.e. endowed with isomorphisms $(X\times Z)+_Z(Y\times Z)\cong(X+Y)\times Z$ for all objects $X,Y,Z$ (cf. Corollary \ref{parexemple}). Since algebraic distributivity is preserved under Birkhoff reflection, the subcategory of $2$-nilpotent objects of an algebraically distributive exact Mal'tsev category is always quadratic (cf. Theorem \ref{acc}).

Algebraic distributivity is a consequence of the \emph{existence of centralisers for subobjects} as shown by Gray and the second author \cite{BGr}. For pointed Mal'tsev categories, it also follows from \emph{algebraic coherence} in the sense of Cigoli-Gray-Van der Linden \cite{CGV}. Our quadraticity result implies that iterated Huq commutator $[X,[X,X]]$ and ternary Higgins commutator $[X,X,X]$ coincide for each object $X$ of an algebraically distributive semi-abelian category (cf. Corollary \ref{Huq=Higgins} and \cite[Corollary 7.2]{CGV}).

There is a remarkable duality for $\sg$-pointed $2$-nilpotent exact Mal'tsev categories: algebraic distributivity amounts to algebraic codistributivity, i.e. to isomorphisms $(X\times Y)+Z\cong(X+Z)\times_Z(Y+Z)$ for all $X,Y,Z$ (cf. Proposition \ref{tau=iota}). Indeed, the difference between $2$-nilpotency and quadraticity is precisely algebraic \emph{co}distributivity (cf. Theorem \ref{quad}). An extension of this duality to all $n\geq 2$ is crucial in relating general nilpotency to identity functors with bounded degree.

The following characterisation is very useful: The identity functor of a $\sg$-pointed exact Mal'tsev category $\EE$ is of degree $\leq n$ if and only if all its objects are $n$-folded (cf. Proposition \ref{nadditive}). An object is \emph{$n$-folded} (cf. Definition \ref{nfolded}) if the $(n+1)$-st \emph{folding map} $\delta^X_{n+1}:X+\cdots+X\to X$ factors through the comparison map $\theta_{X,\dots,X}:X+\cdots+X\onto P_{X,\dots,X}$. In a varietal context this can be expressed in combinatorial terms (cf. Remark \ref{nfoldedgroup}). The full subcategory $\Fld^n(\EE)$ spanned by the $n$-folded objects is a reflective Birkhoff subcategory of $\EE$, and the reflection $J^n:\EE\to\Fld^n(\EE)$ is the universal endofunctor of degree $\leq n$ (cf. Theorem \ref{additivereflection}). Every $n$-folded object is $n$-nilpotent (cf. Proposition \ref{degn>nnil}) while the converse holds if and only if the other Birkhoff reflection $I^n:\EE\to\Nil^n(\EE)$ is also of degree $\leq n$.

In the context of semi-abelian categories, closely related results appear in the work of Hartl and his coauthors \cite{HL, HV, HVe}, although formulated slightly differently. In a semi-abelian category, an object $X$ is $n$-folded if and only if its \emph{Higgins commutator} of length $n+1$ vanishes (cf. Remark \ref{nfoldedgroup}), where the latter is defined as the image of the composite map $cr_{n+1}(X,\dots,X)\to X+\cdots+X\to X$, cf. \cite{HL, HV, MM}. The universal $n$-folded quotient $J^n(X)$ may then be identified with the quotient of $X$ by the Higgins commutator of length $n+1$ in much the same way as the universal $n$-nilpotent quotient $I^n(X)$ may be identified with the quotient of $X$ by the iterated Huq commutator of length $n+1$. It was Hartl's insight that Higgins commutators are convenient for extending the ``polynomial functors'' of Eilenberg-Mac Lane \cite{EM} to a semi-abelian context. Our treatment in the broader context of exact Mal'tsev categories follows more closely Goodwillie's functor calculus \cite{G2}.

In a $\sg$-pointed exact Mal'cev category, abelianisation $I^1$ is the universal endofunctor $J^1$ of degree $\leq 1$ (cf. Mantovani-Metere \cite{MM}). For $n>1$ however, the universal endofunctor $J^n$ of degree $\leq n$ is in general a proper quotient of the $n$-th Birkhoff reflection $I^n$ (cf. Corollary \ref{inclusion}). In order to show that even in a semi-abelian variety the two endofunctors may disagree, we exhibit a Moufang loop of order $16$ (subloop of Cayley's octonions) which is $2$-nilpotent but not $2$-folded (cf. Section \ref{Moufang}). Alternatively, Mostovoy's modified lower central series of a loop \cite{Mos} yields other examples of a similar kind provided the latter agrees with the successive Higgins commutators of the loop (cf. \cite[Example 2.15]{HV} and \cite{SV}).

We did not find a simple categorical structure that would entail the equivalence between $n$-nilpotency and $n$-foldedness for all $n$. As a first step in this direction we show that an $n$-nilpotent semi-abelian category has an identity functor of degree $\leq n$ if and only if its $n$-th cross-effect is \emph{multilinear} (cf. Theorem \ref{multilinear}). We also show that the nilpotency tower has the desired universal property if and only if it is \emph{homogeneous}, i.e. for each $n$, the $n$-th kernel functor is of degree $\leq n$ (cf. Theorem \ref{folklore}). This is preserved under Birkhoff reflection (cf. Theorem \ref{stableHuq=Higgins}). The categories of groups and of Lie algebras have homogeneous nilpotency towers so that a group, resp. Lie algebra is $n$-nilpotent if and only if it is $n$-folded, and the Birkhoff reflection $I^n$ is here indeed the universal endofunctor of degree $\leq n$. The category of \emph{triality groups} \cite{D,G,Ha} also has a homogeneous nilpotency tower although it contains the category of Moufang loops as a full coreflective subcategory, and the latter has an inhomogeneous nilpotency tower (cf. Section \ref{Moufang}).\vspace{1ex}

There are several further ideas closely related to the contents of this article which we hope to address in future work. Let us mention two of them:

The associated graded object of the nilpotency tower $\oplus_{n\geq 1} K[I^n(X)\onto I^{n-1}(X)]$ is a functor in $X$ taking values in graded abelian group objects. For the category of groups this functor actually takes values in graded Lie rings and as such preserves $n$-nilpotent objects and free objects, cf. Lazard \cite{La}. It is likely that for a large class of semi-abelian categories, the associated graded object of the nilpotency tower carries a similar algebraic structure. It would be interesting to establish the relationship between this algebraic structure and the cross-effects of the identity functor.

It follows from \cite[Theorem 4.2]{CKP} and \cite[Theorem IV.4]{Q} that the simplicial objects of a pointed Mal'tsev variety carry a Quillen model structure in which the weak equivalences are the maps inducing a quasi-isomorphism on Moore complexes. Such a model structure also exists for the simplicial objects of a semi-abelian category with enough projectives, cf. \cite{Q,Tim}. In both cases, regular epimorphisms are fibrations, and the trivial fibrations are precisely the regular epimorphisms for which the kernel is homotopically trivial. This implies that Goodwillie's homotopical cross-effects \cite{G2} agree here with our algebraic cross-effects.

Several notions of \emph{homotopical nilpotency} are now available. The first is the least integer $n$ for which the unit $\eta^n_{X_\bullet}:X_\bullet\onto I^n(X_\bullet)$ is a trivial fibration, the second (resp. third) is the least integer $n$ for which $X_\bullet$ is homotopically $n$-folded (resp. the value of an $n$-excisive approximation of the identity). The first is a lower bound for the second, and the second is a lower bound for the third invariant. For \emph{simplicial groups} the first invariant recovers the Berstein-Ganea nilpotency for loop spaces \cite{BG}, the second the cocategory of Hovey \cite{Ho}, and the third the Biedermann-Dwyer nilpotency for homotopy nilpotent groups \cite{BD}. Similar chains of inequalities have recently been studied by Eldred \cite{E} and Costoya-Scherer-Viruel \cite{CSV}.\vspace{1ex}

The plan of this article is as follows.\vspace{1ex}

Section \ref{GL} reviews the notions of central extension and regular pushout. At the end an algebraic Beck-Chevalley condition for pushouts of regular epimorphisms in an exact Mal'tsev category is established.

Section \ref{DF} presents our definition of nilpotency and studies under which conditions the $n$-nilpotent objects form a reflective Birkhoff subcategory.

Section \ref{section3} investigates central reflections, the motivating example being the reflection of the category of $(n+1)$-nilpotent objects into the category of $n$-nilpotent objects. The unit of these central reflections is shown to be pointwise affine.

Section \ref{section4} establishes first aspects of nilpotency. The nilpotency class of a $\sg$-pointed exact Mal'cev category is related to universal properties of the comparison map $\theta_{X,Y}:X+Y\to X\times Y$. This leads to a new family of binary tensor products interpolating between binary sum and binary product.

Section \ref{linearquadratic} studies the $\sg$-pointed exact Mal'tsev categories with quadratic identity functor. They are characterised among the $2$-nilpotent ones as those which are algebraically distributive, resp. algebraically codistributive.

Section \ref{degreen} studies the $\sg$-pointed exact Mal'tsev categories with an identity functor of degree $\leq n$. They are characterised as those in which all objects are $n$-folded. Every $n$-folded object is shown to be $n$-nilpotent. Several sufficient criteria for the converse are given. The semi-abelian varieties of groups, Lie algebras, Moufang loops and triality groups are discussed.\vspace{2ex}

\section{Central extensions and regular pushouts}\label{GL}

In this introductory section we review the notion of \emph{central equivalence relation} and study basic properties of the associated class of \emph{central extensions}, needed for our treatment of nilpotency. By central extension we mean a regular epimorphism with central kernel relation \cite{Bourn4,BG2,BG3}. This algebraic concept of central extension has to be distinguished from the axiomatic concept of Janelidze-Kelly \cite{JK} which is based on a previously chosen admissible Birkhoff subcategory. Nevertheless, it is known that with respect to the Birkhoff subcategory of \emph{abelian group objects}, the two approaches yield the same class of central extensions in any congruence modular variety (cf. \cite{JK2,JK3}) as well as in any exact Mal'tsev category (cf. \cite{BG2,EV2,Gran2}).

We assume throughout that our ambient category is a \emph{Mal'tsev category}, i.e. a finitely complete category in which every reflexive relation is an equivalence relation, cf. \cite{BB,Bourn0,CKP,CLP}. Most of the material of this section is well-known to the expert, and treated in some detail here mainly to fix notation and terminology.

One exception is Section \ref{Beck-Chevalley} which establishes an ``algebraic'' \emph{Beck-Chevalley condition} for pushouts of regular epimorphisms in exact Mal'tsev categories, dual to the familiar Beck-Chevalley condition for pullbacks of monomorphisms in elementary toposes. In recent and independent work, Gran-Rodelo \cite{GR} consider a weaker form of this condition and show that it characterises regular Goursat categories.

\subsection{Smith commutator of equivalence relations}\label{Smithcommutator}--\vspace{1ex}

An \emph{equivalence relation} $R$ on $X$ will be denoted as a reflexive graph $(p_0,p_1):R\dto X$ with section $s_0:X\to R$, but whenever convenient we shall consider $R$ as a subobject of $X\times X$. By a \emph{fibrant} map of equivalence relations $(X,R)\to(Y,S)$ we mean a natural transformation of the underlying reflexive graphs such that the three naturality squares are pullback squares.

A particularly important equivalence relation is the \emph{kernel relation} $R[f]$ of a morphism $f:X\to Y$ which is part of the following diagram:
 $$ \xymatrix@=25pt{
         R[f]   \ar@<-1,ex>[r]_{p_0}\ar@<+1,ex>[r]^{p_1} & X\ar[l]\ar[r]^{f} & Y.
                        }
     $$

The \emph{discrete} equivalence relation $\Delta_X$ on $X$ is the kernel relation $R[1_X]$ of the identity map $1_X:X\to X$. The \emph{indiscrete} equivalence relation $\nabla_X$ on $X$ is the kernel relation $R[\omega_X]$ of the unique map $\omega_X$ from $X$ to a terminal object.

Two equivalence relations $R,S$ on the same object $X$ are said to \emph{centralise each other} if the square
$$ \xymatrix@=25pt{
      R\times_XS \ar@{<-}[r]^>>>>>{(s_0^R,1_S)} \ar@{<-}[d]_{(1_R,s_0^S)}\ar@{.>}[rd]^p &  S\ar[d]^{p_1^S}\\
      R \ar@{->}[r]_{p_0^R} & X
                   }
$$
admits a (necessarily unique) filler which makes the diagram commute, cf. \cite{BG3,Pe}.

In set-theoretical terms such a filler amounts to the existence of a \emph{``partial Mal'tsev operation''} on $X$, namely (considering $R\times_XS$ as a subobject of $X\times X\times X$) a ternary operation $p:R\times_XS\to X$ such that $x=p(x,y,y)$ and $p(x,x,y)=y$. We shall follow Marino Gran and the second author in calling $p:R\times_XS\to X$ a \emph{connector} between $R$ and $S$, cf. \cite{BG2,BG3}.

In a finitely cocomplete regular Mal'tsev category, there exists for each pair $(R,S)$ of equivalence relations on $X$ a smallest effective equivalence relation $[R,S]$ on $X$ such that $R$ and $S$ centralise each other in the quotient $X/[R,S]$. This equivalence relation is the so-called \emph{Smith commutator} of $R$ and $S$, cf. \cite{Bourn4,BG3,Pe,Sm}.

In these terms $R$ and $S$ centralise each other precisely when $[R,S]=\Delta_X$. The Smith commutator is monotone in each variable and satisfies $$[R,S]=[S,R]\text{ and }f([R,S])\subset[f(R),f(S)]$$ for each regular epimorphism $f:X\to Y$, where $f(R)$ denotes the direct image of the subobject $R\subset X\times X$ under the regular epimorphism $f\times f:X\times X\to Y\times Y$. The Mal'tsev condition implies that this direct image represents an equivalence relation on $Y$. Note that equality $f([R,S])=[f(R),f(S)]$ holds if and only if the direct image $f([R,S])$ is an \emph{effective} equivalence relation on $Y$, which is always the case in an \emph{exact} Mal'tsev category.

\subsection{Central equivalence relations and central extensions}\label{centralextensiondef}An equivalence relation $R$ on $X$ is said to be \emph{central} if $[R,\nabla_X]=\Delta_X$. A \emph{central extension} is by definition a regular epimorphism with central kernel relation. An \emph{$n$-fold central extension} is the composite of $n$ central extensions. An \emph{$n$-fold centrally decomposable morphism} is the composite of $n$ morphisms with central kernel relation.

The indiscrete equivalence relation $\nabla_X$ is a central equivalence relation precisely when $X$ admits an \emph{internal Mal'tsev operation} $p:X\times X\times X\to X$. In pointed Mal'tsev categories such a Mal'tsev operation amounts to an \emph{abelian group} structure on $X$, cf. \cite[Proposition 2.3.8]{BB}. An object $X$ of a pointed Mal'tsev category $(\DD,\star_\DD)$ is thus an abelian group object if and only if the map $X\to\star_\DD$ is a central extension.

Central equivalence relations are closed under binary products and inverse image along monomorphisms. In a regular Mal'tsev category, central equivalence relations are closed under direct images, cf. \cite[Proposition 4.2]{BG3} and  \cite[Proposition 2.6.15]{BB}.

\begin{lma}\label{relativecentral}In a regular Mal'tsev category, an $n$-fold centrally decomposable morphism can be written as an $n$-fold central extension followed by a monomorphism.\end{lma}

\proof It suffices to show that a monomorphism $\psi$ followed by a central extension $\phi$ can be rewritten as a central extension $\phi'$ followed by a monomorphism $\psi'$. Indeed, the kernel relation $R[\phi\psi]$ is central, being the restriction $\psi^{-1}(R[\phi])$ of the central equivalence relation $R[\phi]$ along the monomorphism $\psi$; therefore, by regularity, one obtains $\phi\psi=\psi'\phi'$ where $\phi'$ is quotienting by the kernel relation $R[\phi\psi]$.\endproof

\begin{lma}In a Mal'tsev category, morphisms with central kernel relation are closed under pullback. In a regular Mal'tsev category, central extensions are closed under pullback.\end{lma}

\proof It suffices to show the first statement. In the following diagram,

$$\xymatrix@=20pt{
        R[f']  \ar@{->}[d]_{R(x,y)}\ar@<-1,ex>[r]_{p'_0}\ar@<+1,ex>[r]^{p'_1} & X' \ar@{->}[d]^x\ar[l]\ar@{->}[r]^{f'} & Y' \ar@{->}[d]^y\\
        R[f] \ar@<-1,ex>[r]_{p_0} \ar@<+1,ex>[r]^{p_1}  & X \ar[l] \ar@{->}[r]_f & Y
                       }
  $$if the right square is a pullback, then the left square is a fibrant map of equivalence relations. This permits to lift the connector $p:R[f]\times_X\nabla_X\to X$ so as to obtain a connector $p':R[f']\times_{X'}\nabla_{X'}\to X'$.\endproof

\begin{lma}\label{lefterase}Let $X\overset{f}{\to}Y\overset{g}{\to}Z$ be morphisms in a Mal'tsev category.

If $\,gf$ is a morphism with central kernel relation then so is $f$. More generally, if $\,gf$ is $n$-fold centrally decomposable then so is $f$.
\end{lma}
\proof
Since $R[f]\subset R[gf]$, the commutation relation $[R[gf],\nabla_X]=\Delta_X$ implies the commutation relation $[R[f],\nabla_X]=\Delta_X$. Assume now $gf=k_n\cdots k_1$ where each $k_i$ is a morphism with central kernel relation. In the following pullback
$$ \xymatrix@=25pt{
      P  \ar[d]_{\psi} \ar[r]^{\gamma} &  X \ar[d]^{gf} \ar@<1ex>@{.>}[l]^{\phi}\\
      Y \ar[r]_{g} & Z
                   }
$$
$\phi$ is the unique map such that $\gamma\phi=1_X$ and $\psi\phi=f$. If we denote by $h_i$ the morphism with central kernel relation obtained by pulling back $k_i$ along $g$, we get $f=h_n\cdots h_2(h_1\phi)$. Since $\phi$ is a monomorphism, the kernel relation $R[h_1\phi]=\phi^{-1}(R[h_1])$ is central, and hence $f$ is the composite of $n$ morphisms with central kernel relation.\endproof

\begin{prp}[Corollary 3.3 in \cite{Bourn4}]\label{univcentral}In a finitely cocomplete regular Mal'tsev category, each morphism $f:X\rightarrow Y$ factors canonically as in
$$ \xymatrix@=15pt{
      X \ar@{->>}[rr]^{\eta_f} \ar[d]_{f} & &X/[\nabla_X,R[f]]  \ar@{->>}[d]\\
      Y \ar@{<-}[rru]^{\zeta_f}\ar@{<-<}[rr] & & {\;X/R[f]}
                   }
$$
where $\eta_f$ is a regular epimorphism and $\zeta_f$ has a central kernel relation. If $f$ is a regular epimorphism then $\zeta_f$ is a central extension.

This factorisation has the following universal property. Any commutative diagram of undotted arrows as below (with $Z_f=X/[\nabla_X,R[f]]$), such that $\zeta'$ has a central kernel relation, produces a unique dotted map
$$ \xymatrix@=20pt{
  X \ar@{->>}[rd]^{\eta_f} \ar[dd]_{f} \ar[rr]^{x}  && X' \ar[rd]^{\eta'}\\
     & Z_f\ar[dl]^{\zeta_f} \ar@{.>}[rr]^<<<<<<<<{z}  && Z' \ar[dl]^{\zeta'}\\
  Y \ar[rr]_{y} && Y'
 }
$$ rendering the whole diagram commutative.
\end{prp}

\begin{prp}\label{univncentral}Let $n$ be a positive integer. In a finitely cocomplete regular Mal'tsev category, each morphism has a universal factorisation into a regular epimorphism followed by an $n$-fold centrally decomposable morphism. Each $n$-fold central extension has an initial factorisation into $n$ central extensions.\end{prp}
\proof
We proceed by induction, the case $n=1$ being treated in Proposition \ref{univcentral}. Suppose the assertion holds up to level $n-1$ with $f=\zeta^{n-1}_f\zeta^{n-2}_f\cdots\zeta^1_f\eta^1_f$ and $\eta^1_f$ a regular epimorphism. Take the universal factorisation of $\eta^1_f$. The universality of this new factorisation is then a straightforward consequence of the induction hypothesis and the universal property stated in Proposition \ref{univcentral}.

Starting with an $n$-fold central extension $f$, its universal factorisation through a composite of $n-1$ morphisms with central kernel relation makes the regular epimorphism $\eta^1_f$ a central extension by Lemma \ref{lefterase}, and therefore produces a factorisation into $n$ central extensions which is easily seen to be the initial one.
\endproof

\subsection{Regular pushouts}\label{regpushout0}In a regular category, any pullback square of regular epimorphisms is also a pushout square, as follows from the pullback stability of regular epimorphisms. In particular, a commuting square of regular epimorphisms
 $$ \xymatrix@=20pt{
         X \ar@{->>}[d]_{x}\ar@{->>}[r]^{f} & Y \ar@{->>}[d]^{y}\\
          X'   \ar@{->>}[r]_{f'} & Y'
                        }
     $$
is a pushout whenever the comparison map $(x,f):X\to X'\times_{Y'} Y$ to the pullback of $y$ along $f'$ is a regular epimorphism. Such pushouts will be called \emph{regular}, cf. \cite{Bourn3}. A regular pushout induces in particular regular epimorphisms on kernel relations which we shall denote $R(x,y):R[f]\onto R[f']$ and $R(f,f'):R[x]\onto R[y]$.

For a regular \emph{Mal'tsev} category the following more precise result holds.

\begin{prp}[cf. Proposition 3.3 in \cite{Bourn3}]\label{regpushout1}In a regular Mal'tsev category, a commuting square of regular epimorphisms like in (\ref{regpushout0}) is a regular pushout if and only if one of the following three equivalent conditions holds:
\begin{itemize}\item[(a)]the comparison map $X\to X'\times_{Y'}Y$ is a regular epimorphism;
\item[(b)]the induced map $R(x,y):R[f]\to R[f']$ is a regular epimorphism;
\item[(c)]the induced map $R(f,f'):R[x]\to R[y]$ is a regular epimorphism.\end{itemize}
Accordingly, central extensions are closed under regular pushouts.\end{prp}

\proof We already mentioned that (a) implies (b) and (c) in any regular category. It remains to be shown that in a regular Mal'tsev category (b) or (c) implies (a).

For this, it is useful to notice that in a regular category condition (a) holds if and only if the composite relation $R[f]\circ R[x]$ equals the kernel relation of the diagonal of the square by Theorem 5.2 of Carboni-Kelly-Pedicchio \cite{CKP}. Since this kernel relation is given by $x^{-1}(R[f'])$, and condition (b) just means that $x(R[f])=R[f']$, it suffices to establish the identity $R[f]\circ R[x]=x^{-1}(x(R[f]))$. In a regular Mal'tsev category, the composition of equivalence relations is symmetric and coincides with their \emph{join}. The join $R[f]\vee R[x]$ is easily identified with $x^{-1}(x([R[f]))$.

The second assertion follows from (b) resp. (c) and the closure of central kernel relations under direct image in regular Mal'tsev categories.\endproof

\begin{cor}\label{regpush}In a regular Mal'tsev category, any commutative square
 $$ \xymatrix@=20pt{
          X \ar@{->>}[d]_{x}\ar@{->>}[r]^{f} & Y \ar@{->>}[d]^{y} \ar@<1ex>[l]^{s}\\
          X'\ar@{->>}[r]^{f'} & Y' \ar@<1ex>[l]^{s'}
                        }
     $$
with a parallel pair of regular epimorphisms and a parallel pair of split epimorphisms is a regular pushout.\end{cor}

\proof The induced map $R(f,f'):R[x]\to R[y]$ is a split and hence regular epimorphism so that the pushout is regular by Proposition \ref{regpushout1}.\endproof

\begin{cor}\label{stpushout}In an exact Mal'tsev category, pushouts of regular epimorphisms along regular epimorphims exist and are regular pushouts.\end{cor}
\proof
Given a pair $(f,x)$ of regular epimorphisms with common domain, consider the following diagram
$$ \xymatrix@=20pt{
        R[f]  \ar@{->>}[d]_{\tilde x} \ar@<-1,ex>[r]_{p_0}\ar@<+1,ex>[r]^{p_1} & X \ar@{->>}[d]_{x}\ar[l]\ar@{->>}[r]^{f} & Y \ar@{.>>}[d]^{y}\\
        S \ar@<-1,ex>[r]_{q_0} \ar@<+1,ex>[r]^{q_1}  & X'  \ar[l] \ar@{.>>}[r]_{f'} & Y'
                       }
    $$in which $S$ denotes the direct image $x(R[f])$. By exactness, this equivalence relation on $X'$ has a quotient $Y'$. The induced right square is then a regular pushout.\endproof

\begin{rmk}It follows from \cite{Bourn0} that Corollary \ref{regpush} characterises regular Mal'tsev categories among regular categories, while \cite[Theorem 5.7]{CKP} shows that Corollary \ref{stpushout} characterises exact Mal'tsev categories among regular categories.\end{rmk}

\begin{rmk}\label{split}It is worthwhile noting that in \emph{any} category a commuting square of epimorphisms in which one parallel pair admits compatible sections is automatically a pushout square. Dually, a commuting square of monomorphisms in which one parallel pair admits compatible retractions is automatically a pullback square.\end{rmk}

\begin{lma}[cf. \cite{BG2}, Lemma 1.1]\label{kerpush} In a pointed regular category, a regular pushout induces a regular epimorphism on parallel kernels
     $$ \xymatrix@=20pt{
              {K[f]\;}  \ar@{->>}[d]_{K(x,y)} \ar@{>->}[r] & X \ar@{->>}[d]_{x} \ar@{->>}[r]^{f} & Y \ar@{->>}[d]^{y}\\
              {K[f']\;} \ar@{>->}[r]   & X'   \ar@{->>}[r]_{f'} & Y'
                             }
          $$so that the kernel $K[f']$ of $f'$ is the image under $x$ of the kernel $K[f]$ of $f$.\end{lma}

\proof This follows from the fact that, in any pointed category, the induced map on kernels factors as a pullback of the comparison map $(x,f):X\onto X'\times_{Y'}Y$ followed by an isomorphism.\endproof

\begin{prp}\label{nfold1}In a finitely cocomplete regular Mal'tsev category, consider the following diagram of pushouts
 $$ \xymatrix@=15pt{
      X  \ar@{->>}[rr]^{\eta_f} \ar@{->>}[d]_{x} && Z_f \ar@{->>}[d]^{z}\ar@{->>}[rr]^{\zeta_f} &&  Y\ar@{->>}[d]^{y}\\
   X' \ar@{->>}[rr]_{\eta_{f'}}  && Z_{f'}\ar@{->>}[rr]_{\zeta_{f'}} &&  Y'
                       }
    $$in which the upper row represents the universal factorisation of the regular epimorphism $f$ into a regular epimorphism $\eta_f$ followed by a central extension $\zeta_f$.

If the outer rectangle is a regular pushout then the right square as well, and the lower row represents the universal factorisation \ref{univcentral} of its composite $f'=\zeta_{f'}\eta_{f'}$.\end{prp}

\proof Since on vertical kernel relations $R[x]\to R[z]\to R[y]$ we get a regular epimorphism, the second one $R[z]\to R[y]$ is a regular epimorphism as well, and the right square is a regular pushout by Proposition \ref{regpushout1}. Therefore, $\zeta_{f'}$ is a central extension by Corollary \ref{regpush}. It remains to be shown that the lower row fulfills the universal property of factorisation \ref{univcentral}. Consider the following diagram
 $$ \xymatrix@=15pt{
       X  \ar@{->>}[rr]^{\eta_f} \ar@{->>}[d]_{x} && Z_f\ar@{->>}[d]_{z}\ar@{->>}[rr]^{\zeta_f} \ar@<1ex>@{.>}[dd]^<<<<{z''} &&  Y\ar@{->>}[d]^{y}\\
      X' \ar@{->>}[rr]^{\eta_{f'}} \ar[rrd]_{\eta} && Z_{f'}\ar@{->>}[rr]^{\zeta_{f'}} &&  Y'\\
     && Z' \ar@{->>}[rru]_{\zeta}
                        }
     $$with central extension $\zeta$. According to Proposition \ref{univcentral}, there is a unique dotted factorisation $z''$ making the diagram commute. Since the left square is a pushout, $z''$ factors uniquely and consistently through $z':Z_{f'}\to Z'$, showing that the lower row has indeed the required universal property.\endproof

\begin{prp}\label{trucc}
In a finitely cocomplete exact Mal'tsev category, the universal factorisation of a regular epimorphism through an $n$-fold central extension is preserved under pushouts along regular epimorphisms.\end{prp}
\proof
Let us consider the following diagram of pushouts
$$ \xymatrix@=15pt{
       X\ar@{->>}[d]_{x} \ar@{->>}[r]^{\eta} & Z_{n} \ar@{->>}[r]^{\zeta_{n}} \ar@{->>}[d]^{z_{n}} & Z_{n-1} \ar@{->>}[d]^{z_{n}} & ... &  Z_1  \ar@{->>}[d]^{z_1}  \ar@{->>}[r]^{\zeta_1} & Y \ar@{->>}[d]^{y} \\
       X' \ar@{->>}[r]_{\eta'}   & Z'_{n} \ar@{->>}[r]_{\zeta'_{n}} &  Z'_{n-1} & ... &  Z'_1 \ar@{->>}[r]_{\zeta'_1}& Y'
                       }
    $$in which the upper row is the universal factorisation \ref{univncentral} of a regular epimorphism $f:X\to Y$ through an $n$-fold central extension. By Corollary \ref{stpushout}, all pushouts are regular. Therefore, the morphisms $\zeta'_{k}:Z'_k\to Z'_{k-1}$ are central extensions for all $k$. It remains to be shown that the lower row satisfies the universal property of the factorisation \ref{univncentral} of $f':X'\to Y'$ through an $n$-fold central extension. This follows induction on $n$ beginning with the case $n=1$ proved in Proposition \ref{nfold1}.\endproof

\begin{prp}\label{nfoldslice}Let $\,\DD$ be an exact Mal'tsev category. Consider the following diagram of pushouts
    $$ \xymatrix@=15pt{
          X  \ar@{->>}[r]^{f_n} \ar@{->>}[d]_{x_n} & X_{n-1} \ar@{->>}[d]_{x_{n-1}} \ar@{->>}[r]^{f_{n-1}} & X_{n-2} \ar@{->>}[d]_{x_{n-2}} & ... & X_1 \ar@{->>}[r]^{f_1} \ar@{->>}[d]^{x_1} &  X_0 \ar@{->>}[d]^{x_0} \ar@{->>}[r]^{f_0} &  Y\ar@{->>}[d]^{y}
           \\X'  \ar@{->>}[r]_{f'_n} & X'_{n-1} \ar@{->>}[r]_{f'_{n-1}} & X'_{n-2} & ... & X'_1 \ar@{->>}[r]_{f'_1} &  X'_0 \ar@{->>}[r]_{f'_0} & Y'
           }
        $$in which $x_n:X\onto X'$ is a regular epimorphism.

If the upper row represents an $n$-fold central extension of the regular epimorphism $f_0:X_0\onto Y$ in the slice category $\DD/Y$ then the lower row represents an $n$-fold central extension of $f'_0:X'_0\onto Y'$ in the slice category $\DD/Y'$.\end{prp}
    \proof
Let us set $\phi_i=f_0\cdot f_1\cdots f_i$ and $\phi_i'=f_0'\cdot f'_1\cdots f'_i$. Since the indiscrete equivalence relation $\nabla_{f_0}$ on the object $f_0:X_0\onto Y$ of the slice category $\DD/Y$ is given by $R[f_0]$, our assumption on the upper row translates into the conditions $$[R[f_i],R[\phi_i]]=\Delta_{X_i}\text{ for }1\leq i \leq n.$$ Since any of the rectangles is a regular pushout by Corollary \ref{stpushout}, we get $x_i(R[f_i])=R[f'_i]$ and $x_i(R[\phi_i])=R[\phi'_i]$, and consequently $[R[f'_i],R[\phi'_i]]=\Delta_{X'_i}$ for all $i$.\endproof

\subsection{Regular pushouts in pointed Mal'tsev categories with binary sums}\label{pointedexactMalcev}In a pointed category with binary sums and binary products, each pair of objects $(X_1,X_2)$ defines a canonical comparison map $\theta_{X_1,X_2}:X_1+X_2\to X_1\times X_2$, uniquely determined by the requirement that the composite morphism$$\xymatrix@=15pt{X_i\,\ar@{>->}[r]&X_1+X_2\ar@{->}[rr]^{\theta_{X_1,X_2}}&&X_1\times X_2\ar@{->>}[r]&X_j}$$ is the identity (resp. the null morphism) if $i=j$ (resp. $i\not=j$), where $i,j\in\{1,2\}$.

Recall that $\theta_{X_1,X_2}$ is a strong epimorphism for all objects $X_1,X_2$ precisely when the category is \emph{unital} in the sense of the second author, and that every pointed Mal'tsev category is unital, cf. \cite{BB,Bourn0}. In a regular category strong and regular epimorphisms coincide.

Note also that an exact Mal'tsev category has coequalisers for reflexive pairs, so that an exact Mal'tsev category with binary sums has all finite colimits. In order to shorten terminology, we call \emph{$\sg$-pointed} any pointed category with binary sums.

\vspace{1ex}

Later we shall need the following two examples of regular pushouts.
\begin{prp}\label{regpushout2}
For any regular epimorphism $f:X\onto Y$ and any object $Z$ of a $\sg$-pointed regular Mal'tsev category, the following square
$$ \xymatrix@=20pt{  X+Z \ar@{->>}[r]^{\theta_{X,Z}} \ar@{->>}[d]_{f+Z} & X\times Z \ar@{->>}[d]^{f\times Z}\\
        Y+Z \ar@{->>}[r]_{\theta_{Y,Z}} & Y\times Z
                       }
    $$
is a regular pushout.\end{prp}

\proof

The regular epimorphism $\theta_{R[f],Z}:R[f]+Z\onto R[f]\times Z$ factors as below
$$ \xymatrix@=20pt{
        R[f]+Z \ar@{->>}[rr]^{\theta_{R[f],Z}} \ar[rd] &&  R[f] \times Z=R[f\times Z] \\
         & R[f+Z] \ar@{->>}[ur]
                       }
    $$
inducing a regular epimorphism $R[f+Z]\to R[f\times Z]$ on the vertical kernel relations of the square above. Proposition \ref{regpushout1} allows us to conclude.\endproof

\begin{cor}\label{corbr}
For any objects $X,Y,Z$ of a $\sg$-pointed regular Mal'tsev category, the following square
 $$ \xymatrix@=15pt{
       (X+Y)+Z  \ar@{->>}[rr]^{\theta_{X+Y,Z}} \ar@{->>}[d]_{\theta_{X,Y}+ Z} && (X+Y)\times Z\ar@{->>}[d]^{\theta_{X,Y}\times Z} \\
    (X\times Y)+Z  \ar@{->>}[rr]_{\theta_{X\times Y,Z}}  && (X\times Y) \times Z
                        }
     $$ is a regular pushout.
\end{cor}

\subsection{Central subobjects, centres and centralisers}\label{cooperate}\label{protomodular}\label{centralizer}In a pointed Mal'tsev category $(\DD,\star_\DD)$, two morphisms with common codomain $f:X\to Z$ and $g:Y\to Z$ are said to \emph{commute} \cite{Bourn5,Huq} if the square
$$ \xymatrix@=25pt{
      X\times Y \ar@{<-}[r]^<<<<<{(0_{YX},1_Y)} \ar@{<-}[d]_{(1_X,0_{XY})}\ar@{.>}[rd]^{\phi_{f,g}}&  Y\ar[d]^{g}\\
      X \ar@{->}[r]_{f} & Z
                   }
$$admits a (necessarily unique) filler $\phi_{f,g}:X\times Y\to Z$ making the whole diagram commute, where $0_{XY}:X\to\star_\DD\to Y$ denotes the zero morphism. A monomorphism $Z\into X$ which commutes with the identity $1_X:X\to X$ is called \emph{central}, and the corresponding subobject is called a \emph{central subobject} of $X$.

Every regular epimorphism $f:X\onto Y$ with central kernel relation $R[f]$ has a central kernel $K[f]$. In pointed protomodular categories, the converse is true: the centrality of $K[f]$ implies the centrality of $R[f]$, so that central extensions are precisely the regular epimorphisms with central kernel, cf. \cite[Proposition 2.2]{GL}.

Recall \cite{BB,Bourn} that a pointed category is \emph{protomodular} precisely when the category has pullbacks of split epimorphisms, and for each split epimorphism, section and kernel-inclusion form a strongly epimorphic cospan. Every finitely complete protomodular category is a Mal'tsev category \cite[Proposition 3.1.19]{BB}. The categories of groups and of Lie algebras are pointed protomodular. Moreover, in both categories, each object possesses a \emph{centre}, i.e. a maximal central subobject. Central group (resp. Lie algebra) extensions are thus precisely regular epimorphisms $f:X\onto Y$ with kernel $K[f]$ contained in the centre of $X$. This is of course the classical definition of a central extension in group (resp. Lie) theory.

In these categories, there exists more generally, for each subobject $N$ of $X$, a so-called \emph{centraliser}, i.e. a subobject  $Z(N\!\into\! X)$ of $X$ which is maximal among subobjects commuting with $N\!\into\!X$. The existence of centralisers has far-reaching consequences, as shown by James Gray and the second author, cf. \cite{BGr,Gr,Gr1}. Since they are useful for our study of nilpotency, we discuss some of them here.

Following \cite{Bourn0}, we denote by $\Pt_Z(\DD)$ the category of split epimorphisms (with chosen section) in $\DD$ over a fixed codomain $Z$, cf. Section \ref{pointfibration}. For each $f:Z\to Z'$ pulling back along $f$ defines a functor $f^*:\Pt_{Z'}(\DD)\to\Pt_Z(\DD)$ which we call \emph{pointed base-change} along $f$. In particular, the terminal map $\omega_Z:Z\to 1_\DD$ defines a functor $(\omega_Z)^*:\DD\to\Pt_Z(\DD)$. Since in a pointed regular Mal'tsev category $\DD$, morphisms commute if and only if their images commute, morphisms in $\Pt_Z(\DD)$ of the form
$$ \xymatrix@=15pt{
        X\times Z \ar@{->}[rr]^{\phi_{f,f'}}\ar[dr]_{p_Z} &&  Y \ar[dl]_s \\
        & Z\ar@<-1ex>[ul]\ar@<+1ex>[ur]_r
                       }
    $$correspond bijectively to morphisms $f:X\to K[r]$ such that $X\overset{f}{\longrightarrow}K[r]\into Y$ commutes with $s:Z\into Y$ in $\DD$. Therefore, if split subobjects have centralisers in $\DD$, then for each object $Z$, the functor $(\omega_Z)^*:\DD\to\Pt_Z(\DD):X\mapsto X\times Z$ admits a right adjoint $(\omega_Z)_*:\Pt_Z(\DD)\to\DD:(r,s)\mapsto K[r]\cap Z(s)$.

A category with the property that for each object $Z$, the functor $(\omega_Z)^*$ has a right adjoint is called \emph{algebraically cartesian closed} \cite{BGr}. Algebraic cartesian closedness implies canonical isomorphisms $(X\times Z)+_Z(Y\times Z)\cong (X+Y)\times Z$ for all objects $X,Y,Z$, a property we shall call \emph{algebraic distributivity}, cf. Section \ref{extensive}.

\subsection{An algebraic Beck-Chevalley condition}\label{Beck-Chevalley}--\vspace{1ex}

The dual of an elementary topos is an exact Mal'tsev category, cf. \cite[Remark 5.8]{CKP}. This suggests that certain diagram lemmas for elementary toposes admit a dual version in our algebraic setting. Supporting this analogy we establish here an ``algebraic dual'' of the well-known \emph{Beck-Chevalley condition}. As a corollary we get a diagram lemma which will be used several times in Section \ref{degreen}.

Another instance of the same phenomenon is the \emph{cogluing lemma} for regular epimorphisms in exact Mal'tsev categories (cf. proof of Theorem \ref{folklore}a and Corollary \ref{stpushout}) which is dual to a gluing lemma for monomorphisms in an elementary topos.

\begin{lma}[cf. Lemma 1.1 in \cite{Gran1}]\label{propbr}Consider a commutative diagram
$$ \xymatrix@=25pt{ X \ar@{->>}[d]\ar@{->>}[r]^x&
        X' \ar@{->>}[d]\ar@{->}[r] & X'' \ar@{->}[d]\\
    Y\ar@{->>}[r]_y &     Y' \ar@{->}[r]\ & Y''
                       }
    $$in a regular category.

If the outer rectangle is a pullback and the left square is a regular pushout (\ref{regpushout0}) then left and right squares are pullbacks.\end{lma}
\proof The whole diagram contains three comparison maps: one for the outer rectangle, denoted $\phi:X\to Y\times_{Y''}X''$, one for the left and one for the right square, denoted respectively $\phi_{l}:X\to Y\times_{Y'}X'$ and $\phi_{r}:X'\to Y'\times_{Y''}X''$. We get the identity $\phi=y^*(\phi_{r})\circ\phi_{l}$ where $y^*$ denotes base-change along $y$. Since the outer rectangle is a pullback, $\phi$ is invertible so that $\phi_l$ is a section and $y^*(\phi_r)$ a retraction.

Since the left square is a regular pushout, the comparison map $\phi_l$ is a regular epimorphism and hence $\phi_l$ and $y^*(\phi_r)$ are both invertible. Since $y$ is a regular epimorphism in a regular category, base-change $y^*$ is conservative so that $\phi_r$ is invertible as well, i.e. both squares are pullbacks.\endproof

\begin{prp}[Algebraic Beck-Chevalley]\label{B-C1}
Let $\,\DD$ be an exact Mal'tsev category with pushouts of split monomorphisms along regular epimorphisms.

Any pushout of regular epimorphisms
$$ \xymatrix@=20pt{
        \bar U \ar@{->>}[d]_{u} \ar@{->>}[r]^{\bar g}   & \bar V \ar@{->>}[d]^{v}  \\
        U \ar@{->>}[r]_{g} & V
                       }
    $$yields a functor isomorphism $\bar g_!u^*\cong v^*g_!$ from the fibre $\Pt_U(\DD)$ to the fibre $\Pt_V(\DD)$.
\end{prp}
\proof
We have to show that for any point $(r,s)$ over $U$, the following diagram
$$ \xymatrix@=12pt{
  &    \bar U' \ar@{->>}[rrr]^{\bar g'}  \ar@{->>}[dr]^{u'} \ar[ddl]_<<<<{\bar r} &&& \bar V' \ar[ddl]_<<<<{\bar r'}  \ar@{->>}[dr]^<<<{v'}\\
  &&   U' \ar[ddl]_<<<<{r} \ar@{->>}[rrr]^<<<<<<{g'} &&&  V' \ar[ddl]_{r'} \\
   \bar U \ar@{->>}[dr]_<<<{u} \ar@{->>}[rrr]_>>>>>>{\bar g}\ar@<-1ex>[ruu] &&&  \bar V \ar@{->>}[dr]_<<<{v}\ar@<-1ex>[ruu] \\
  &   U \ar@{->>}[rrr]_{g} \ar@<-1ex>[ruu] &&& V \ar@<-1ex>[ruu]
                   }
$$
in which $(\bar r,\bar s)=u^*(r,s)$ and $g_!(r,s)=(r',s')$ and  $\bar g_!(\bar r,\bar s)=(\bar r',\bar s')$, has a right face which is a downward-oriented pullback; indeed, this amounts to the required identity $v^*(r',s')=(\bar{r}',\bar{s}')$.

Since bottom face and the upward-oriented front and back faces are pushouts, the top face is a pushout as well, which is regular by Corollary \ref{stpushout}. Taking pullbacks in top and bottom faces induces a split epimorphism $U'\times_{V'}\bar{V}'\onto U\times_V\bar{V}$ through which the left face of the cube factors as in the following commutative diagram
$$ \xymatrix@=25pt{ \bar{U'} \ar@{->>}[d]_{\bar{r}}\ar@{->>}[r]&
        U'\times_{V'}\bar{V}' \ar@{->>}[d]\ar@{->>}[r] & U' \ar@{->>}[d]^{r}\\
    \bar{U}\ar@{->>}[r] &     U\times_V\bar{V} \ar@{->>}[r]\ & U
                       }
    $$in which the left square is a regular pushout by Corollary \ref{regpush}. Lemma \ref{propbr} shows then that the right square is a pullback. Therefore, we get the following cube

$$ \xymatrix@=12pt{
  &    U'\times_{V'}\bar{V}' \ar@{->>}[rrr] \ar@{->>}[dr] \ar[ddl]_<<<<{} &&& \bar V' \ar[ddl]_<<<<{\bar r'}  \ar@{->>}[dr]^<<<{v'}\\
  &&   U'\ar[ddl]_<<<<{r} \ar@{->>}[rrr]^<<<<<<{} &&&  V' \ar[ddl]_{r'} \\
   U\times_V\bar{V} \ar@{->>}[dr]_<<<{} \ar@{->}[rrr]_>>>>>>{\bar g}\ar@<-1ex>[ruu] &&&  \bar V \ar@{->>}[dr]_<<<{v}\ar@<-1ex>[ruu] \\
  &   U \ar@{->>}[rrr]_{g} \ar@<-1ex>[ruu] &&& V \ar@<-1ex>[ruu]
                   }
$$in which the downward-oriented left face and the bottom face are pullbacks. Therefore, the composite of the top face followed by the downward-oriented right face is a pullback. Moreover, as above, the top face is a regular pushout. It follows then from Lemma \ref{propbr} that the downward-oriented right face is a pullback as required.\endproof

\begin{cor}\label{B-C3}
In an exact Mal'tsev category with pushouts of split monomorphisms along regular epimorphisms, each commuting square of natural transformations of split epimorphisms
$$ \xymatrix@=12pt{
  &    \bar U' \ar@{->>}[rrr]^{\bar g'}  \ar@{->>}[dr]^{u'} \ar[ddl]_<<<<{\bar f} &&& \bar V' \ar[ddl]_<<<<{\bar f'}  \ar@{->>}[dr]^<<<{v'}\\
  &&   U' \ar[ddl]_<<<<{f} \ar@{->>}[rrr]^<<<<<<{g'}  &&&  V' \ar[ddl]_{f'}  \\
   \bar U \ar@<-1ex>@{->>}[dr]_<<<{u} \ar@{->>}[rrr]_>>>>>>{\bar g}\ar@<-1ex>[ruu] &&&  \bar V \ar@<-1ex>@{->>}[dr]^<<<{v}\ar@<-1ex>[ruu] \\
  &   U \ar@{->>}[rrr]_{g} \ar@<-1ex>[ruu]  &&& V \ar@<-1ex>[ruu]
                   }
$$
such that all horizontal arrows are regular epimorphisms, and front and back faces are upward-oriented pushouts, induces the following upper pushout square
$$ \xymatrix@=10pt{
  \bar U' \ar@{->>}[rrr]^{\bar g'}  \ar@{->>}[ddrr]^{(\bar{f},u')} \ar@<-1ex>[ddddr]_{\bar f} &&& \bar V' \ar@<-1ex>[ddddr]_<<<<{\bar f'}  \ar@{->>}[ddrr]^{(\bar{f}',v')}\\
  &&&&&\\
  &&  \bar{U}\times_UU' \ar@<1ex>[ddl]_<<<<{} \ar@{->>}[rrr]^<<<<<<{\bar{g}\times g'}  &&& \bar{V}\times_VV' \ar@<1ex>[ddl]_<<<<{}  \\
  &&&&&\\
  &  \bar U \ar@{->>}[rrr]_{\bar g} \ar@<-2ex>[ruu] \ar[uuuul] &&& \bar V \ar@<-2ex>[ruu] \ar[uuuul]
                   }
$$
in which the kernel relation of the regular epimorphism $(\bar{f}',v'):\bar{V}'\to \bar{V}\times_VV'$ may be identified with the intersection $R[\bar f']\cap R[v']$.
\end{cor}
\proof
Taking downward-oriented pullbacks in left and right face of the first diagram yields precisely a diagram as studied in the proof of Proposition \ref{B-C1}. This implies that the front face of the second diagram is an upward-oriented pushout. Since the back face of the second diagram is also an upward-oriented pushout, the upper square is a pushout as well, as asserted. The kernel relation of the comparison map $(\bar{f}',v'):\bar{V}'\to \bar{V}\times_VV'$ is the intersection of the kernel relations of $\bar{f}'$ and of $v'$.\endproof

\section{Affine objects and nilpotency}\label{DF}

\subsection{Affine objects}

\begin{dfn}\label{affineobject}Let $\CC$ be a full subcategory of a Mal'tsev category $\DD$.

An object $X$ of $\,\DD$ is said to be \emph{$\CC$-affine} if there exists a morphism $f:X\to Y$ in $\DD$ with central kernel relation $R[f]$ and with codomain $Y$ in $\CC$.

The morphism $f$ is called a \emph{$\CC$-nilindex} for $X$.
\end{dfn}

We shall write $\Aff_\CC(\DD)$ for the full replete subcategory of $\DD$ spanned by the $\CC$-affine objects of $\DD$. Clearly $\Aff_\CC(\DD)$ contains $\CC$.

When $\CC$ consists only of a terminal object $1_\DD$ of $\DD$, we call the $\CC$-affine objects simply the \emph{affine objects} of $\DD$ and write $\Aff_\CC(\DD)=\Aff(\DD)$. Recall that the unique morphism $X\to 1_\DD$ has a central kernel relation precisely when the indiscrete equivalence relation on $X$ centralises itself, which amounts to the existence of a (necessarily unique associative and commutative) \emph{Mal'tsev operation} on $X$.

When $\DD$ is \emph{pointed}, such a Mal'tsev operation on $X$ induces (and is induced by) an \emph{abelian group} structure on $X$. For a pointed Mal'tsev category $\DD$, the category $\Aff(\DD)$ of affine objects is thus the category $\Ab(\DD)$ of abelian group objects of $\DD$.

\begin{rmk}\label{centralextension}When $\DD$ is a \emph{regular} Mal'tsev category, any nilindex $f:X\to Y$ factors as a regular epimorphism $\tilde{f}:X\onto f(X)$ followed by a monomorphism $f(X)\into Y$ with codomain in $\CC$; therefore, if $\CC$ is \emph{closed under taking subobjects in $\DD$}, this defines a strongly epimorphic nilindex $\tilde{f}$ for $X$ with same central kernel relation as $f$. In other words, for regular Mal'tsev categories $\DD$ and subcategories $\CC$ which are closed under taking subobjects in $\DD$, the $\CC$-affine objects of $\DD$ are precisely the objects which are obtained as central extensions of objects of $\CC$.\end{rmk}

\begin{prp}\label{epireflective}For any full subcategory $\CC$ of a Mal'tsev category $\DD$, the subcategory $\Aff_\CC(\DD)$ is closed under taking subobjects in $\DD$. If $\,\CC$ is closed under taking binary products in $\DD$ then $\Aff_\CC(\DD)$ as well, so that $\Aff_\CC(\DD)$ is finitely complete.\end{prp}
\proof
Let $m:X\into X'$ be a monomorphism with $\CC$-affine codomain $X'$. If $f':X'\to Y'$ is a nilindex for $X'$, then $fm:X\to Y'$ is a nilindex for $X$, since central equivalence relations are stable under pointed base-change along monomorphisms and we have $R[fm]=m^{-1}(R[f])$.

If $X$ and $Y$ are  $\CC$-affine with nilindices $f$ and $g$ then $f\times g$ is a nilindex for $X\times Y$ since maps with central kernel relations are stable under products.
\endproof

A \emph{Birkhoff subcategory} \cite{JK} of a regular category $\DD$ is a subcategory $\CC$ which is closed under taking subobjects, products and quotients in $\DD$. A Birkhoff subcategory of an exact (resp. Mal'tsev) category is exact (resp. Mal'tsev), and regular epimorphisms in $\CC$ are those morphisms in $\CC$ which are regular epimorphisms in $\DD$.

If $\DD$ is a variety (in the single-sorted monadic sense) then Birkhoff subcategories of $\DD$ are precisely subvarieties of $\DD$, cf. the proof of Lemma \ref{variety} below.

\begin{prp}\label{inheritance}Let $\CC$ be a full subcategory of an exact Mal'tsev category $\DD$.

If $\,\CC$ is closed under taking subobjects and quotients in $\DD$ then $\Aff_\CC(\DD)$ as well. In particular, if $\,\CC$ is a Birkhoff subcategory of $\DD$, then $\Aff_\CC(\DD)$ as well.
\end{prp}
\proof
Let $X$ be a $\CC$-affine object of $\DD$ with nilindex $f:X\to Y$. We can suppose $f$ is a regular epimorphism, cf. Remark \ref{centralextension}. Thanks to Proposition \ref{epireflective} it remains to establish closure under quotients. Let $g:X\onto X'$ be a regular epimorphism in $\DD$. Since $\DD$ is exact, the pushout of $f$ along $g$ exists in $\DD$
$$ \xymatrix@=20pt{
         X \ar@{->>}[d]_{g}\ar@{->>}[r]^{f} & Y \ar@{->>}[d]^{h}\\
        X'   \ar@{->>}[r]_{f'} & Y'
                       }
    $$and $f'$ is a central extension since $f$ is, cf. Corollary \ref{stpushout}. By hypothesis $\CC$ is stable under quotients. Therefore the quotient $Y'$ belongs to $\CC$, and $f'$ is a nilindex for $X'$ so that $X'$ is $\CC$-affine as required.\endproof

\subsection{The $\CC$-lower central sequence}--\vspace{1ex}

Definition \ref{affineobject} is clearly the beginning of an iterative process. We write $\CC=\Nil^0_\CC(\DD)$ and define inductively $\Nil^n_\CC(\DD)$ to be the category $\Aff_{\Nil_\CC^{n-1}(\DD)}(\DD)$. The objects of this category $\Nil^n_\CC(\DD)$ are called the \emph{$\CC$-nilpotent objects of order $n$ of $\,\DD$}, and we get the following diagram
$$ 
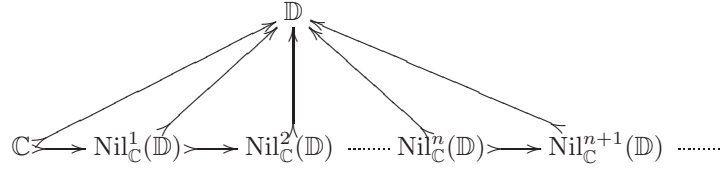
\xymatrix@=15pt{
&  & \mathbb D && \\
&&&&\\
{\mathbb C\;}\ar@{>->}[r]\ar@{>->}[rruu]  & {\Nil^1_\CC(\DD)\;}\ar@{>->}[r]\ar@{>->}[ruu]  & {\Nil^2_\CC(\DD)\;}\ar@{.}[r] \ar@{>->}[uu] & {\Nil_\CC^n(\DD)\;}\ar@{>->}[r]   \ar@{>->}[luu]& {\Nil_\CC^{n+1}(\DD)\;} \ar@{>->}[lluu] \ar@{.}[r] &
                  }
  $$
which we call the \emph{$\CC$-lower central sequence} of $\DD$.\vspace{1ex}

If $\CC=\{1_\DD\}$, we obtain the (absolute) lower central sequence of $\DD$:
  $$ \xymatrix@=15pt{
  &  & \DD && \\
  &&&&\\
  {\{1_\DD\}\;}\ar@{>->}[r] \ar@{>->}[rruu] & {\Nil^1(\DD)\;}\ar@{>->}[r]\ar@{>->}[ruu]  & {\Nil^2(\DD)\;}\ar@{.}[r] \ar@{>->}[uu] & {\Nil^n(\DD)\;}\ar@{>->}[r] \ar@{>->}[luu]  & {\Nil^{n+1}(\DD)\;} \ar@{>->}[lluu] \ar@{>->}[lluu] \ar@{.}[r] &
                    }
    $$

\begin{rmk}It follows from Remark \ref{centralextension} and an iterative application of Proposition \ref{inheritance} that for an exact Mal'tsev category $\DD$, the nilpotent objects of order $n$ are precisely those which can be obtained as an \emph{$n$-fold central extension} of the terminal object $1_\DD$ and that moreover $\Nil^n(\DD)$ is a Birkhoff subcategory of $\DD$.

If $\DD$ is the category of groups (resp. Lie algebras) then $\Nil^n(\DD)$ is precisely the full subcategory spanned by nilpotent groups (resp. Lie algebras) of class $\leq n$. Indeed, it is well-known that a group (resp. Lie algebra) is nilpotent of class $\leq n$ precisely when it can be obtained as an $n$-fold ``central extension'' of the trivial group (resp. Lie algebra), and we have seen in Section \ref{protomodular} that the group (resp. Lie) theorist's definition of central extension agrees with ours. We will see in Proposition \ref{iteratedSmith} below that the equivalence between the central extension definition and the iterated commutator definition of nilpotency carries over to our general context of finitely cocomplete exact Mal'tsev categories. Huq \cite{Huq} had foreseen a long time ago that a categorical approach to nilpotency was possible. Everaert-Van der Linden \cite{EV} recast Huq's approach in modern language in the context of semi-abelian categories.\end{rmk}

\begin{dfn}
A Mal'tsev category $\DD$ with full subcategory $\CC$ is called \emph{$\CC$-nilpotent of order $n$} (resp. \emph{of class $n$}) if $\,\DD=\Nil^n_\CC(\DD)$ (resp. if $n$ is the least such integer).
\end{dfn}

When $\CC=\{1_\DD\}$ the prefix $\CC$ will be dropped, and instead of ``nilpotent of order $n$'' we also just say ``$n$-nilpotent''.

\begin{prp}\label{nilpotentcategory}A Mal'tsev category is $n$-nilpotent if and only if each morphism is $n$-fold centrally decomposable.

A regular Mal'tsev category is $n$-nilpotent if and only if each morphism factors as an $n$-fold central extension followed by a monomorphism.\end{prp}
\proof The second statement follows from the first by Lemma \ref{relativecentral}.
If each morphism is $n$-fold centrally decomposable, then this holds for terminal maps $\omega_X:X\to 1_{\DD}$, so that all objects are $n$-nilpotent. Conversely, assume that all objects are $n$-nilpotent, i.e. that for all objects $X$, the terminal map $\omega_X$ is $n$-fold centrally decomposable. Then, for each morphism $f:X\to Y$, the identity $\omega_X=\omega_Yf$ together with Lemma \ref{lefterase} imply that $f$ is $n$-fold centrally decomposable as well.\endproof

\subsection{Epireflections, Birkhoff reflections and central reflections}--\vspace{1ex}

We shall see that if $\CC$ is a reflective subcategory of $\DD$, then the categories $\Nil_\CC^n(\DD)$ are again reflective subcategories of $\DD$, provided $\DD$ and the reflection fulfill suitable conditions. In order to give precise statements we need to fix some terminology.

A full replete subcategory $\CC$ of $\DD$ is called \emph{reflective} if the inclusion $\CC\inc\DD$ admits a left adjoint functor $I:\DD\to\CC$, called \emph{reflection}. The unit of the adjunction at an object $X$ of $\DD$ will be denoted by $\eta_X:X\to I(X)$. Reflective subcategories $\CC$ are stable under formation of limits in $\DD$. In particular, reflective subcategories of Mal'tsev categories are Mal'tsev categories.

A reflective subcategory $\CC$ of $\DD$ is called \emph{strongly epireflective} and the reflection $I$ is called a \emph{strong epireflection} if the unit $\eta_X:X\to I(X)$ is pointwise a \emph{strong epimorphism}. Strongly epireflective subcategories are characterised by the property that $\CC$ is closed under taking subobjects in $\DD$. In particular, strongly epireflective subcategories of regular categories are regular categories. %Indeed, the inclusion of an epireflective subcategory preserves and reflects relations and difunctional relations (cf. \cite{Ga}), and a Mal'tsev category can also be defined as a category in which every relation is difunctional, cf. \cite[Proposition 2.2.8]{BB}.

A \emph{Birkhoff reflection} (cf. \cite{BR}) is a strong epireflection $I:\DD\to\CC$ such that for each regular epimorphism $f:X\to Y$ in $\DD$, the following naturality square
$$ \xymatrix@=20pt{
        X \ar@{->>}[d]_{f}  \ar@{->>}[r]^{\eta_X} & I(X) \ar@<-1,ex>@{->>}[d]^{I(f)}\\
        Y  \ar@{->>}[r]_{\eta_Y} & I(Y)
                       }
    $$is a \emph{regular pushout} (see Section \ref{regpushout0} and Proposition \ref{regpushout1}).

A subcategory of $\DD$ defined by a Birkhoff reflection is a Birkhoff subcategory of $\DD$, and is thus exact whenever $\DD$ is. It follows from Corollary \ref{stpushout} that a reflective subcategory of an exact Mal'tsev category is a Birkhoff subcategory \emph{if and only if} the reflection is a Birkhoff reflection.

A \emph{central reflection} is a strong epireflection $I:\DD\to\CC$ with the property that the unit $\eta_X:X\onto I(X)$ is pointwise a central extension.\vspace{2ex}

The following exactness result will be used at several places. In the stated generality, it is due to Diana Rodelo and the second author \cite{BR}, but the interested reader can as well consult \cite{JK,Gran2,GR} for closely related statements.

\begin{prp}\label{Diana}In a regular Mal'tsev category, strong epireflections preserve pullback squares of split epimorphisms, and Birkhoff reflections preserve pullbacks of split epimorphisms along regular epimorphisms.\end{prp}
\proof See Proposition 3.4 and Theorem 3.16 in \cite{BR}.\endproof

\begin{lma}\label{lemma1.1}
Let $\CC$ be a reflective subcategory of $\DD$ with reflection $I$, and assume that $\eta_X:X\to I(X)$ factors through an epimorphism $f:X\onto Y$ as in:
$$ \xymatrix@=15pt{
      X \ar@{->}[rr]^{\eta_X} \ar@{->>}[d]_{f} && I(X)\\
      Y  \ar@{->}[rru]_{\eta}            }
  $$
Then $I(f)$ is an isomorphism and we have $\eta=I(f)^{-1}\eta_Y$.
\end{lma}
\proof
Consider the following diagram
$$ \xymatrix@=15pt{
      X \ar@{->}[rr]^{\eta_X} \ar@{->>}[d]_{f} && I(X) \ar[d]^{I(f)} \\
      Y  \ar@{->}[rr]_{\eta_Y} \ar@{->}[rru]_{\eta} && I(Y)         }
  $$
  where the lower triangle commutes because $f$ is an epimorphism. If we apply the reflection $I$ to the whole diagram we get two horizontal isomorphisms $I(\eta_X)$ and $I(\eta_Y)$. It follows that $I(\eta)$ is an isomorphism as well, hence so is $I(f)$, and $\eta=I(f)^{-1}\eta_Y$.
\endproof

\begin{lma}\label{lemma1.2}
For any reflective subcategory $\CC$ of a Mal'tsev category $\DD$ the $\CC$-affine objects of $\DD$ are those $X$ for which the unit $\eta_X$ has a central kernel relation.

\end{lma}
\proof
If $\eta_X:X\rightarrow I(X)$ has a central kernel relation then $X$ is $\CC$-affine. Conversely, let $X$ be $\CC$-affine with nilindex $f:X\rightarrow Y$. Then $Y$ is an object of the reflective subcategory $\CC$ so that $f$ factors through $\eta_X:X\to I(X)$.  Accordingly, we get $R[\eta_X] \subset R[f]$, and hence $R[\eta_X]$ is central because $R[f]$ is.
\endproof

\begin{cor}\label{cor1.3}A reflection $I:\DD\to\CC$ of a regular Mal'tsev category $\,\DD$ is central if and only if $\,\DD$ is $\CC$-nilpotent of order $1$ (i.e. all objects of $\,\DD$ are $\CC$-affine).\end{cor}

\begin{thm}[cf. \cite{Bourn4}, Sect. 4.3 in \cite{Huq}, Prp. 7.8 in \cite{EV} and Thm. 3.6 in \cite{JP}]\label{noBirkhoff}--

For a reflective subcategory $\,\CC$ of a finitely cocomplete regular Mal'tsev category $\,\DD$, the category $\Aff_\CC(\DD)$ is a strongly epireflective subcategory of $\,\DD$.

The associated strong epireflection $I_\CC^1:\DD\to\Aff_\CC(\DD)$ is obtained by factoring the unit $\eta_X:X\to I(X)$ universally through a map with central kernel relation.

If $\,\CC$ is a reflective Birkhoff subcategory of a finitely cocomplete exact Mal'tsev category $\DD$, then the reflection $I_\CC^1:\DD\to\Aff_\CC(\DD)$ is a Birkhoff reflection.
\end{thm}
\proof
Proposition \ref{univcentral} yields the following factorisation of the unit:
$$ \xymatrix@=15pt{
      X \ar[rr]^{\eta_X} \ar@{->>}[d]_{\eta_X^1} && I(X)\\
      I_\CC^1(X)  \ar[rru]_{\bar{\eta}_X}            }
  $$
Since $\bar{\eta}_X$ has a central kernel relation and $I(X)$ is an object of $\mathbb C$, the object $I_\CC^1(X)$ belongs to $\Aff_\CC(\DD)$. We claim that the maps $\eta^1_X:X\onto I_\CC^1(X)$ have the universal property of the unit of an epireflection $I^1_\CC:\DD\to\Aff_\CC(\DD)$.

Let $f:X\rightarrow T$ be a map with $\CC$-affine codomain $T$ which means that $\eta_T$ has a central kernel relation. Then consider the following diagram:
$$ \xymatrix@=20pt{
 X \ar@{->>}[rd]^{\eta_X^1} \ar[dd]_{\eta_X} \ar[rr]^{f}  && T \ar[dd]^{} \ar@{=}[rd]\\
  & I_\CC^1(X)  \ar[dl]^{\bar{\eta}_X} \ar@{.>}[rr]^<<<<<<{\bar f}  && T \ar[dl]^{\eta_T}\\
 I(X) \ar[rr]_{I(f)} && I(T)
 }
$$
According to Proposition \ref{univcentral}, there is a unique factorisation $\bar f$ making the diagram commute. If $\DD$ is exact and $\CC$ a Birkhoff subcategory, the subcategory $\Aff_\CC(\DD)$ is closed under taking subobjects and quotients by Proposition \ref{inheritance}. The reflection $I^1_\CC$ is thus a Birkhoff reflection in this case.\endproof

\begin{rmk}\label{Huq}A reflective Birkhoff subcategory $\CC$ of a semi-abelian category $\DD$ satisfies all hypotheses of the preceding theorem. In this special case, the Birkhoff reflection $I^1_\CC:\DD\to\Aff_\CC(\DD)$ is given by the formula $$I^1_\CC(X)=X/[X,K[\eta_X]]$$where $[X,K[\eta_X]]$ is the Huq commutator of $X$ and $K[\eta_X]$, cf. Section \ref{cooperate}.

Indeed, the pointed protomodularity of $\DD$ implies (cf. \cite[Proposition 2.2]{GL}) that the kernel of the quotient map $X\to X/[\nabla_X,R[\eta_X]]$ is canonically isomorphic to the Huq commutator $[X,K[\eta_X]]$ so that the formula follows from Proposition \ref{univcentral}.\end{rmk}

\subsection{The Birkhoff nilpotency tower}\label{Birkhoffcentral}--\vspace{1ex}

According to Theorem \ref{noBirkhoff}, any reflective Birkhoff subcategory $\CC$ of a finitely cocomplete exact Mal'tsev category $\DD$ produces iteratively the following commutative diagram of Birkhoff reflections:

$$ \xymatrix@=15pt{
&  & \DD \ar[lldd]_>>>>>>{I} \ar[ldd]_>>>>{I_\CC^1} \ar[dd]_>>>>{I_\CC^2}  \ar[ddr]^>>>>{I_\CC^{n}} \ar[ddrr]^>>>>{I_\CC^{n+1}}&& \\
&&&&\\
{\CC\;}\ar@{<-}[r]  & {\Nil^1_\CC(\DD)\;}\ar@{<-}[r]  & {\Nil^2_\CC(\DD)\;}\ar@{.}[r]  & {\Nil^n_\CC(\DD)\;}\ar@{<-}[r]  & {\Nil^{n+1}_\CC(\DD)\;} \ar@{.}[r]  &
                  }
  $$

A Birkhoff subcategory of an exact Mal'tsev category is an exact Mal'tsev category so that the subcategories $\Nil^n_\CC(\DD)$ are all exact Mal'tsev categories, and the horizontal reflections $\Nil^{n+1}_\CC(\DD)\to\Nil^n_\CC(\DD)$ are \emph{central} reflections by Corollary \ref{cor1.3}.

In the special case $\CC=\{1_\DD\}$ we get the following commutative diagram of Birkhoff subcategories and Birkhoff reflections:

$$ \xymatrix@=15pt{
&  & \DD \ar[lldd]_>>>>>{I} \ar[ldd]_>>>>{I^1} \ar[dd]_>>>>{I^2}  \ar[ddr]^>>>>{I^n} \ar[ddrr]^>>>>{I^{n+1}}&& \\
&&&&\\
{\{1_\DD\}\;}\ar@{<-}[r]  & {\Nil^1(\DD)\;}\ar@{<-}[r]  & {\Nil^2(\DD)\;}\ar@{.}[r]  & {\Nil^n(\DD)\;}\ar@{<-}[r]  & {\Nil^{n+1}(\DD)\;} \ar@{.}[r]  &
                  }
  $$
If $\DD$ is pointed, then the first Birkhoff reflection $I^1=I^1_{\{\star_\DD\}}:\DD\to\Nil^1(\DD)$ can be identified with the classical abelianisation functor $\DD\to\Ab(\DD)$. In particular, the abelian group objects of $\DD$ are precisely the nilpotent objects of order $1$.

When $\CC$ is a reflective Birkhoff subcategory of a finitely cocomplete exact Mal'tsev category $\DD$, then $\DD$ is $\CC$-nilpotent of class $n$ if and only if $n$ is the least integer such that either the unit of the $n$-th Birkhoff reflection $I_\CC^n$ is invertible, or equivalently, the $(n-1)st$ Birkhoff reflection $I_\CC^{n-1}$  is a \emph{central} reflection, see Corollary \ref{cor1.3}.

\begin{prp}\label{iteratedSmith}For an exact Mal'tsev category $\DD$ with binary sums, the unit of the $n$-th Birkhoff reflection  $\eta_X^n:X\to I^n(X)$ is given by quotienting out the iterated Smith commutator $[\nabla_X,[\nabla_X,[\nabla_X,\cdots,\nabla_X]]]$ of length $n+1$.

If $\,\DD$ is semi-abelian, this unit is also given as the quotient of $X$ by the iterated Huq commutator $[X,[X,[X,\dots,X]]]$ of length $n+1$.\end{prp}

\proof The second statement follows from the first by Remark \ref{Huq}. The first statement follows from the inductive construction of $\eta_X^n$ in the proof of Theorem \ref{noBirkhoff} together with Proposition \ref{univcentral}.\endproof

A finite limit and finite colimit preserving functor is called \emph{exact}. A functor between exact Mal'tsev categories with binary sums is exact if and only if it preserves finite limits, regular epimorphisms and binary sums, cf. Section \ref{pointedexactMalcev}.

\begin{lma}\label{exact}Any exact functor $F:\DD\to\EE$ between exact Mal'tsev categories with binary sums commutes with the $n$-th Birkhoff reflections, i.e. $I^n_\EE\circ F\cong F_{|\Nil^n(\DD)}\circ I^n_\DD$.\end{lma}
\proof According to Theorem \ref{noBirkhoff} and Proposition \ref{univcentral} it suffices to show that $F$ takes the canonical factorisation of $f:X\onto Y$ through the central extension $\zeta_f:X/[\nabla_X,R[f]]\onto Y$ to the factorisation of $F(f):F(X)\onto F(Y)$ through the central extension $\zeta_{F(f)}:F(X)/[\nabla_{F(X)},R[F(f)]]\onto F(Y)$. Since $F$ is left exact, we have $F(\nabla_X)=\nabla_{F(X)}$ and $F(R[f])=R[F(f)]$, and since $F$ preserves regular epimorphisms, we have $F(X/[\nabla_X,R[f]])=F(X)/F([\nabla_X,R[f]])$. It remains to be shown that $F$ preserves Smith commutators. This follows from exactness of $F$ and the fact that in a finitely cocomplete exact Mal'tsev category the Smith commutator is given by an explicit formula involving only finite limits and finite colimits.\endproof

\section{Affine morphisms and central reflections}\label{section3}

We have seen that the nilpotency tower of a $\sg$-pointed exact Mal'tsev category is a tower of \emph{central reflections}. In this section we establish a useful general property of central reflections in exact Mal'tsev categories, namely that the unit of a central reflection is pointwise affine. Since this property might be useful in other contexts as well, we first discuss possible weakenings of the notion of exactness.

\subsection{Quasi-exact and efficiently regular Mal'tsev categories}\label{effective}--\vspace{1ex}

An exact category is a regular category in which equivalence relations are \emph{effective}, i.e. arise as kernel relations of some morphism. In general, effective equivalence relations $R$ on $X$ have the property that the inclusion $R\into X\times X$ is a strong monomorphism. Equivalence relations with this property are called \emph{strong}. A regular category in which strong equivalence relations are effective is called \emph{quasi-exact}.

Any quasi-topos (cf. Penon \cite{Pen}) is quasi-exact so that there are plenty examples of quasi-exact categories which are not exact. There are also quasi-exact Mal'tsev categories which are not exact, as for instance the category of topological groups and continuous group homomorphisms. Further weakenings of exactness occur quite naturally as shown in the following chain of implications:
$$ \xymatrix@=20pt{
        \textrm{exact}\ar@{=>}[r] & \textrm{quasi-exact}\ar@{=>}[r]&\textrm{efficiently regular}\ar@{=>}[r]&\textrm{fibrational kernel relations} }$$

A category is \emph{efficiently regular} \cite{BR} if every equivalence relation $(X,S)$, which is a \emph{regular refinement} of an effective equivalence relation $(X,R)$, is itself effective. By regular refinement we mean any map of equivalence relations $(X,S)\to(X,R)$ inducing the identity on $X$ and a regular monomorphism $S\to R$.

We call a kernel relation $(X,R)$ \emph{fibrational} if for each fibrant map of equivalence relations $(Y,S)\to(X,R)$ the domain is an effective equivalence relation as well. According to Janelidze-Sobral-Tholen \cite{JST} a kernel relation $(X,R)$ is fibrational precisely when its quotient map $X\onto X/R$ has \emph{effective descent}, i.e. base-change along $X\onto X/R$ is a monadic functor. A regular category has thus fibrational kernel relations precisely when all regular epimorphisms have effective descent. The careful reader will observe that in all proofs of this section where we invoke efficient regularity we actually just need that the considered kernel relations are fibrational.

The second implication above follows from the fact that any regular monomorphism is a strong monomorphism, while the third implication follows from the facts that for any fibrant map of equivalence relations $f:(Y,S)\to(X,R)$ the induced map on relations $S\to f^*(R)$ is a regular (even split) monomorphism, and that in any regular category, effective equivalence relations are closed under inverse image.

\subsection{Fibration of points and essentially affine categories}\label{pointfibration}--\vspace{1ex}

Recall \cite{Bourn0} that for any category $\DD$, we denote by $\Pt(\DD)$ the category whose objects are split epimorphisms with chosen section (``genereralised points'') of $\DD$ and whose morphisms are natural transformations between such (compatible with the chosen sections), and that $\P_{\DD}:\Pt(\DD) \to \DD$ denotes the functor associating to a split epimorphism its codomain.

The functor $\P_{\DD}:\Pt(\DD)\to\DD$ is a fibration (the so-called \emph{fibration of points}) whenever $\DD$ has pullbacks of split epimorphisms. The $\P_\DD$-cartesian maps are precisely pullbacks of split epimorphisms. Given any morphism $f:X\rightarrow Y$ in $\DD$, base-change along $f$ with respect to the fibration $\P_\DD$ is denoted by $f^*:\Pt_{Y}(\DD) \rightarrow \Pt_X(\DD)$, and will be called \emph{pointed base-change} in order to distinguish it from the classical base-change $\DD/Y\to\DD/X$ on slices.

Pointed base-change $f^*$ has a \emph{left adjoint} pointed cobase-change $f_!$ if and only if pushouts along $f$ of split monomorphisms with domain $X$ exist in $\DD$. In this case pointed cobase-change along $f$ is given by precisely this pushout, cf. \cite{Bourn}. Accordingly, the unit (resp. counit) of the $(f_!,f^*)$-adjunction is an isomorphism precisely when for each natural transformation of split epimorphisms
$$ \xymatrix@=25pt{
       X' \ar[r]^{f'}   \ar@{->>}[d]^{r}  & Y' \ar@{->>}[d]^{r'}\\
       X \ar[r]_{f}  \ar@<1ex>[u]^{s} & Y \ar@<1ex>[u]^{s'}
                       }
  $$
the downward square is a pullback as soon as the upward square is a pushout (resp. upward square is a pushout as soon as the downward square is a pullback).

It is important that in a regular Mal'tsev category pointed base-change along a regular epimorphism is \emph{fully faithful}: in a pullback square like above with regular epimorphism $f:X\onto Y$, the upward-oriented square is automatically a pushout. This follows from fact that the induced morphism on kernel relations $R(s,s'):R[f]\to R[f']$ together with the diagonal $X'\to R[f']$ forms a strongly epimorphic cospan because the kernel relation $R[f']$ is the product of $R[f]$ and $X'$ in the fibre $\Pt_X(\DD)$, and all fibres are \emph{unital} in virtue of the Mal'tsev condition, cf. \cite{BB,Bourn0}.

%\begin{lma}\label{fullyfaithful}pointed base-change along a regular epimorphism is a fully faithful functor in any regular protomodular category.\end{lma}

%\proof Assume that, in the commutative square above, $f:X\onto Y$ is a regular epimorphism, and that the downward-oriented square is a pullback. We have to show that the upward-oriented square is a pushout. We know that the upward-oriented square is a pullback. It is a pushout if and only if the kernel relation $R[f]\dto X$ induces by composition with $s:X\into X'$ a pair of morphisms $R[f]\dto X'$ admitting $f':X'\onto Y'$ as coequaliser. This follows readily from the fact that the induced morphism $R(s,s'):R[f]\to R[f']$ and the section $s_0^{f'}:X'\to R[f']$ form a strongly epimorphic cospan by protomodularity, cf. \cite[Lemma 3.1.22]{BB}.\endproof

Recall \cite{Bourn0} that a category is called \emph{essentially affine} if pushouts of split monomorphisms and pullbacks of split epimorphisms exist, and moreover for any morphism $f:X\to Y$ the pointed base-change adjunction $(f_!,f^*)$ is an adjoint equivalence.

Additive categories with pullbacks of split epimorphisms are essentially affine. This follows from the fact that in this kind of category every split epimorphism is a projection, and every split monomorphism is a coprojection. Conversely, a \emph{pointed} essentially affine category is an additive category with pullbacks of split epimorphisms. Any slice or coslice category of an additive category with pullbacks of split epimorphisms is an example of an essentially affine category that is not pointed, and hence not additive. Therefore, the property of a morphism $f:X\to Y$ in $\DD$ to induce an adjoint equivalence $f^*:\Pt_{Y}(\DD) \rightarrow \Pt_X(\DD)$ expresses somehow a ``relative additivity'' of $f$. This motivates the following definition:

\begin{dfn}In a category $\DD$ with pullbacks of split epimorphisms, a morphism $f:X\to Y$ will be called \emph{$\P_{\DD}$-affine} when the induced pointed base-change functor $f^*:\Pt_\DD(Y)\to\Pt_\DD(X)$ is an equivalence of categories.

An \emph{affine extension} in $\DD$ is any regular epimorphism which is $\P_\DD$-affine.\end{dfn}

Clearly any isomorphism is $\P_{\DD}$-affine. It follows from the analogous property of equivalences of categories that for any composable morphisms $f,g$ in $\DD$, if two among $f$, $g$ and $gf$ are $\P_{\DD}$-affine then so is the third, i.e. $\P_{\DD}$-affine morphisms fulfill the so-called \emph{two-out-of-three} property.

It might be confusing that we use the term ``affine'' in two different contexts, namely as well for objects as well for morphisms. Although their respective definitions seem unrelated at first sight, this isn't the case. We will see in Proposition \ref{central} that every affine extension is a central extension, and it will follow from Theorem \ref{pointwiseaffine} that for a reflective subcategory $\CC$ of an efficiently regular Mal'tsev category $\DD$, every object of $\DD$ is $\CC$-affine if and only if every object of $\DD$ is an affine extension of an object of $\CC$. We hope that this justifies our double use of the term ``affine''.

\begin{prp}\label{central}
In any Mal'tsev category $\DD$, the kernel relation of a $\P_{\DD}$-affine morphism is central. In particular, each affine extension is a central extension.
\end{prp}
\proof
Let $f:X\to Y$ be an $\P_{\DD}$-affine morphism with kernel relation $$(p_0,p_1):R[f]\dto X$$and let$$(p_0^X,p_1^X):X\times X\dto X$$be the indiscrete equivalence relation on $X$ with section $s_0^X:X\to X\times X$. Since pointed base-change $f^*$ is an equivalence of categories, there is a split epimorphism $(r,s):Y'\rightleftarrows Y$ such that $f^*(r,s)=(p_0^X,s_0^X)$ and we get the right hand pullback of diagram
 $$ \xymatrix@=30pt{
        R[\check{f}]  \ar@<-1,ex>[d]_{R(p_0^X,r)} \ar@<-1,ex>[r]_{\check{p}_0}\ar@<+1,ex>[r]^{\check{p}_1}
       & X\times X \ar@<-1,ex>[d]_{p_0^X}\ar@{.>}@<1,ex>[d]^{p_1^X}\ar[l] \ar[r]^{\check{f}} & Y' \ar@<-1,ex>[d]_{r}\\
         R[f] \ar@<-1,ex>[r]_{p_0} \ar@<+1,ex>[r]^{p_1} \ar[u]_{R(s_0^X,s)} & X \ar[u]_{} \ar[l] \ar[r]_{f} & Y \ar[u]_{s}
                       }
    $$
in which the left hand side are the respective kernel relations. Therefore the left hand side consists of two pullbacks, and the map $p_1^X\check{p}_0$ produces the required connector between $R[f]$ and the indiscrete equivalence relation $\nabla_X$ on $X$.
\endproof

In particular, if an epireflection $I:\DD\to\CC$ of a Mal'tsev category $\DD$ has a pointwise $\P_{\DD}$-affine unit $\eta_X:X\to I(X)$, then it is a central reflection. The following converse will be essential in understanding nilpotency.

\begin{thm}\label{pointwiseaffine}A central reflection $I$ of an efficiently regular Mal'tsev category $\DD$ has a unit which is pointwise an affine extension. In particular, morphisms $f:X\to Y$ with invertible image $I(f):I(X)\to I(Y)$ are necessarily $\P_{\DD}$-affine.\end{thm}

\proof The second assertion follows from the first since $\P_{\DD}$-affine morphisms fulfill the two-out-of-three property and isomorphisms are $\P_{\DD}$-affine.

For the first assertion note that in a regular Mal'tsev category pointed base-change along a regular epimorphism is fully faithful, and hence $\eta_Y^*:\Pt_{I(Y)}(\DD)\to\Pt_Y(\DD)$ is a fully faithful functor. Corollary \ref{prop2} below shows that $\eta_Y^*$ is essentially surjective, hence $\eta_Y^*$ is an equivalence of categories for all objects $Y$ in $\DD$.\endproof

\subsection{Centralising double relations}Given a pair $(R,S)$ of equivalence relations on $Y$, we denote by $R\square S$ the inverse image of the reflexive relation $S\times S$ under $(p_0^R,p_1^R):R\rightarrowtail Y\times Y$. This defines a double relation
  $$ \xymatrix@=35pt{
      R \square S \ar@<-1,ex>[d]_{q_0^R}\ar@<+1,ex>[d]^{q_1^R} \ar@<-1,ex>[r]_{q_0^S}\ar@<+1,ex>[r]^{q_1^S}
     & S \ar@<-1,ex>[d]_{p_0^S}\ar@<+1,ex>[d]^{p_1^S} \ar[l]\\
       R \ar@<-1,ex>[r]_{p_0^R} \ar@<+1,ex>[r]^{
       p_1^R} \ar[u]_{} & Y
  \ar[u]_{} \ar[l]
                     }
  $$
actually the largest double relation relating $R$ and $S$. In set-theoretical terms, this double relation $R \square S$ corresponds to the subset of elements $(u,v,u',v')$ of $Y^4$ such that the relations $uRu',vRv',uSv,u'Sv'$ hold.
%$$ \xymatrix@=10pt{
%     u \ar@{.}[r]^S \ar@{.}[d]_R & v\ar@{.}[d]^R\\
%     u' \ar@{.}[r]_S & v'
%                    }
% $$

\begin{lma}\label{lemma2}Any split epimorphism $(r,s): X\rightleftarrows Y$ of a regular Mal'tsev category with epireflection $I$ and unit $\eta$ induces the following diagram
$$ \xymatrix@=35pt{
      R[\eta_{R[r]}] \ar@<-1,ex>[d]_{R(p_0^r,Ip_0^r)}\ar@<+1,ex>[d]^{R(p_1^r,Ip_1^r)} \ar@<-1,ex>[r]\ar@<+1,ex>[r]
     & R[r] \ar@<-1,ex>[d]_{p_0^r}\ar@<+1,ex>[d]^{p_1^r} \ar[l] \ar@{->>}[r]^{\eta_{R[r]}} & I(R[r]) \ar@<-1,ex>[d]_{Ip_0^r}\ar@<+1,ex>[d]^{Ip_1^r}\\
       R[\eta_X] \ar[d]_{R(r,Ir)} \ar@<-1,ex>[r]\ar@<+1,ex>[r] \ar[u]_{}  & X \ar@{->>}[r]_{\eta_X} \ar[u]_{}  \ar[d]_{r} \ar[l] & IX \ar[u]_{} \ar[d]_{Ir}\\
       R[\eta_Y] \ar@<-1,ex>[r]\ar@<+1,ex>[r] \ar@<-1ex>[u]_{R(s,Is)}  & Y \ar@{->>}[r]_{\eta_Y} \ar[l] \ar@<-1ex>[u]_{s} & IY \ar@<-1ex>[u]_{Is}
                       }
  $$
in which the rows and the two right vertical columns represent kernel relations. The left most column represents then the kernel relation of the induced map $R(r,Ir):R[\eta_X]\to R[\eta_Y]$, and we have $R[\eta_{R[r]}]=R[\eta_X]\square R[r]=R[R(r,Ir)]$.
\end{lma}
\proof
This follows from Proposition \ref{Diana} which shows that $I(R[r])$ may be identified with the kernel relation $R[I(r)]$ of the split epimorphism $Ir:IX\to IY$.
\endproof

For sake of simplicity a split epimorphism $(r,s):X\rightleftarrows Y$ is called a \emph{$\CC$-affine point over $Y$} whenever its domain $X$ is $\CC$-affine.
\begin{prp}\label{prop1}
Let $\CC$ be a reflective subcategory of an efficiently regular Mal'tsev category $\DD$ with reflection $I$ and unit $\eta$.

Any $\CC$-affine point $(r,s)$ over $Y$ is the image under $\eta_Y^*$ of a $\CC$-affine point $(\bar{r},\bar{s})$ over $IY$ such that both points have isomorphic reflections in $\CC$.\end{prp}
\proof
Since by Lemma \ref{lemma1.2} $\eta_X$ has a central kernel relation, the kernel relations $R[\eta_X]$ and $R[r]$ centralise each other. This induces a centralising double relation
 $$ \xymatrix@=35pt{
      R[\eta_X] \times_X R[r] \ar@<-1,ex>[d]\ar@{.>}@<+1,ex>[d] \ar@<-1,ex>[r]\ar@<+1,ex>[r]
     & R[r] \ar@<-1,ex>[d]_{p_0^r}\ar@{.>}@<+1,ex>[d]^{p_1^r} \ar[l]\\
       R[\eta_X] \ar@<-1,ex>[r] \ar@<+1,ex>[r] \ar[u]_{} & X
  \ar[u]_{} \ar[l]
                     }
  $$
which we consider as a fibrant split epimorphism of equivalence relations (disregarding the dotted arrows). Pulling back along the monomorphism
  $$ \xymatrix@=35pt{
       R[\eta_X]   \ar@<-1,ex>[r] \ar@<+1,ex>[r] & X  \ar[l]\\
        R[\eta_Y] \ar@<-1,ex>[r]\ar@<+1,ex>[r] \ar[u]^{R(s,Is)} & Y \ar[u]_{s} \ar[l]
                      }
   $$
yields on the left hand side of the following diagram
  $$ \xymatrix@=35pt{
        R_I[(r,s)]  \ar@<-1,ex>[d] \ar@<-1,ex>[r]\ar@<+1,ex>[r]
       & X \ar@<-1,ex>[d]_{r} \ar[l] \ar@{.>>}[r]^{q} & \bar X \ar@{.>}@<-1,ex>[d]_{\bar r}\\
         R[\eta_Y] \ar@<-1,ex>[r]\ar@<+1,ex>[r] \ar[u] & Y \ar[u]_{s} \ar[l] \ar@{->>}[r]_{\eta_Y} & IY \ar@{.>}[u]_{\bar s}
                       }
    $$another fibrant split epimorphism of equivalence relations. Since $\DD$ is efficiently regular, the upper equivalence relation is effective with quotient $q:X\onto{\bar X}$. We claim that the induced point $(\bar r,\bar s):\bar X\rightleftarrows IY$ has the required properties.

Indeed, the right square is a pullback by a well-known result of Barr and Kock, cf. \cite[Lemma A.5.8]{BB}, so that $\eta_Y^*(\bar{r},\bar{s})=(r,s)$. The centralising double relation $R[\eta_X] \times_X R[r]$ is coherently embedded in the double relation $R[\eta_X]\square R[r]$ of Lemma \ref{lemma2}, cf. \cite[Proposition 2.6.13]{BB}.  This induces an inclusion $ R_I[(r,s)]\rightarrowtail R[\eta_X]$ and hence a morphism $\phi: \bar X\twoheadrightarrow IX$ such that $\phi q=\eta_X$.

According to Lemma \ref{lemma1.1}, we get an isomorphism $Iq:IX\cong I\bar{X}$ compatible with the units. The kernel relation $R[\eta_{\bar{X}}]$ is thus the direct image of the central kernel relation $R[\eta_X]$ under the regular epimorphism $q:X\to\bar{X}$ and as such central as well. In particular, $(\bar{r},\bar{s})$ is a $\,\CC$-affine point with same reflection in $\CC$ as $(r,s)$.
\endproof

\begin{cor}\label{prop2}
For an efficiently regular Mal'tsev category $\DD$ with central reflection $I$ and unit $\eta$, pointed base-change $\eta_Y^*:\Pt_{I(Y)}(\DD)\to\Pt_Y(\DD)$ is essentially surjective.\end{cor}
\proof
Since the reflection is central, Corollary \ref{cor1.3} shows that Proposition \ref{prop1} applies to the whole fibre $\Pt_Y(\DD)$ whence essential surjectivity of $\eta_Y^*$.\endproof

\subsection{Affine extensions in efficiently regular Mal'tsev categories}--\vspace{1ex}

A functor $G:\EE\to\EE'$ is called \emph{saturated on quotients} if for each object $A$ in $\EE$ and each strong epimorphism $g':G(A)\to B'$ in $\EE'$, there exists a strong epimorphism $g:A\to B$ in $\EE$ such that $G(g)$ and $g'$ are isomorphic under $G(A)$. Note that a right adjoint functor $G:\EE\to\EE'$ is essentially surjective whenever it is saturated on quotients and each object $B'$ of $\EE'$ is the quotient of an object $B''$ for which the unit $\eta_{B''}:B''\to GF(B'')$ is invertible.

\begin{lma}\label{lem1}In an efficiently regular Mal'tsev category, pointed base-change along a regular epimorphism is saturated on quotients.\end{lma}
\proof Let $f:X\onto Y$ be a regular epimorphism, let $(r,s)$ be a point over $Y$, and $l:f^*((r,s))\onto(r',s')$ be a quotient map over $X$. Consider the following diagram
$$ \xymatrix@=25pt{
       R[f'] \ar@{->>}[rd] \ar@{->>}[dd] \ar@<-1,ex>[rr]_{p_0^{f'}} \ar@<+1,ex>[rr]^{p_1^{f'}}   && X' \ar@{->>}[rd]^{l} \ar@{->>}[rr]^{f'}   \ar@{->>}[dd]_>>>>>>{f^*(r)} \ar[ll] && Y' \ar@{->>}[dd]_>>>>>>{r} \ar@{.>>}[rd]^{\bar l} \\
        & S \ar@{->>}[ld] \ar@<-1,ex>[rr]_{} \ar@<+1,ex>[rr]^{} &&  X'' \ar@{.>>}[rr]^<<<<<<{f''} \ar[ll] \ar@{->>}[ld]^{r'} && Y'' \ar@{.>>}[ld]^{\rho}\\
       R[f] \ar@<-1,ex>[rr]_{p_0^f} \ar@<+1,ex>[rr]^{p_1^f} \ar@<-1ex>[uu]_{}  && X \ar@{->>}[rr]_{f} \ar[ll] \ar@<-1ex>[uu]_<<<<<<{f^*(s)} && Y \ar@<-1ex>[uu]_<<<<<<s
                       }
  $$in which the right hand side square is a pullback, and the left hand side is defined by factoring the induced map on kernel relations $R[f']\to R[f]$ through the direct image $S=l(R[f'])$ under $l$. Since the right square is a pullback, the left square represents a fibrant split epimorphism of equivalence relations. The factorisation of this fibrant morphism induced by $l$ yields two fibrant maps of equivalence relations, cf. Lemma \ref{propbr}. Note that the second $(X'',S)\to(X,R[f])$ is a fibrant \emph{split} epimorphism. Efficient regularity implies then that the equivalence relation $S$ is effective with quotient $f'':U''\onto V''$, defining a point $(\rho,\sigma)$ over $Y$.

  This induces (by a well-known result of Barr and Kock, cf. \cite[Lemma A.5.8]{BB}) a decomposition of the right pullback into two pullbacks. The induced regular epimorphism $\bar{l}:(r,s)\onto(\rho,\sigma)$ has the required properties, namely $f^*(\bar{l})=l$.\endproof

\begin{prp}\label{caraffine}
In an efficiently regular Mal'tsev category with binary sums, a regular epimorphism $f:X\onto Y$ is an affine extension if and only if for each object $Z$ either of the following two diagrams
$$ \xymatrix@=25pt{X+Z \ar@<-1,ex>@{->>}[d]_{\pi^Z_X} \ar@{->>}[r]^{f+Z} & Y+Z \ar@<-1,ex>@{->>}[d]_{\pi^Z_Y}&&
        X+Z \ar@{->>}[d]_{\theta_{X,Z}} \ar@{->>}[r]^{f+Z} & Y+Z \ar@<-1,ex>@{->>}[d]^{\theta_{Y,Z}}\\
    X \ar@{->>}[r]_{f}\ar@{->}[u]_{\iota^Z_X} & Y \ar@{->}[u]_{\iota^Z_Y}&&     X\times Z \ar@{->>}[r]_{f\times Z}\ & Y\times Z
                       }
    $$is a downward-oriented pullback square.\end{prp}
\proof
If $f$ is an affine extension then the downward-oriented left square is a pullback because the upward-oriented left square is a pushout, cf. Section \ref{pointfibration}. Moreover, the outer rectangle of the following diagram

    $$ \xymatrix@=20pt{
          X+Z \ar[rr]^{f+Z} \ar[d]_{\theta_{X,Z}} && Y+Z \ar[d]^{\theta_{Y,Z}}\\
          X\times Z \ar[rr]_{f\times Z} \ar[d]_{p_X}  && Y\times Z \ar[d]^{p_Y}\\
         X \ar[rr]_{f} && Y
                           }
     $$
is a pullback if and only if the upper square is a pullback, because the lower square is always a pullback.

Assume conversely that the downward oriented left square is a pullback. In a regular Mal'tsev category, pointed base-change along a regular epimorphism is fully faithful so that $f$ is affine whenever $f^*$ is essentially surjective. Lemma \ref{lem1} shows that in an efficiently regular Mal'tsev category $f^*$ is saturated on quotients. It suffices thus to show that in the fibre over $X$ each point is the quotient of a point for which the unit of the pointed base-change adjunction is invertible.

Since for each object $Z$, the undotted downward-oriented square

$$ \xymatrix@=35pt{
        X+Z \ar@<-1,ex>[d]_{\pi_{X}^Z} \ar@{.>}@<+1,ex>[d]^{\langle 1_X,r\rangle} \ar@{->>}[r]^{f+Z} & Y+Z \ar@<-1,ex>[d]_{\pi_{Y}^Z}\ar@{.>}@<+1,ex>[d]^{\langle 1_Y,fr\rangle} \\
        X \ar@{->>}[r]_{f}\ar[u] & Y \ar[u]
                       }
    $$
 is a pullback, the dotted downward-oriented square (which is induced by an arbitrary morphism $r:Z\to X$) is a pullback as well. This holds in any regular Mal'tsev category, since the whole diagram represents a natural transformation of reflexive graphs, cf. \cite{Bourn0}. It follows that the point $(\langle 1_X,r\rangle,\iota^Z_{X}):X+Z\rightleftarrows X$ has an invertible unit with respect to the pointed base-change adjunction $(f_!,f^*)$.

Now, an arbitrary point $(r,s):Z\rightleftarrows X$ can be realised as a quotient
   $$ \xymatrix@=25pt{
           Z \ar[dr]^{r} && X+Z \ar[dl]^{\langle 1_X,r\rangle} \ar@{->>}[ll]_{\langle s,1_Z\rangle} \\
           & X \ar@<1,ex>[ul]^{s}  \ar@<1,ex>[ur]^<<<<<<{\iota^Z_{X}}
                          }
       $$of the latter point $(\langle 1_X,r\rangle,\iota^Z_{X}):X+Y\rightleftarrows X$ with invertible unit.\endproof
 We end this section with several properties of affine extensions in semi-abelian categories. They will only be used in Section \ref{degreen}.

\begin{prp}\label{caraffine2}In a semi-abelian category, a regular epimorphism $f:X\onto Y$ is an affine extension if and only if either of the following conditions is satisfied:\begin{itemize}\item[(a)] for each object $Z$, the induced map $f\diamond Z:X\diamond Z\to Y\diamond Z$ is invertible, where $X\diamond Z$ stands for the kernel of $\,\theta_{X,Z}:X+Z\onto X\times Z$;\item[(b)]every pushout of $f$ along a split monomorphism is a pullback.\end{itemize}\end{prp}

\proof That condition (a) characterises affine extensions follows from Proposition \ref{caraffine} and protomodularity. The necessity of condition (b) follows from Section \ref{pointfibration}. The sufficiency of condition (b) follows from the ``pullback cancellation property'' in semi-abelian categories, cf. \cite[Proposition 4.1.4]{BB}.\endproof

\begin{rmk}\label{diamond}This product $\,X\diamond Z$ is often called the \emph{co-smash product} of $X$ and $Z$, since it is the dual of the \emph{smash product} as investigated by Carboni-Janelidze \cite{CJ} in the context of lextensive categories. The co-smash product $X\diamond Z$ coincides in semi-abelian categories with the second cross-effect $cr_2(X,Z)$ of the identity functor, cf. Definition \ref{ncube} and \cite{MM,HL,HV}. Since the co-smash product is in general not associative (cf. \cite{CJ}), parentheses should be used with care.\end{rmk}

\begin{prp}\label{affinequotient}Let $Y\lto W\otl Z$ and $\bar{Y}\lto\bar{W}\otl\bar{Z}$ be cospans in the fibre $\Pt_X(\DD)$ of a semi-abelian category $\DD$. Let $f:Y\onto\bar{Y},\,g:Z\onto\bar{Z},\,h:W\onto\bar{W}$ be affine extensions in $\DD$ inducing a map of cospans in $\Pt_X(\DD)$. Assume furthermore that the first cospan realises $W$ as the binary sum of $\,Y$ and $Z$ in $\Pt_X(\DD)$.

Then the second cospan realises $\bar{W}$ as the binary sum of $\,\bar{Y}$ and $\bar{Z}$ in $\Pt_X(\DD)$ if and only if the kernel cospan $K[f]\lto K[h]\otl K[g]$ is strongly epimorphic in $\DD$.\end{prp}

\proof Let us consider the following commutative diagram
$$
\xymatrix@=20pt{X\;\ar@{>->}[r]^i\ar@{>->}[d]_j&
Y\ar@{->>}[r]^f\ar@{>->}[d] & \bar{Y} \ar@{>->}[d]\\Z\ar@{->>}[d]_g\;\ar@{>->}[r]&W\ar@{->>}[r]^{h_1}\ar@{->>}[d]_{h_2}\ar@{->>}[rd]^h&W_1\ar@{->>}[d]\\\bar{Z}\ar@{>->}[r]&W_2\ar@{->>}[r]&\bar{W}
              }
$$in which $i$ (resp. $j$) denotes the section of the point $Y$ (resp. $Z$) over $X$, and all little squares except the lower right one are pushouts. It follows that the outer square is a pushout (i.e. $\bar{W}=\bar{Y}+_X\bar{Z}$) if and only if the lower right square is a pushout. According to Carboni-Kelly-Pedicchio \cite[Theorem 5.2]{CKP} this happens if and only if the kernel relation $R[h]$ is the join of the kernel relations $R[h_1]$ and $R[h_2]$. In a semi-abelian category this is the case if and only if the kernel $K[h]$ is generated as normal subobject of $W$ by the kernels $K[h_1]$ and $K[h_2]$, resp. (since $h_1$ and $h_2$ are affine extensions) by the kernels $K[f]$ and $K[g]$, cf. Proposition \ref{caraffine2}b.

Now, $h$ is also an affine extension so that by Proposition \ref{central}, the kernel $K[h]$ is a central subobject of $W$. In particular, any subobject of $K[h]$ is central and normal in $W$ (cf. the characterisation of normal subobjects in semi-abelian categories by Mantovani-Metere \cite[Theorem 6.3]{MM}). Therefore, generating $K[h]$ as normal subobject of $W$ amounts to the same as generating $K[h]$ as subobject of $W$.\endproof

\section{Aspects of nilpotency}\label{section4}

Recall that a morphism is called $n$-fold centrally decomposable if it is the composite of $n$ morphisms with central kernel relation.

For consistency, a monomorphism is called $0$-fold centrally decomposable, and an isomorphism a $0$-fold central extension.

\begin{prp}\label{discrepancy}For all objects $X,Y$ of a $\sg$-pointed $n$-nilpotent Mal'tsev category, the comparison map $\theta_{X,Y}:X+Y\to X\times Y$ is $(n-1)$-fold centrally decomposable.
\end{prp}
\proof
In a pointed $n$-nilpotent Mal'tsev category, each object maps to an abelian group object through an $(n-1)$-fold centrally decomposable morphism, cf. Proposition \ref{nilpotentcategory}. Since the codomain of such a morphism $\phi_{X,Y}:X+Y\to A$ is an abelian group object, the restrictions to the two summands commute and $\phi_{X,Y}$ factors
$$ \xymatrix@=15pt{
        X+Y \ar@{->>}[rr]^{\theta_{X,Y}} \ar[rd]_{\phi_{X,Y}} &&  X \times Y \ar[ld]^{\psi_{X,Y}} \\
         & A
                       }
    $$
so that $\theta_{X,Y}$ is $(n-1)$-fold centrally decomposable by Lemma \ref{lefterase}.\endproof

\begin{prp}\label{abelianization}For any finitely cocomplete regular pointed Mal'tsev category, the following pushout square
$$ \xymatrix@=20pt{
               X+X \ar@{->>}[d]_{\langle 1_X, 1_X\rangle} \ar@{->>}[r]^{\theta_{X,X}} & X\times X \ar@{->>}[d]\\
               X   \ar@{->>}[r] & A(X)
                             }
          $$defines the abelianisation $A(X)$ of $X$. In particular, the lower row can be identified with the unit $\eta_X^1:X\to I^1(X)$ of the strong epireflection of Theorem \ref{noBirkhoff}.\end{prp}
\proof The first assertion follows by combining \cite[Proposition 1.7.5, Theorems 1.9.5 and 1.9.11]{BB} with the fact that pointed Mal'tsev categories are strongly unital in the sense of the second author, cf. \cite[Corollary 2.2.10]{BB}. The second assertion expresses the fact that $X\onto A(X)$ and $X\onto I^1(X)$ share the same universal property.\endproof

\begin{thm}\label{ndiscrepancy}A $\sg$-pointed exact Mal'tsev category is $n$-nilpotent if and only if for all objects $X,Y$ the comparison map $\theta_{X,Y}:X+Y\to X\times Y$ is an $(n-1)$-fold central extension.\end{thm}
\proof By Proposition \ref{discrepancy} $n$-nilpotency implies that $\theta_{X,Y}$ is an $(n-1)$-fold central extension. For the converse, consider the pushout square of Proposition \ref{abelianization}, which is regular by Corollary \ref{stpushout}. The unit $\eta_X^1:X\to I^1(X)$ is thus an $(n-1)$-fold central extension by Proposition \ref{trucc} so that all objects are $n$-nilpotent.\endproof

\begin{cor}\label{abelian}For a $\sg$-pointed exact Mal'tsev category, the following three properties are equivalent:\begin{itemize}\item[(a)]the category is $1$-nilpotent;\item[(b)]the category is linear (cf. Definition \ref{ncube});\item[(c)]the category is abelian.\end{itemize}\end{cor}

\proof The equivalence of (a) and (b) follows from Theorem \ref{ndiscrepancy}. The equivalence of (b) and (c) follows from the fact that a $\sg$-pointed Mal'tsev category is additive if and only if it is linear (cf. \cite[Theorem 1.10.14]{BB}) together with the well-known fact (due to Miles Tierney) that abelian categories are precisely the additive categories among exact categories.\endproof

\begin{thm}\label{2discrepancy}For a $\sg$-pointed exact Mal'tsev category, the following five properties are equivalent:
\begin{itemize}\item[(a)]all objects are $2$-nilpotent;
\item[(b)]for all $X$, abelianisation $\eta^1_X:X\to I^1(X)$ is a central extension;
\item[(b$'$)]for all $X$, abelianisation $\eta^1_X:X\to I^1(X)$ is an affine extension;
\item[(c)]for all $X,Y$, the map $\theta_{X,Y}:X+Y\to X\times Y$ is a central extension;
\item[(c$'$)]for all $X,Y$, the map $\theta_{X,Y}:X+Y\to X\times Y$ is an affine extension.\end{itemize}\end{thm}

\proof Properties (a) and (b) are equivalent by definition of $2$-nilpotency. Theorem \ref{ndiscrepancy} shows that (b) and (c) are equivalent. Theorem \ref{pointwiseaffine} and Proposition \ref{central} imply that (b) and (b$'$) are equivalent.

Finally, since the first Birkhoff reflection $I^1$ preserves binary sums and binary products (cf. Proposition \ref{Diana}), we get $I^1(\theta_{X,Y})=\theta_{I^1(X),I^1(Y)}$ which is invertible in the subcategory of $1$-nilpotent objects by Proposition \ref{discrepancy}. It follows that under assumption (a), the map $\theta_{X,Y}$ is an affine extension by Theorem \ref{pointwiseaffine}, which is property (c$'$). Conversely, (c$'$) implies $(c)$ by Proposition \ref{central}.\endproof

\subsection{Niltensor products}

In order to extend Proposition \ref{2discrepancy} to higher $n$ we introduce here a new family of binary tensor products, called \emph{niltensor products}.

For any finitely cocomplete pointed regular Mal'tsev category $(\DD,\star_\DD)$ the \emph{$n$-th niltensor product} $X\otimes_nY$ is defined by factorizing the comparison map $\theta_{X,Y}$ universally into a regular epimorphism $\theta_{X,Y}^n:X+Y\to X\otimes_nY$ followed by an $(n-1)$-fold central extension
 $$ \xymatrix@=4pt{
 && X\otimes_nY \ar@{->>}[rr]^{\omega^{n-1}_{X,Y}}& & X\otimes_{n-1}Y  \ar@{.}[rrrr] && & & X\otimes_3Y \ar@{->>}[rr]^{\omega^2_{X,Y}}  & & X\otimes_2Y \ar@{->>}[rrrrddd]^{\omega^1_{X,Y}} &&\\
 &&&&&&&\\
 &&&&&&&\\
X+Y   \ar@{->>}[rrrrrrrrrrrrrr]_{\theta_{X,Y}} \ar@{->>}[rruuu]^{\theta^n_{X,Y}}  \ar@{->>}[rrrrrrrrrruuu]_{\theta^{2}_{X,Y}}\ar@{->>}[rrrruuu]_>>>>>>{\theta^{n-1}_{X,Y}}&&&&&&&&&&&&&&      X\times Y
                    }
    $$
as provided by Proposition \ref{univncentral}. This $n$-th niltensor product is symmetric and has $\star_\DD$ as unit, but it does not seem to be associative in general.%$X^{\frac{n}{\star}}\,\star_\DD=X=\star_\DD\,\,{}^{\frac{n}{\star}}X$
\begin{prp}\label{trucccc}
In a $\sg$-pointed exact Mal'tsev category, the following diagram is an iterated pushout diagram
 $$ \xymatrix@=3pt{
 && X\otimes_nX \ar@{->>}[rr]^{\omega^{n-1}_{X,X}} \ar@{->>}[dddddddd] & & X\otimes_{n-1}X  \ar@{.}[rrrr] \ar@{->>}[dddddddd] && & & X\otimes_3X \ar@{->>}[rr]^{\omega^2_{X,X}}\ar@{->>}[dddddddd]  & & X\otimes_2X \ar@{->>}[rrrrddd]^{\omega^1_{X,X}}\ar@{->>}[dddddddd] &&\\
 &&&&&&&\\
 &&&&&&&\\
X+X \ar@{->>}[dd]  \ar@{->>}[rruuu]^{\theta^n_{X,X}} \ar@{.>>}[rrrrrrrrrrrrrr]_{\theta_{X,X}} &&&&&&&&&&&&&& \ar@{->>}[dd]  X\times X\\
&&&&&&&\\
X  \ar@{.>>}[rrrrrrrrrrrrrr]_{\eta_X} \ar@{->>}[rrddd]_{\eta^n_X} &&&&&&&&&&&&&& I^1(X)\\
 &&&&&&&\\
 &&&&&&&\\
&& I^{n}(X) \ar@{->>}[rr]^{} & & I^{n-1}(X) \ar@{.}[rrrr]  &&&& I^{3}(X) \ar@{->>}[rr]_{}  & & I^{2}(X)\ar@{->>}[rrrruuu]_{} &&\\
                    }
    $$
where the left vertical map is the folding map $\langle 1_X,1_X\rangle:X+X\to X.$
\end{prp}
\proof
This follows from Corollary \ref{stpushout} and Propositions \ref{trucc} and \ref{abelianization}.
\endproof

\begin{thm}\label{nnill}For a $\sg$-pointed exact Mal'tsev category, the following five properties are equivalent:
\begin{itemize}\item[(a)]all objects are $n$-nilpotent;\item[(b)]for all $X$, the $(n-1)$th unit $\eta^{n-1}_X:X\onto I^{n-1}(X)$ is a central extension;\item[(b$'$)]for all $X$, the $(n-1)$th unit $\eta^{n-1}_X:X\onto I^{n-1}(X)$ is an affine extension;\item[(c)]for all $X,Y$, the map $\theta^{n-1}_{X,Y}:X+Y\to X\otimes_{n-1}Y$ is a central extension;\item[(c$'$)]for all $X,Y$, the map $\theta^{n-1}_{X,Y}:X+Y\to X\otimes_{n-1}Y$ is an affine extension.
\end{itemize}
\end{thm}
\proof Properties (a) and (b) are equivalent by definition of $n$-nilpotency. Theorem \ref{ndiscrepancy} shows that (a) implies (c) while Proposition \ref{trucccc} shows that (c) implies (b). Therefore, (a), (b) and (c) are equivalent. Proposition \ref{central} and Theorem \ref{pointwiseaffine} imply that (b) and (b$'$) are equivalent.

Finally, Proposition \ref{trucc} implies that the Birkhoff reflection $I^{n-1}$ takes the comparison map $\theta^{n-1}_{X,Y}$ to the corresponding map $\theta^{n-1}_{I^{n-1}(X),I^{n-1}(Y)}$ for the $(n-1)$-nilpotent objects $I^{n-1}(X)$ and $I^{n-1}(Y)$. Since the $(n-1)$-nilpotent objects form an $(n-1)$-nilpotent Birkhoff subcategory, Theorem \ref{ndiscrepancy} shows that the latter map must be invertible; therefore, (a) implies (c$'$) by Theorem \ref{pointwiseaffine}. Conversely, (c$'$) implies (c) by Proposition \ref{central}.\endproof

\begin{dfn}\label{pseudoadditive}A $\sg$-pointed Mal'tsev category is said to be \emph{pseudo-additive} (resp. \emph{pseudo-$n$-additive}) if for all $X,Y,$ the map $\theta_{X,Y}:X+Y\onto X\times Y$ (resp. $\theta^n_{X,Y}:X+Y\onto X\otimes_nY$) is an affine extension.\end{dfn}

\begin{prp}\label{ps} A $\sg$-pointed exact Mal'tsev category is pseudo-additive (i.e. $2$-nilpotent) if and only if the following diagram
$$ \xymatrix@=25pt{
          (X+Y)+Z \ar@<-1,ex>@{->>}[d]_{\theta_{X+Y,Z}} \ar@{->>}[r]^{\theta_{X,Y}+Z} & (X\times Y)+Z \ar@<-1,ex>@{->>}[d]^{\theta_{X\times Y,Z}}\\
        (X+Y)\times Z \ar@{->>}[r]_{\theta_{X,Y}\times Z}\ & (X\times Y)\times Z
                       }
    $$is a pullback for all objects $X,Y,Z$.
\end{prp}
\proof This follows from Theorem \ref{2discrepancy} and Proposition \ref{caraffine}.\endproof

\begin{prp}\label{psn}A $\sg$-pointed exact Mal'tsev category is pseudo-$(n-1)$-additive (i.e. $n$-nilpotent) if and only if the following diagram
$$ \xymatrix@=20pt{
      (X+Y)+Z  \ar@{->>}[rr]^{\theta_{X+Y,Z}} \ar[d]_{\theta^{n-1}_{X,Y}+ Z} && (X+Y)\times Z\ar[d]^{\theta^{n-1}_{X,Y}\times Z} \\
   (X\otimes_{n-1}Y)+Z  \ar@{->>}[rr]_{\theta_{X\otimes_{n-1}Y,Z}}  && (X\otimes_{n-1}Y) \times Z
                       }
    $$is a pullback for all objects $X,Y,Z$.
\end{prp}
\proof This follows from Theorem \ref{nnill} and Proposition \ref{caraffine}.\endproof

We end this section with a general remark about the behaviour of $n$-nilpotency under slicing and passage to the fibres. Note that any left exact functor between Mal'tsev categories preserves central equivalence relations, morphisms with central kernel relation, and consequently $n$-nilpotent objects.

\begin{prp}\label{slice}If $\,\DD$ is an $n$-nilpotent Mal'tsev category, then so are any of its slice categories $\DD/Y$ and of its fibres $\Pt_Y(\DD)$.
\end{prp}
\proof The slices $\DD/Y$ of a Mal'tsev category $\DD$ are again Mal'tsev categories. Moreover, base-change
$\omega_Y^*:\DD\to\DD/Y$ is a left exact functor so that the objects of $\DD/Y$ of the form $\omega_Y^*(X)=p_Y:Y\times X\to Y$ are $n$-nilpotent provided $\DD$ is an $n$-nilpotent Mal'tsev category. We can conclude with Proposition \ref{epireflective} by observing that \emph{any} object $f:X\to Y$ of $\DD/Y$ may be considered as a subobject
$$ \xymatrix@=20pt{
         {X\;} \ar[rd]_{f}\ar@{>->}[rr]^{(f,1_X)} && Y\times X \ar[dl]^{p_Y}\\
         & Y
                       }
    $$of $\omega_Y^*(X)$ in $\,\DD/Y$.
The proof for the fibres is the same as for the slices, since any object $(r,s)$ of $\Pt_Y(\DD)$ may be considered as a subobject
$$
\xymatrix@=30pt{
{X\;} \ar@{>->}[rr]^{(r,1_X)}\ar@<-1ex>[rd]_{r} && Y\times X  \ar@<-1ex>[dl]_{p_Y}\\
 & Y \ar[lu]_>>>>>>s \ar[ur]_{(1_Y,s)}
              }
$$
of the projection $p_Y:Y\times X\to Y$ splitted by $(1_Y,s):Y\to Y\times X.$\endproof

\section{Quadratic identity functors}\label{linearquadratic}

We have seen that $1$-nilpotency has much to do with linear identity functors (cf. Corollary \ref{abelian}). We now investigate the relationship between $2$-nilpotency and quadratic identity functors, and below in Section $6$, the relationship between $n$-nilpotency and identity functors of degree $n$. While a linear functor takes binary sums to binary products, a quadratic functor takes certain cubes constructed out of triple sums to limit cubes. This is the beginning of a whole hierarchy assigning degree $\leq n$ to a functor whenever the functor takes certain $(n+1)$-dimensional cubes constructed out of iterated sums to limit cubes.

This definition of \emph{degree} of a functor is much inspired by Goodwillie \cite{G2} who described polynomial approximations of a homotopy functor in terms of their behaviour on certain cubical diagrams. Eilenberg-Mac Lane \cite{EM} defined the degree of a functor with values in an \emph{abelian} category by a vanishing condition of so-called \emph{cross-effects}. Our definition of degree does not need cross-effects. Yet, a functor with values in a semi-abelian (or homological \cite{BB}) category is of degree $\leq n$ precisely when all its cross-effects of order $n+1$ vanish, cf. Corollary \ref{crosseffects}. Our cubical cross-effects agree up to isomorphism with those of Hartl-Loiseau \cite{HL} and Hartl-Van der Linden \cite{HV}, which are defined as kernel intersections.

There are several other places in literature where degree $n$ functors, especially quadratic functors, are studied in a non-additive context, most notably Baues-Pirashvili \cite{BP}, Johnson-McCarthy \cite{JMc} and Hartl-Vespa \cite{HVe}. In all these places, the definition of a degree $n$ functor is based on a vanishing condition of cross-effects, closely following the original approach of Eilenberg-Mac Lane. It turned out that for us Goodwillie's cubical approach to functor calculus was more convenient. Our Definition \ref{nfolded} of an $n$-folded object and the resulting characterisation of degree $n$ identity functors in terms of $n$-folded objects (cf. Proposition \ref{nadditive}) rely in an essential way on cubical combinatorics.

\subsection{Degree and cross-effects of a functor}--\vspace{1ex}

An \emph{$n$-cube} in a category $\EE$ is given by a functor $\Xi:[0,1]^n\to\EE$ with domain the $n$-fold cartesian product $[0,1]^n$ of the arrow category $[0,1]$.

The category $[0,1]$ has two objects $0$,$1$ and exactly one non-identity arrow $0\to 1$. Thus, an $n$-cube in $\EE$ is given by objects $\Xi(\eps_1,\dots,\eps_n)$ in $\EE$ with $\eps_i\in\{0,1\}$, and arrows $\xi_{\eps_1,\dots,\eps_n}^{\eps'_1,\dots,\eps'_n}:\Xi(\eps_1,\dots,\eps_n)\to\Xi(\eps'_1,\dots,\eps'_n)$ in $\EE$, one for each arrow in $[0,1]^n$, which compose in an obvious way.

To each $n$-cube $\Xi$ we associate a \emph{punctured} $n$-cube $\check{\Xi}$ obtained by restriction of $\Xi$ to the full subcategory of $[0,1]^n$ spanned by the objects $(\eps_1,\dots,\eps_n)\not=(0,\dots,0)$.

\begin{dfn}\label{ncube}Let $(\EE,\star_\EE)$ be a $\sg$-pointed category. For each $n$-tuple of objects $(X_1,\dots,X_n)$ of $\,\EE$ we denote by $\Xi_{X_1,\dots,X_n}$ the following $n$-cube:\begin{itemize}\item $\Xi_{X_1,\dots,X_n}(\eps_1,\dots,\eps_n)=X_1(\eps_1)+\cdots+X_n(\eps_n)$\\ with $X(0)=X$ and $X(1)=\star_\EE$;\item $\xi_{\eps_1,\dots,\eps_n}^{\eps'_1,\dots,\eps'_n}=j_{\eps_1}^{\eps'_1}+\cdots+j_{\eps_n}^{\eps'^n}$\\ where $j_\eps^{\eps'}$ is the identity if $\eps=\eps'$, resp. the null morphism if $\eps\not=\eps'$.\end{itemize}

\noindent A functor of $\sg$-pointed categories $F:(\EE,\star_\EE)\to(\EE',\star_{\EE'})$ is called \emph{of degree $\leq n$} if\begin{itemize}\item $F(\star_\EE)\cong\star_{\EE'}$;\item for each $(n+1)$-cube $\Xi_{X_1,\dots,X_{n+1}}$ in $\EE$, the image-cube $F\circ\Xi_{X_1,\dots,X_{n+1}}$ is a \emph{limit-cube} in $\EE'$, i.e. $F(X_1+\dots+X_{n+1})$ may be identified with the limit of the punctured image-cube $F\circ\check{\Xi}_{X_1,\dots,X_{n+1}}$.\end{itemize}

%\noindent An object $Z$ in $\EE$ is said to be of degree $\leq n$ when the pointed cobase-change $\alpha_Z^*:\EE\to Pt_Z\EE$ is polynomial of degree $\leq n$.

\noindent A functor of degree $\leq 1$ (resp. $\leq 2$) is called \emph{linear} (resp. \emph{quadratic}).

\noindent A $\sg$-pointed category is called \emph{linear} (resp. \emph{quadratic}) if its identity functor is.

\noindent If $\,\EE'$ has pullbacks, the limit over the punctured image-cube is denoted$$P^F_{X_1,\dots,X_{n+1}}=\varprojlim_{[0,1]^{n+1}-{(0,\dots,0)}} F\circ\check{\Xi}_{{X_1},\dots,X_{n+1}}$$and the associated comparison map is denoted $$\theta^F_{X_1,\dots,X_{n+1}}:F(X_1+\cdots+X_{n+1})\to P^F_{X_1,\dots,X_{n+1}}.$$

The \emph{$(n+1)$-st cross-effects} of the functor $F:\EE\to\EE'$ are the \emph{total kernels} of the image-cubes, i.e. the kernels of the comparison maps $\theta^F_{X_1,\dots,X_{n+1}}$: $$cr_{n+1}^F(X_1,\dots,X_{n+1})=K[\theta^F_{X_1,\dots,X_{n+1}}].$$\end{dfn}

If $F$ is the identity functor, the symbol $F$ will be dropped from the notation. A functor $F:\EE\to\EE'$ has degree $\leq n$ if and only if for all $(n+1)$-tuples $(X_1,\dots,X_{n+1})$ of objects of $\EE$, the comparison maps $\theta^F_{X_1,\dots,X_{n+1}}$ are invertible.\vspace{1ex}

Our \emph{total-kernel} definition of the cross-effect $cr_{n+1}^F(X_1,\dots,X_{n+1})$ is directly inspired by Goodwillie \cite[pg. 676]{G2} but agrees up to isomorphism with the \emph{kernel intersection} definition of Hartl-Loiseau \cite{HL} and Hartl-Van der Linden \cite{HV}. Their kernel intersection is dual to the $(n+1)$-fold smash product of Carboni-Janelidze \cite{CJ}, cf. Remark \ref{diamond} and also Remark \ref{history}, where the duality between cross-effects and smash-products is discussed in more detail.\vspace{1ex}

Indeed, each of the $n+1$ ``contraction morphisms''$$\pi^F_{\hat{X}_i}:F(X_1+\cdots+X_{n+1})\to F(X_1+\cdots+\widehat{X_i}+\cdots+X_{n+1}),\,1\leq i\leq n+1,$$factors through $\theta^F_{X_1,\dots,X_{n+1}}$ so that we get a composite morphism$$r_F:F(X_1+\cdots+X_{n+1})\overset{\theta^F_{X_1,\dots,X_{n+1}}}{\longrightarrow} P^F_{X_1,\dots,X_{n+1}}\into\prod_{i=1}^{n+1}F(X_1+\cdots+\widehat{X_i}+\cdots+X_{n+1})$$
embedding the limit construction $P^F_{X_1,\dots,X_{n+1}}$ into an $(n+1)$-fold cartesian product. Therefore, the kernel of $\theta^F_{X_1,\dots,X_{n+1}}$ coincides with the kernel of $r_F$$$K[r_F]=\bigcap_{i=1}^{n+1}K[\pi^F_{\hat{X}_i}],$$which is precisely the kernel intersection of \cite{HL,HV} serving as their definition for the $(n+1)$st cross-effect $cr_{n+1}^F(X_1,\dots,X_{n+1})$ of the functor $F$.\vspace{1ex}

For $n=1$ and $F=id_\EE$ we get the following $2$-cube $\Xi_{X,Y}$
$$ \xymatrix@=15pt{
          X+Y \ar[rr]^{\pi_{\hat X}} \ar[d]_{\pi_{\hat Y}} && Y \ar[d]\\
          X \ar[rr]   && \star_\EE
                           }
     $$so that the limit $P_{X,Y}$ of the punctured $2$-cube is $X\times Y$ and the comparison map$$\theta_{X,Y}:X+Y\to X\times Y$$ is the one already used before. In particular, the just introduced notion of \emph{linear category} is the usual one.

%\emph{Note that we use the symbol $\pi_{\hat{X}}$ generically to denote a map defined on a sum which is the identity on all summands except $X$ on which it is the null morphism.}\vspace{1ex}

For $n=2$ and $F=id_\EE$ we get the following $3$-cube $\Xi_{X,Y,Z}$

$$ \xymatrix@=15pt{&Y+Z\ar[rr]^{\pi_{\hat{Y}}}\ar[dd]_>>>>>{\pi_{\hat{Z}}}&&Z\ar[dd]\\
          X+Y+Z \ar[dd]_>>>>>{\pi_{\hat{Z}}}\ar[rr]^>>>>>{\pi_{\hat{Y}}}\ar[ru]^{\pi_{\hat{X}}}&& X+Z \ar[dd]_>>>>>{\pi_{\hat{Z}}}\ar[ru]^{\pi_{\hat{X}}}&\\&Y\ar[rr]\ar@{.>}@<-1ex>[uu]&&\star_\EE\ar@<-1ex>@{.>}[uu]\\
          X+Y\ar[rr]^{\pi_{\hat{Y}}}\ar[ru]^{\pi_{\hat{X}}} \ar@<-1ex>@{.>}[uu]  && X\ar[ru]\ar@<-1ex>@{.>}[uu]&
                           }$$
which induces a split natural transformation of $2$-cubes: $$\Xi_{X,Y}+Z\rightleftarrows \Xi_{X,Y}$$ For sake of simplicity, we denote by $+Z$ the functor $\EE\to\Pt_Z(\EE)$ which takes an object $X$ to $X+Z\to Z$ with obvious section, and similarly, we denote by $\times Z$ the functor $\EE\to\Pt_Z(\EE)$, which takes an object $X$ to $X\times Z\to Z$ with obvious section.

The previous split natural transformation of 2-cubes induces a natural transformation of split epimorphisms
$$
\xymatrix@=15pt{
          X+Y+Z \ar@{}[rrd]|*+[o][F-]{2}\ar@{->}[rrr]^<<<<<<<<<<{\theta_{X,Y}^{+Z}}\ar[d]_{\pi_{\hat Z}} &&& (X+Z)\times_Z(Y+Z) \ar[d] \\
          X+Y\ar@{->}[rrr]_{\theta_{X,Y}}\ar@<-1ex>[u]  &&& X\times Y\ar@<-1ex>[u]
                           }
$$
the comparison map of which may be identified with the comparison map $$\theta_{X,Y,Z}:X+Y+Z\to P_{X,Y,Z}$$ of $\Xi_{X,Y,Z}$. In particular, the category $\EE$ is quadratic if and only if square (2) is a pullback square. Notice that  in a regular Mal'tsev category, the downward-oriented square (2) is necessarily a regular pushout by Corollary \ref{regpush}.

\begin{prp}\label{totalkernel}In a $\sg$-pointed regular Mal'tsev category, the comparison map $\theta_{X,Y,Z}$ is a regular epimorphism with kernel relation $R[\pi_{\hat X}]\cap R[\pi_{\hat Y}]\cap R[\pi_{\hat Z}]$.
\end{prp}
\proof
The first assertion expresses the regularity of pushout (2), the second follows from identities $R[\theta^{+Z}_{X,Y}]=R[\pi_{\hat X}]\cap R[\pi_{\hat Y}]$ and $R[\theta_{X,Y,Z}]=R[\theta^{+Z}_{X,Y}]\cap R[\pi_{\hat Z}]$ which hold because both, $\theta^{+Z}_{X,Y}$ and $\theta_{X,Y,Z}$, are comparison maps.\endproof

\begin{lma}\label{quad0}
A $\sg$-pointed category with pullbacks is quadratic if and only if square (2$\,'$) of the following diagram
$$ \xymatrix@=20pt{
      X+Y+Z \ar@{}[rrd]|*+[o][F-]{2'} \ar@{->}[rr]^{\theta_{X+Y,Z}} \ar[d]_{\theta_{X,Y}^{+ Z}} && (X+Y)\times Z\ar@{}[rrd]|*+[o][F-]{2''}\ar[d]^{\theta_{X,Y}^{\times Z}}\ar@{->}[rr] &&  X+Y\ar[d]^{\theta_{X,Y}}\\
   (X+Z)\times_Z(Y+Z) \ar@{->}[rr]_{\theta_{X,Z}\times_Z\theta_{Y,Z}}  && (X\times Z)\times_Z(Y\times Z)\ar@{->}[rr] &&  X\times Y
                       }
    $$
is a pullback square.
   \end{lma}

\proof Composing squares (2$'$) and (2$''$) yields square $(2)$ above. Square (2$''$) is a pullback since $(X\times Z)\times_Z(Y\times Z)$ is canonically isomorphic to $(X\times Y)\times Z$.\endproof

\subsection{The main diagram}\label{maindiagram}We shall now give several criteria for quadraticity. For this we consider the following diagram
$$ \xymatrix@=15pt{
     {(X+Y)\diamond Z\;} \ar@{>->}[rr]\ar[d]_{\theta_{X,Y}\diamond Z}  && (X+Y)+Z \ar@{}[rrd]|*+[o][F-]{a} \ar@{->>}[rr]^{\theta_{X+Y,Z}} \ar[d]^{\theta_{X,Y}+ Z} && (X+Y)\times Z\ar[d]^{\theta_{X,Y}\times Z} \\
     {(X\times Y)\diamond Z\;} \ar@{>->}[rr] \ar[d]_{\varphi^Z_{X,Y}} && (X\times Y)+Z \ar@{}[rrd]|*+[o][F-]{b} \ar@{->>}[rr]^{\theta_{X\times Y,Z}} \ar[d]^{\phi_{X,Y}^Z} && (X\times Y) \times Z \ar[d]^{\mu_{X,Y}^Z} \\
     {(X\diamond Z)\times (Y\diamond Z)\;} \ar@{>->}[rr] && (X+Z)\times_Z (Y+Z) \ar@{->>}[rr]_{\theta_{X,Z}\times_Z \theta_{Y,Z}} && (X\times Z)\times_Z (Y\times Z)
                       }
    $$
in which the vertical composite morphisms from left to right are $\theta_{X,Y}^{\diamond Z},\,\theta_{X,Y}^{+Z},\theta_{X,Y}^{\times Z}$, the horizontal morphisms on the left are the kernel-inclusions of the horizontal regular epimorphisms on their right, and $\mu_{X,Y}^Z$ is the canonical isomorphism.

Observe that square (b) is a pullback if and only if the canonical map $$\phi^Z_{X,Y}:(X\times Y)+Z\to (X+Z)\times_Z(Y+Z)$$ is invertible. This is the case if and only if for each $Z$ pointed cobase-change $$(\alpha_Z)_!:\DD\to\Pt_Z(\DD)$$along the initial map $\alpha_Z:\star_\DD\to Z$ preserves binary products, cf. Section \ref{pointfibration}.

\begin{prp}\label{quad1}
A $\sg$-pointed regular Mal'tsev category is quadratic if and only if squares $(a)$ and $(b)$ of the main diagram are pullback squares.\end{prp}
\proof
Since composing squares (a) and (b) yields square (2$'$) of Lemma \ref{quad0}, the condition is sufficient. Lemma \ref{propbr} and Corollary \ref{corbr} imply that the condition is necessary as well.\endproof
\begin{thm}\label{quad}
A $\sg$-pointed exact Mal'tsev category is quadratic if and only if it is $2$-nilpotent and pointed cobase-change along initial maps preserves binary products.\end{thm}
\proof
By Proposition \ref{quad1}, the category is quadratic if and only if the squares (a) and (b) are pullback squares, i.e. the category is $2$-nilpotent by Proposition \ref{ps}, and pointed cobase-change along initial maps preserves binary products.\endproof

\begin{cor}\label{bilinear}A semi-abelian category is quadratic\footnote{The authors of \cite{CGV,HL,HV} call a semi-abelian category \emph{two-nilpotent} if each object has a vanishing ternary Higgins commutator, cf. Remark \ref{nfoldedgroup}. By Proposition \ref{nadditive} this means that the identity functor is quadratic. Corollaries \ref{bilinear} and \ref{parexemple} describe how to enhance a $2$-nilpotent semi-abelian category in our sense so as to get a two-nilpotent semi-abelian category in their sense.} if and only if either of the following three equivalent conditions is satisfied:
\begin{itemize}\item[(a)]all objects are $2$-nilpotent, and the comparison maps $$\varphi_{X,Y}^Z:(X\times Y)\diamond Z\to(X\diamond Z)\times(Y\diamond Z)$$ are invertible for all objects $X,Y,Z$;\item[(b)]the third cross-effects of the identity functor $cr_3(X,Y,Z)=K[\theta_{X,Y,Z}]$ vanish for all objects $X,Y,Z;$\item[(c)]the co-smash product is linear, i.e. the canonical comparison maps$$\,\theta_{X,Y}^{\diamond Z}:(X+Y)\diamond Z\to (X\diamond Z)\times (Y\diamond Z)$$ are invertible for all objects $X,Y,Z.$\end{itemize}\end{cor}

\proof Theorem \ref{quad} shows that condition (a) amounts to quadraticity.

For condition (b) note that by protomodularity the cross-effect $K[\theta_{X,Y,Z}]$ vanishes if and only if the regular epimorphism $\theta_{X,Y,Z}$ is invertible.

The equivalence of conditions (b) and (c) follows from the isomorphism of kernels $K[\theta_{X,Y,Z}]\cong K[\theta_{X,Y}^{\diamond Z}]$. The latter is a consequence of the $3\times 3$-lemma which, applied to main diagram \ref{maindiagram} and square (2), yields the chain of isomorphisms $$K[\theta_{X,Y}^{\diamond Z}]\cong K[K[\theta^{+Z}_{X,Y}]\onto K[\theta^{\times Z}_{X,Y}]]\cong K[\theta_{X,Y,Z}].$$\endproof

\subsection{Algebraic distributivity and algebraic extensivity}\label{extensive}--\vspace{1ex}

We shall see that in a pseudo-additive regular Mal'tsev category $(\DD,\star_\DD)$, pointed cobase-change along initial maps $\alpha_Z:\star_\DD\to Z$ preserves binary products if and only if pointed base-change along terminal maps $\omega_Z:Z\to \star_\DD$ preserves binary sums. The latter condition means that for all objects $X,Y,Z$, the following square$$ \xymatrix@=20pt{
        Z\ar@{->}[d] \ar@{->}[r] & Y\times Z \ar@{->}[d]\\
        X\times Z \ar@{->}[r] & (X+Y)\times Z
                       }
    $$is a pushout, inducing thus for all objects $X,Y,Z$ an isomorphism $$(X\times Z)+_Z(Y\times Z)\cong(X+Y)\times Z$$ which can be considered as an algebraic distributivity law. This suggests the following definitions, where the added adjective ``algebraic'' means here that the familiar definition has to be modified by replacing the \emph{slices} of the category with the \emph{fibres} of the fibration of points, cf. Section \ref{pointfibration} and Carboni-Lack-Walters \cite{CLW}.

\begin{dfn}A category with pullbacks of split epimorphisms is \emph{algebraically distributive} if pointed base-change along terminal maps preserves binary sums.

A category with pullbacks of split epimorphisms and pushouts of split monomorphisms is \emph{algebraically extensive} if any pointed base-change preserves binary sums.\end{dfn}

We get the following implications between several in literature studied ``algebraic'' notions, where we assume that pullbacks of split epimorphisms and (whenever needed) pushouts of split monomorphisms exist:
$$ 
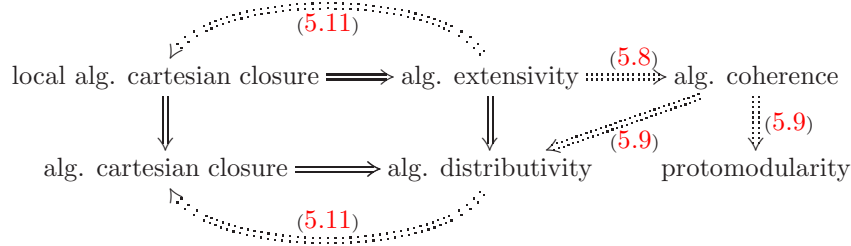
\xymatrix@=20pt{
        \textrm{local alg. cartesian closure}\ar@{=>}[d] \ar@{=>}[r] \ar@(u,u)@{<:}[r]_{(\ref{variety})}& \textrm{alg. extensivity}\ar@{=>}[d]\ar@{:>}[r]^>>>>>{(\ref{lad>ac})}&\textrm{alg. coherence}\ar@{:>}[d]^{(\ref{coherent})}\ar@{:>}[ld]^{(\ref{coherent})}\\
        \textrm{alg. cartesian closure}\ar@{=>}[r]\ar@{=>}[r] \ar@(d,d)@{<:}[r]^{(\ref{variety})} & \textrm{alg. distributivity}
                        &\textrm{protomodularity}}$$

The existence of centralisers implies \emph{algebraic cartesian closure} \cite{BGr} and hence \emph{algebraic distributivity}, cf. Section \ref{centralizer}. The categories of groups and of Lie algebras are not only algebraically cartesian closed, but also \emph{locally algebraically cartesian closed} \cite{Gr,Gr1}, which means that any pointed base-change admits a right adjoint. \emph{Algebraic coherence}, cf. Cigoli-Gray-Van der Linden \cite{CGV}, requires any pointed base-change to be \emph{coherent}, i.e. to preserve strongly epimorphic cospans.
%Algebraic coherence is closed under epireflections, cf. \cite[Proposition 2.19]{CGV}.

\begin{lma}\label{lad>ac}An algebraically extensive regular category is algebraically coherent.\end{lma}

\proof In a regular category, pointed base-change preserves regular epimorphisms. Henceforth, if the fibres have binary sums and pointed base-change preserves them, pointed base-change also preserves strongly epimorphic cospans.\endproof

\begin{lma}[cf. \cite{Bourn6}, Theorem 3.10 and \cite{CGV}, Theorem 6.1]\label{coherent}An algebraically coherent pointed Mal'tsev category is protomodular and algebraically distributive.\end{lma}

\proof To any split epimorphism $(r,s):Y\rightleftarrows X$ we associate the split epimorphism $(\bar{r}=r\times 1_X,\bar{s}=s\times 1_X):Y\times X\rightleftarrows X\times X$

 $$
\xymatrix@=30pt{
{Y\times X\;} \ar@<+1ex>[rr]^{\bar{r}}\ar@<-1ex>[rd]_{p_2} && X\times X \ar@{->}[ll]^{\bar{s}} \ar@<+1ex>[dl]^{p_2}\\
 & X \ar[lu]_{(s,1_X)} \ar[ur]^{(1_X,1_X)}
              }
$$ in the fibre over $X$. The kernel of $(\bar{r},\bar{s})$ may be identified with the given point $(r,s)$ over $X$ where the kernel-inclusion is defined by $(1_Y,r):Y\into Y\times X$. Kernel-inclusion and section strongly generate the point $Y\times X$ over $X$, cf. \cite[Proposition 3.7]{Bourn6}. Pointed base-change along $\alpha_X:\star\to X$ takes $(\bar{r},\bar{s})$ back to $(r,s)$, so that by algebraic coherence, section and kernel-inclusion of $(r,s)$ strongly generate $Y$. In a pointed category this amounts to protomodularity.

For the second assertion observe that if $F$ and $G$ are composable coherent functors such that $G$ is conservative and $GF$ preserves binary sums, then $F$ preserves binary sums as well; indeed, the isomorphism $GF(X)+GF(Y)\to GF(X+Y)$ decomposes into two isomorphisms $GF(X)+GF(Y)\to G(F(X)+F(Y))\to GF(X+Y)$. This applies to $F=\omega_Z^*$ and $G=\alpha_Z^*$ (where $\alpha_Z:\star\to Z$ and $\omega_Z:Z\to\star$) because $\omega_Z\alpha_Z=id_Z$ and $\alpha_Z^*$ is conservative, so that $\omega_Z^*$ preserves binary sums for all $Z$.\endproof

\begin{lma}\label{semiabelian}A $\sg$-pointed exact Mal'tsev category is algebraically extensive if and only if it is a semi-abelian category with exact pointed base-change functors.\end{lma}

%Every algebraically extensive semi-abelian category is peri-abelian (cf. \cite{BGr}, 4.4).\end{lma}

\proof This follows from Lemmas \ref{lad>ac} and \ref{coherent} and the fact that a left exact and regular epimorphism preserving functor between semi-abelian categories is right exact if and only if it preserves binary sums, cf. Section \ref{pointedexactMalcev}.\endproof

%The second assertion follows from the first in virtue of Lemma \ref{Huq2} and Proposition \ref{abelianization}.\endproof

For the following lemma a \emph{variety} means a category equipped with a forgetful functor to sets which is monadic and preserves filtered colimits. Every variety is bicomplete, cowellpowered, and has finite limits commuting with filtered colimits.

\begin{lma}[cf. \cite{Gr}, Theorem 2.9]\label{variety}A semi-abelian variety is (locally) algebraically cartesian closed if and only if it is algebraically distributive (extensive).\end{lma}

\proof Since the fibres of a semi-abelian category are semi-abelian, the pointed base-change functors preserve binary sums if and only if they preserve finite colimits, cf. Section \ref{pointedexactMalcev}. Since any colimit is a filtered colimit of finite colimits, and pointed base-change functors of a variety preserve filtered colimits, they preserve binary sums if and only if they preserve all colimits. It follows then from Freyd's special adjoint functor theorem that a pointed base-change functor of a semi-abelian variety preserves binary sums if and only if it has a right adjoint.\endproof

A pointed category $\DD$ with pullbacks is called \emph{fibrewise} algebraically cartesian closed (resp. distributive) if for all objects $Z$ of $\DD$ the fibres $\Pt_Z(\DD)$ are algebraically cartesian closed (resp. distributive). This is the case if and only if pointed base-change along every \emph{split epimorphism} has a right adjoint (resp. preserves binary sums). Any algebraically coherent pointed Mal'tsev category is fibrewise algebraically distributive, cf. the proof of Lemma \ref{coherent}.

\begin{prp}\label{accclosure}For a pointed regular Mal'tsev category, fibrewise algebraic cartesian closure is preserved under strong epireflections.\end{prp}

\proof
Let $(r,s):X\rightleftarrows Y$ be a point in a strongly epireflective subcategory $\CC$ of a fibrewise algebraically cartesian closed regular Mal'tsev category $\DD$. Let $f:Y\to Z$ be a split epimorphism in $\CC$. We shall denote $f_*:\Pt_Y(\DD)\to\Pt_Z(\DD)$ the \emph{right} adjoint of pointed base-change $f^*:\Pt_Z(\DD)\to\Pt_Y(\DD)$. Consider the following diagram
$$ \xymatrix@=20pt{
X\ar@<-1ex>[dd]_<<<<<{r}\ar@{=}[rd] && \bar{\bar X}\ar[ll]_{\varepsilon} \ar@{->>}[rd]^{\eta_{\bar{\bar X}}} \ar@<-1ex>[dd]_<<<<<{\bar{\bar r}} \ar[rr]^{\bar{f}}  && \bar X \ar@<1ex>[ll]^{} \ar@<-1ex>[dd]_<<<<<{\bar r} \ar@{->>}[rd]^{\eta_{\bar X}}\\
& X\ar[dl]^{r} && I(\bar{\bar X}) \ar[ll]_>>>>{I(\varepsilon)} \ar[dl]^<<<<{I(\bar{\bar r})} \ar[rr]^<<<{I(\bar{f})}  && I(\bar X) \ar[dl]^{I(\bar r)}\ar@(u,u)@{.>}[ul]_{\bar{\phi}}\\
Y \ar@{=}[rr] \ar[uu]_>>>>>{s} && Y \ar[rr]^{f}\ar[uu]_>>>>>{\bar{\bar s}} && Z \ar@<1ex>[ll]^{} \ar[uu]_>>>>>{\bar s}
}
$$
where $(\bar r,\bar s)=f_*(r,s)$, the downward-oriented right square is a pullback, and $\varepsilon:(\bar{\bar{r}},\bar{\bar{s}})=f^*f_*(r,s)\to (r,s)$ is the counit at $(r,s)$. Since $\DD$ is a regular Mal'tsev category, the strong epireflection $I$ preserves this pullback of split epimorphisms (cf. Proposition \ref{Diana}) so that $I(\bar{\bar{X}})$ is isomorphic to $f^*(I(\bar{X}))$. Since by adjunction, maps $I(\bar{X})\to f_*(X)=\bar{X}$ correspond bijectively to maps $I(\bar{\bar{X}})=f^*(I(\bar{X}))\to X$ there is a unique dotted map $\bar{\phi}: I(\bar X)\to \bar X$ such that $\varepsilon\circ f^*(\bar{\phi})=I(\varepsilon)$.

Accordingly we get $\bar{\phi}\eta_{\bar X}=1_{\bar X}$ so that $\eta_{\bar X}$ is invertible and hence $\bar{X}$ belongs to $\CC$. This shows that the right adjoint $f_*:\Pt_Y(\DD)\to\Pt_Z(\DD)$ restricts to a right adjoint $f_*:\Pt_Y(\CC)\to\Pt_{Z}(\CC)$ so that $\CC$ is fibrewise algebraically cartesian closed.\endproof

For regular Mal'tsev categories, algebraic cartesian closure amounts to the existence of centralisers for all (split) subobjects, see Section \ref{centralizer}. Part of Proposition \ref{accclosure} could thus be reformulated by saying that in this context the existence of centralisers is preserved under strong epireflections, which can also be proved directly. In a varietal context, Proposition \ref{accclosure} also follows from Lemmas \ref{variety} and \ref{algdis}.

\begin{lma}\label{fibreextensive}If an algebraically extensive semi-abelian (or homological \cite{BB}) category $\DD$ has an identity functor of degree $\leq n$, then all its fibres $\Pt_Z(\DD)$ as well.\end{lma}

\proof The kernel functors $(\alpha_Z)^*:\Pt_Z(\DD)\to \DD$ are conservative and preserve binary sums. Therefore, the kernel functors preserve the limits $P^{\Pt_Z(\DD)}_{X_1,\dots,X_{n+1}}$ and the comparison maps $\theta^{\Pt_Z(\DD)}_{X_1,\dots,X_{n+1}}$. Accordingly, if the identity functor of $\DD$ is of degree $\leq n$, then $\alpha_Z^*(\theta_{X_1,\dots,X_{n+1}})$ is invertible, hence so is $\theta^{\Pt_Z(\DD)}_{X_1,\dots,X_{n+1}}$, for all objects $X_1,\dots,X_{n+1}$ of the fibre $\Pt_Z(\DD)$.\endproof

It should be noted that in general neither algebraic extensivity nor local algebraic cartesian closure is preserved under Birkhoff reflection. This is in neat contrast to (fibrewise) algebraic distributivity and algebraic coherence, which are preserved under strong epireflections, cf. Proposition \ref{accclosure} and \cite[Proposition 3.7]{CGV}.

\subsection{Duality for pseudo-additive regular Mal'tsev categories}
\begin{lma}\label{duality}For any pointed category $\DD$ with binary sums and binary products consider the following commutative diagram $$ \xymatrix@=10pt{
       (X\times Y)+ Z \ar@{.>}[dd]^{p_X^Y+Z} \ar[rr]^{\rho_{X,Y,Z}} && X\times (Y+Z) \ar@{.>}[dd]^{X\times \pi_Z^Y}\\
        && \\
       X+Z\ar[rr]_{\theta_{X,Z}} \ar@<1ex>[uu]^{j_X^Y+Z} && X\times Z \ar@<1ex>[uu]^{X\times \iota_Z^Y}
                           }
     $$
     in which $\rho_{X,Y,Z}$ is induced by the pair $X\times \iota_Y^Z:X\times Y \to X\times (Y+Z)$ and $\alpha_X\times \iota_Z^Y:Z\to X\times (Y+Z)$.
     \begin{itemize}\item[(1)]pointed base-change $(\omega_X)^*:\DD \to\Pt_X(\DD)$ preserves binary sums if and only if the upward-oriented square is a pushout for all objects $Y,Z$;
\item[(2)]pointed cobase-change $(\alpha_Z)_!:\Pt_Z(\DD)\to\DD$ preserves binary products if and only if the downward-oriented square is a pullback for all objects $X,Y$.\end{itemize}
\end{lma}\label{tau}
\proof
The left upward-oriented square of the following diagram
$$\xymatrix@=10pt{
   X\times Y \ar@{.>}[dd]^{p_X^Y}\ar[rr]_{\iota_{X\times Y}^Z}\ar@<2ex>[rrrr]^{X\times \iota_{Y}^Z}  && (X\times Y)+ Z \ar@{.>}[dd]^{p_X^Y+Z} \ar[rr]_{\rho_{X,Y,Z}} && X\times (Y+Z) \ar@{.>}[dd]^{X\times \pi_Z^Y}\\
    &&    && \\
   X \ar@<1ex>[uu]^{j_X^Y}\ar[rr]^{\iota_{X}^Z} \ar@<-2ex>[rrrr]_{j_X^Z}  &&  X+Z\ar[rr]^{\theta_{X,Z}} \ar@<1ex>[uu]^{j_X^Y+Z} && X\times Z \ar@<1ex>[uu]^{X\times \iota_Z^Y}
                           }
     $$is a pushout so that the whole upward-oriented rectangle is a pushout (i.e. $(\omega_X)^*$ preserves binary sums) if and only if the right upward-oriented square is a pushout.

The right downward-oriented square of the following diagram
     $$\xymatrix@=10pt{
        (X\times Y)+ Z \ar@{.>}[dd]^{p_X^Y+Z} \ar[rr]_{\rho_{X,Y,Z}}\ar[rr]_{\rho_{X,Y,Z}} \ar@<2ex>[rrrr]^{p_Y^X+Z} && X\times (Y+Z) \ar[rr]_{p_{Y+Z}^X} \ar@{.>}[dd]^{X\times \pi_Z^Y} && Y+Z\ar@{.>}[dd]^{\pi_Z^Y}\\
         &&    && \\
       X+Z\ar[rr]^{\theta_{X,Z}} \ar@<1ex>[uu]^{j_X^Y+Z}\ar@<-2ex>[rrrr]_{\pi_Z^X} && X\times Z \ar@<1ex>[uu]^{X\times \iota_Z^Y} \ar[rr]^{p_Z^X} && Z \ar@<1ex>[uu]^{\iota_Z^Y}
                                }
          $$is a pullback so that the whole downward-oriented rectangle is a pullback (i.e. $(\alpha_Z)_!$ preserves binary products) if and only if the left downward-oriented square is a pullback.
\endproof

\begin{prp}\label{tau=iota}In a $\sg$-pointed regular Mal'tsev category $\DD$, pointed base-change $(\omega_Z)^*:\DD \to\Pt_Z(\DD)$ preserves binary sums for all objects $Z$ as soon as pointed cobase-change $(\alpha_Z)_!:\DD\to\Pt_Z(\DD)$ preserves binary products for all objects $Z$. The converse implication holds if $\,\DD$ is pseudo-additive (cf. Definition \ref{pseudoadditive}).
\end{prp}
\proof
According to the previous lemma pointed cobase-change $(\alpha_X)_!$ preserves binary products if and only if the downward-oriented square is a pullback which implies that the upward-oriented square is a pushout, and hence pointed base-change $(\omega_Z)^*$ preserves binary sums. If $\DD$ is pseudo-additive, the comparison map $\theta_{X,Z}$ is an affine extension. Therefore, the downward-oriented square is a pullback if and only if the upward-oriented square is a pushout, whence the converse.
\endproof

\begin{cor}\label{parexemple}
A $\sg$-pointed exact Mal'tsev category is quadratic if and only if it is $2$-nilpotent and algebraically distributive.
\end{cor}

\begin{lma}\label{algdis}In a $\sg$-pointed regular Mal'tsev category, (fibrewise) algebraic distributivity is preserved under strong epireflections.\end{lma}

\proof This follows from Proposition \ref{Diana} which shows that strong epireflections preserve besides pushouts and binary sums also binary products (in the fibres).\endproof

\begin{thm}\label{acc}
For any algebraically distributive, $\sg$-pointed exact Mal'tsev category, the Birkhoff subcategory of $2$-nilpotent objects is quadratic.\end{thm}
\proof The Birkhoff subcategory is pointed, exact, $2$-nilpotent, and algebraically distributive by Lemma \ref{algdis}, and hence quadratic by Corollary \ref{parexemple}.\endproof

\begin{cor}[cf. Cigoli-Gray-Van der Linden \cite{CGV}, Corollary 7.2]\label{Huq=Higgins}For each object $X$ of an algebraically distributive semi-abelian category, the iterated Huq commutator $[X,[X,X]]$ coincides with the ternary Higgins commutator $[X,X,X]$.\end{cor}
\proof
Recall (cf. \cite{HL,HV}) that $[X,X,X]$ is the direct image of $K[\theta_{X,X,X}]$ under the ternary folding map $X+X+X\to X$. In general, the iterated Huq commutator $[X,[X,X]]$ is contained in $[X,X,X]$, cf. Corollary \ref{inclusion}. In a semi-abelian category, the unit of second Birkhoff reflection $I^2$ takes the form $\eta^2_X:X\to X/[X,[X,X]]$, cf. Remark \ref{Huq}. Since in the algebraically distributive case, the subcategory of $2$-nilpotent objects is quadratic by Theorem \ref{acc}, the image of $[X,X,X]$ in $X/[X,[X,X]]$ is trivial by Corollaries \ref{bilinear}b and \ref{Higgins}, whence $[X,[X,X]]=[X,X,X]$.\endproof

\begin{rmk}The category of groups (resp. Lie algebras) has centralisers for subobjects and is thus algebraically distributive. Therefore, the category of $2$-nilpotent groups (resp. Lie algebras) is a quadratic semi-abelian variety.

The reader should observe that although on the level of $2$-nilpotent objects there is a perfect symmetry between the property that pointed base-change along terminal maps preserves binary sums and the property that pointed cobase-change along initial maps preserves binary products (cf. Proposition \ref{tau=iota}), only the algebraic distributivity carries over to the category of all groups (resp. Lie algebras) while the algebraic ``codistributivity'' fails in either of these categories. ``Codistributivity'' is a quite restrictive property, which is rarely satisfied without assuming $2$-nilpotency.\end{rmk}

\section{Identity functors with bounded degree}\label{degreen}

In the previous section we have seen that quadraticity is a slightly stronger property than $2$-nilpotency, insofar as it also requires a certain compatibility between binary sum and binary product (cf. Theorem \ref{quad} and Proposition \ref{tau=iota}). In this last section, we relate $n$-nilpotency to identity functors of degree $\leq n$.

\subsection{Degree $n$ functors and $n$-folded objects}--\vspace{1ex}

Any $(n+1)$-cube $\Xi_{X_1,\dots,X_{n+1}}$ (cf. Definition \ref{ncube}) defines a split natural transformation of $n$-cubes inducing a natural transformation of split epimorphisms

$$
\xymatrix@=20pt{
          X_1+\cdots+X_{n+1}\ar@{}[rrd]|*+[o][F-]{n}\ar@{->}[rrrr]^{\theta_{X_1,\dots,X_n}^{+X_{n+1}}}\ar[d]_{\pi_{\hat{X}_{n+1}}}&& && P^{+X_{n+1}}_{X_1,\dots,X_n} \ar[d]_{P^{+\omega_{X_{n+1}}}_{X_1,\dots,X_n}} \\
          X_1+\cdots+X_{n }\ar@{->}[rrrr]_{\theta_{X_1,\dots,X_n}}\ar@<-1ex>[u] && && P_{X_1,\dots,X_n} \ar@<-1ex>[u]_{P^{+\alpha_{X_{n+1}}}_{X_1,\dots,X_n}}
                           }
$$
the comparison map of which may be identified with the comparison map $$\theta_{X_1,\dots,X_{n+1}}:X_1+\cdots+X_{n+1}\to P_{X_1,\dots,X_{n+1}}$$ of the given $(n+1)$-cube. In particular, our category has an identity functor of degree $\leq n$ if and only if square (n) is a pullback square for all objects $X_1,\dots,X_{n+1}$.

\begin{prp}\label{intersection}In a $\sg$-pointed regular Mal'tsev category, the comparison map $\theta_{X_1,\dots,X_{n+1}}$ is a regular epimorphism with kernel relation $R[\pi_{\hat X_1}]\cap\cdots\cap R[\pi_{\hat{X}_{n+1}}]$.\end{prp}

\proof This follows by induction on $n$ like in the proof of Proposition \ref{totalkernel}.\endproof

\begin{rmk}\label{history}The previous proposition shows that in a $\sg$-pointed regular Mal'tsev category, the intersection of the kernel relations of the contraction maps may be considered as the ``total kernel relation'' of the cube. This parallels the more elementary fact that the total-kernel definition of the cross-effects $cr_{n+1}(X,\dots,X_{n+1})$ coincides with the kernel-intersection definition of Hartl-Loiseau \cite{HL} and Hartl-Van der Linden \cite{HV}. In particular, in any $\sg$-pointed regular Mal'tsev category, the image of the morphism $r_{id}:X_1+\cdots+X_{n+1}\to\prod_{i=1}^{n+1}X_1+\cdots+\widehat{X_i}+\cdots+X_{n+1}$ coincides with the limit $P_{X_1,\dots,X_{n+1}}$ of the punctured $(n+1)$-cube.

We already mentioned that these kernel intersections are dual to the $(n+1)$-fold smash products of Carboni-Janelidze \cite{CJ}. An alternative way to describe the duality between cross-effects and smash-products is to consider the limit construction $P_{X_1,\dots,X_{n}}$ as the dual of the so-called \emph{fat wedge} $T^{X_1,\dots,X_{n}}$, cf. Hovey \cite{Ho}. Set-theoretically, the fat wedge is the subobject of the product $X_1\times\cdots\times X_n$ formed by the $n$-tuples having at least one coordinate at a base-point. If base-point inclusions behave ``well'' with respect to cartesian product, the fat wedge is given by a colimit construction, strictly dual to the limit construction defining $P_{X_1,\dots,X_n}$. The $n$-fold smash-product $X_1\wedge\cdots\wedge X_n$ is then the cokernel of the monomorphism $T^{X_1,\dots,X_n}\into X_1\times\cdots\times X_n$ while the $n$-th cross-effect $cr_n(X_1,\dots,X_n)$ is the kernel of the regular epimorphism $X_1+\cdots+X_n\onto P_{X_1,\dots,X_n}$.

The cubical cross-effects are just the algebraic version of Goodwillie's homotopical cross-effects \cite[pg. 676]{G2}. Nevertheless, for functors taking values in abelian categories, the cubical cross-effects agree with the original cross-effects of Eilenberg-Mac Lane \cite[pg. 77]{EM}. Indeed, by \cite[Theorems 9.1 and 9.6]{EM}, for a based functor $F:\DD\to\EE$ from a $\sg$-pointed category $(\DD,+,\star_\DD)$ to an abelian category $(\EE,\oplus,0_\EE),$ the latter are completely determined by the following \emph{decomposition formula}$$F(X_1+\cdots+ X_n)\cong \bigoplus_{1\leq k\leq n}\bigoplus_{1\leq i_1<\cdots<i_k\leq n}cr_k^F(X_{i_1},\dots,X_{i_k})$$for any objects $X_1,\dots,X_n$ in $\DD$. It suffices thus to show that the cubical cross-effects satisfy the decomposition formula if values are taken in an abelian category.

For $n=2$ we get $P^F_{X_1,X_2}=F(X_1)\oplus F(X_2)$ from which it follows that $\theta^F_{X_1,X_2}:F(X_1+ X_2)\onto F(X_1)\oplus F(X_2)$ is a \emph{split} epimorphism. Henceforth, we get the asserted isomorphism $F(X_1+ X_2)\cong F(X_1)\oplus F(X_2)\oplus cr_2^F(X_1,X_2)$.

The $3$-cube $F(\Xi_{X_1,X_2,X_3})$ induces a natural transformation of split epimorphisms

$$
\xymatrix@=20pt{
          F(X_1+X_2+X_3)\ar@{}[rrd]\ar@{->}[rrrr]\ar[d]&& && P_3 \ar[d] \\
          F(X_1+X_2)\ar@{->>}[rrrr]_{\theta^F_{X_1,X_2}}\ar@<-1ex>[u] && && F(X_1)\oplus F(X_2) \ar@<-1ex>[u]
                           }
$$in which $P_3$ is isomorphic to $F(X_1)\oplus F(X_2)\oplus F(X_3)\oplus cr_2^F(X_1,X_3)\oplus cr_2^F(X_2,X_3)$. From this, we get for $P^F_{X_1,X_2,X_3}=F(X_1+ X_2)\times_{F(X_1)\oplus F(X_2)} P_3$ the formula $$P^F_{X_1,X_2,X_3}\cong F(X_1)\oplus F(X_2)\oplus F(X_3)\oplus cr_2^F(X_1,X_2)\oplus cr_2^F(X_1,X_3)\oplus cr_2^F(X_2,X_3)$$so that $\theta^F_{X_1,X_2,X_3}$ is again a \emph{split} epimorphism inducing the asserted decomposition of $F(X_1+X_2+X_3)$. The same scheme keeps on for all positive integers $n$.\qed

\end{rmk}

\begin{dfn}\label{nfolded}
An object $X$ of a $\sg$-pointed regular Mal'tsev category will be called \emph{$n$-folded} if the folding map $\delta_{n+1}^X$ factors through the comparison map $\theta_{X,\dots,X}$
$$ \xymatrix@=10pt{
     X+\cdots+X  \ar@{->>}[rr]^{\theta_{X,\dots,X}} \ar@{->>}[rdd]_{\delta^X_{n+1}} && P_{X,\dots,X} \ar@{.>}[ddl]^{m_X}  \\
      && \\
     & X
                            }
  $$i.e. if the folding map $\,\delta^X_{n+1}$ annihilates the kernel relation $R[\theta_{X,\dots,X}]$.\end{dfn}

An object $X$ is \emph{$1$-folded} if and only if the identity of $X$ commutes with itself. In a $\sg$-pointed regular Mal'tsev category this is the case if and only if $X$ is an abelian group object, cf. the proof of Proposition \ref{abelianization}.

\begin{rmk}\label{nfoldedgroup}In a semi-abelian (or homological \cite{BB}) category, an object $X$ is $n$-folded if and only if the image of the kernel $K[\theta_{X,\dots,X}]$ under the folding map $\delta_{n+1}^X:X+\cdots+X\onto X$ is trivial. Recall \cite{HL,HV,Hi,MM} that this image is by definition the so-called \emph{Higgins commutator} $[X,\dots,X]$ of length $n+1$. Therefore, an object of a semi-abelian category is $n$-folded precisely when its Higgins commutator of length $n+1$ vanishes. Under this form $n$-folded objects have already been studied by Hartl emphasising their role in his theory of polynomial approximation.

%We will see below that $n$-folded objects are closely related to functors of degree $\leq n$. In a semi-abelian category, the Higgins commutators share a similar relationship with degree $n$ functors, as the Huq commutators do with $n$-nilpotency. This was starting point of Hartl's ...

In a varietal context, $n$-foldedness can be expressed in more combinatorial terms. For instance, a group $X$ is \emph{$n$-folded} if and only if \emph{$(n+1)$-reducible} elements of $X$ are trivial. An element $w\in X$ is called $(n+1)$-reducible if there is an element $v$ in the free group $\Ff(X\sqcup\cdots\sqcup X)$ on $n+1$ copies of $X$ (viewed as a set) such that\begin{itemize}\item[(a)]$w$ is the image of $v$ under the composite map$$\Ff(\overbrace{X\sqcup\cdots\sqcup X}^{n+1})\cong \overbrace{\Ff(X)+\cdots+\Ff(X)}^{n+1}\onto \overbrace{X+\cdots+X}^{n+1}\overset{\delta_{n+1}^X}{\onto}X$$
\item[(b)]for each of the $n+1$ contractions $\pi_i^{\Ff(X)}:\Ff(X)^{+(n+1)}\onto\Ff(X)^{+n}$, cf. Section \ref{Higginscommutator}, the image $\pi_i^{\Ff(X)}(v)$ maps to the neutral element of $X$ under$$\overbrace{\Ff(X)+\cdots+\Ff(X)}^n\onto\overbrace{X+\cdots+X}^n\overset{\delta_{n}^X}{\onto}X.$$\end{itemize}
Indeed, since the evaluation map $\Ff(X)\onto X$ is a regular epimorphism, and the evaluation maps in (a) and (b) are compatible with the contraction maps $\pi_i:X^{+ (n+1)}\onto X^{+ n}$, Proposition \ref{intersection}, Section \ref{Higginscommutator} and Corollary \ref{Higgins} imply that we get in this way the image of the kernel $K[\theta_{X,\dots,X}]$ under the folding map $\delta^X_{n+1}$.

Any product of commutators $\prod_{i=1}^k[x_i,y_i]=\prod_{i=1}^kx_iy_ix_i^{-1}y_i^{-1}$ in $X$ is $2$-reducible by letting the $x_i$ (resp. $y_i$) belong to the first (resp. second) copy of $X$. Conversely, a direct computation shows that any $2$-reducible element of $X$ can be rewritten as a product of commutators of $X$. This recovers in a combinatorial way the aforementioned fact that $X$ is abelian (i.e. $1$-nilpotent) if and only if $X$ is $1$-folded.

The relationship between $n$-nilpotency and $n$-foldedness is more subtle, closely related to the cross-effects of the identity functor (cf. Theorem \ref{additivereflection}). For groups and Lie algebras the two concepts coincide (cf. Theorem \ref{folklore}c) but, though any $n$-folded object is $n$-nilpotent, the converse is wrong in general (cf. Section \ref{Moufang}).\end{rmk}

\begin{prp}\label{nadditive}Let $F:\DD\to\EE$ be a based functor between $\sg$-pointed categories and assume that $\EE$ is a regular Mal'tsev category.\begin{itemize}\item[(a)]If $F$ is of degree $\leq n$ then $F$ takes values in $n$-folded objects of $\,\EE$;\item[(b)]If $F$ preserves binary sums and takes values in $n$-folded objects of $\,\EE$ then $F$ is of degree $\leq n$;\item[(c)]The identity functor of $\,\EE$ is of degree $\leq n$ if and only if all objects of $\,\EE$ are $n$-folded.\end{itemize}\end{prp}
\proof Clearly, (c) follows from (a) and (b). For (a) note that $\delta^{F(X)}_{n+1}$ factors through $F(\delta^X_{n+1})$, and that by definition of a functor of degree $\leq n$, the comparison map $\theta^F_{X,\dots,X}$ is invertible so that $F(\delta^X_{n+1})$ gets identified with $m_{F(X)}$.

For (b) observe first that preservation of binary sums yields the isomorphisms $P^F_{X_1,\dots,X_{n+1}}\cong P_{F(X_1),\dots,F(X_{n+1})}$ and $\theta^F_{X_1,\dots,X_{n+1}}\cong\theta_{F(X_1),\dots,F(X_{n+1})}$. We shall show that if moreover $F$ takes values in $n$-folded objects of $\EE$ then $\theta^F_{X_1,\dots,X_{n+1}}$ is invertible for all $(n+1)$-tuples $(X_1,\dots,X_{n+1})$ of objects of $\DD$.

Consider any family $(f_i:X_i\to T)_{1\leq i\leq n+1}$ of morphisms in $\EE$, and let $\phi=\delta^T_{n+1}\circ(f_1+\cdots+f_{n+1}):X_1+\cdots+X_{n+1}\to T$ be the induced map.  We have the following factorisation of $F(\phi)$ through $\theta^F_{X_1,\dots,X_{n+1}}$:
$$ \xymatrix@=10pt{
      F(X_1)+\cdots+F(X_{n+1}) \ar@{->>}[rr]^{\theta_{F(X_1),\dots,F(X_{n+1})}} \ar[dd]_{F(f_1)+\cdots+F(f_{n+1})} && P_{F(X_1),\dots,F(X_{n+1})} \ar[dd]^{P_{F(f_1),\dots,F(f_{n+1})}} \\
      && \\
     F(T)+\cdots+F(T)  \ar@{->>}[rr]^{\theta_{F(T),\dots,F(T)}} \ar[rdd]_{\delta^T_{n+1}} && P_{F(T),\dots,F(T)} \ar[ddl]^{m_{F(T)}}  \\
      && \\
     & F(T)
                            }
  $$
In particular, if $T=X_1+\cdots+X_{n+1}$ and $f_i$ is the inclusion of the $i$th summand, we get a retraction of $\theta_{F(X_1),\dots,F(X_{n+1})}$ which accordingly is a monomorphism. Since $\theta_{F(X_1),\dots,F(X_{n+1})}$ is also a regular epimorphism, it is invertible.\endproof

\begin{prp}\label{nfoldBirkhoff}
The full subcategory $\Fld^n(\EE)$ of $n$-folded objects of a $\sg$-pointed regular Mal'tsev category $\,\EE$ is closed under products, subobjects and quotients.
\end{prp}
\proof
For any two $n$-folded objects $X$ and $Y$ the following diagram
$$ \xymatrix@=10pt{
      (X\times Y)+\cdots+(X\times Y) \ar@{->>}[rr]^>>>>>>>>>>{\theta_{X\times Y,\dots,X\times Y}} \ar[dd] && P_{X\times Y,\dots,X\times Y} \ar[dd] \\
      && \\
     (X+\cdots+X)\times (Y+\cdots+Y)  \ar@{->>}[rr]^>>>>>>>>>>{\theta_{X,\dots,X}\times \theta_{Y,\dots,Y}} \ar[rdd]_{\delta^X_{n+1}\times \delta^Y_{n+1}} && P_{X,\dots,X} \times P_{Y,\dots,Y} \ar[ddl]^{m_X\times m_Y}  \\
      && \\
     & X\times Y
                            }
  $$induces the required factorisation of $\delta_{n+1}^{X\times Y}$ through $\theta_{X\times Y,\dots,X\times Y}$.

For a subobject $n:U\rightarrowtail X$ of an $n$-folded object $X$ consider the diagram
 $$ \xymatrix@=10pt{
       U+\cdots+U \ar@{->>}[rrrrd]^{\theta_{U,\cdots,U}} \ar@{.>}[rd]^>>>>{\theta} \ar[dddd]_{\delta^U_{n+1}} \ar@(d,r)[rrddd]_{n+\cdots+n}\\
        & {W\;} \ar@{>.>}[rrr]^{\nu} \ar@{.>}[dddl] &&& P_{U,\dots,U} \ar[dd]^{P_{n,\dots,n}} \\
       &&&&&& \\
      &&  X+\cdots+X  \ar@{->>}[rr]^{\theta_{X,\dots,X }} \ar[d]_{\delta^X_{n+1}} && P_{X,\dots,X} \ar[dll]^{m_X}  \\
     {U\;} \ar@{>->}[rr]_n && X
                             }
   $$in which the dotted quadrangle is a pullback. The commutatitvity of the diagram induces a morphism $\theta$ such that $\nu\theta= \theta_{U,\dots,U}$. Since $\theta_{U,\dots,U}$ is a regular epimorphism, the monomorphism $\nu$ is invertible, whence the desired factorisation.

Finally, for a regular epimorphism $f:X\onto Y$ with $n$-folded domain $X$ consider the following diagram
$$ \xymatrix@=20pt{
  X+\cdots+X \ar@{->>}[rd]^{\theta_{X,\dots,X}} \ar[dd]_{\delta^X_{n+1}} \ar@{->>}[rr]^{f+\cdots+f}  && Y+\cdots+Y \ar@{->>}[rd]^{\theta_{Y,\dots,Y}}\ar[dd]^<<<<<<{\delta^Y_{n+1}}\\
     & P_{X,\dots,X}  \ar[dl]^{m_X} \ar@{->>}[rr]^<<<<<<{P_{f,\dots,f}}  && P_{Y,\dots,Y} \ar@{.>}[dl]^{m_Y}\\
  X \ar@{->>}[rr]_{f} && Y
 }
$$in which the existence of the dotted arrow has to be shown. According to Lemma \ref{regulartheta} the induced morphism on kernel relations $$R(f+\cdots+f,P_{f,\dots,f}):R[\theta_{X,\dots,X}]\to R[\theta_{Y,\dots,Y}]$$  is a regular epimorphism. A diagram chase shows then that $\delta^Y_{n+1}$ annihilates $R[\theta_{Y,\dots,Y}]$ whence the required factorisation of $\delta^Y_{n+1}$ through $\theta_{Y,\dots,Y}$.\endproof

\begin{lma}\label{regulartheta}In a $\sg$-pointed regular Mal'tsev category, any finite family of regular epimorphisms $f_i:X_i\onto Y_i\quad(i=1,\dots,n)$ induces a regular epimorphism on kernel relations $R(f_1+\cdots+f_n,P_{f_1,\dots,f_n}):R[\theta_{X_1,\dots,X_n}]\onto R[\theta_{Y_1,\dots,Y_n}]$.\end{lma}

\proof Since regular epimorphisms compose (in any regular category) it suffices to establish the assertion under the assumption $f_i=1_{X_i}$ for $i=2,\dots,n$. Moreover we can argue by induction on $n$ since for $n=1$ the comparison map is the terminal map $\theta_X:X\onto\star$ and a binary product of regular epimorphisms is a regular epimorphism. Assume now that the statement is proved for $n-1$ morphisms. Using the isomorphism of kernel relations $$R[\theta_{X_1,\dots,X_{n-1},X_n}]\cong R[R[\theta^{+X_n}_{X_1,\dots,X_{n-1}}]\onto R[\theta_{X_1,\dots,X_{n-1}}]]$$and Propositon \ref{regpushout1} it suffices then to show that the following by $f_1:X_1\onto Y_1$ induced commutative square$$ \xymatrix@=25pt{R[\theta^{+X_n}_{X_1,X_2,\dots,X_{n-1}}] \ar@<-1,ex>@{->>}[d] \ar@{->>}[r]& R[\theta^{+X_n}_{Y_1,X_2,\dots,X_{n-1}}]\ar@<-1,ex>@{->>}[d]\\
    R[\theta_{X_1,X_2,\dots,X_{n-1}}] \ar@{->>}[r]\ar@{->}[u] & R[\theta_{Y_1,X_2,\dots,X_{n-1}}]\ar@{->}[u]
                       }
    $$is a downward-oriented \emph{regular} pushout. This in turn follows from Corollary \ref{regpush} since the vertical arrows above are split epimorphisms by construction and the horizontal arrows are regular epimorphisms by induction hypothesis.\endproof

Special instances of the following theorem have been considered in literature: if $\EE$ is the category of \emph{groups} and $n=2$, the result can be deduced from Baues-Pirashvili \cite{BP} and Hartl-Vespa \cite{HVe}; if $\EE$ is a \emph{semi-abelian} category, under the identification of $n$-folded objects given in Remark \ref{nfoldedgroup} and with the kernel intersection definition of degree, the result has been announced by Manfred Hartl in several of his talks.

\begin{thm}\label{additivereflection}The full subcategory $\Fld^n(\EE)$ of $\,n$-folded objects of a $\,\sg$-pointed exact Mal'tsev category $\EE$ is a reflective Birkhoff subcategory. The associated Birkhoff reflection $J^n:\EE\to\Fld^n(\EE)$ is the universal endofunctor of $\,\EE$ of degree $\leq n$.\end{thm}
\proof Observe that the second assertion is a consequence of the first and of Proposition \ref{nadditive}a-b. In virtue of Proposition \ref{nfoldBirkhoff} it suffices thus to construct the reflection $J^n:\EE\to\Fld^n(\EE)$. The  latter is obtained by the following pushout
$$ \xymatrix@=10pt{
     X+\cdots+X  \ar@{->>}[rr]^{\theta_{X,\cdots,X}} \ar[dd]_{\delta^X_{n+1}} && P_{X,\dots,X} \ar[dd]^{\mu_nX}  \\
      && \\
     X \ar@{->>}[rr]_{\epsilon_nX}& & J^n(X)
                            }
  $$which is regular by Corollary \ref{stpushout} so that $J^n(X)=X/H_{n+1}[X]$ where $H_{n+1}[X]$ is the direct image of $R[\theta_{X,\dots,X}]$ under the folding map $\delta_{n+1}^X$. We will show that $J^n(X)$ is $n$-folded and that any morphism $X\to T$ with $n$-folded codomain $T$ factors uniquely through $J^n(X)$. For this, consider the following diagram
$$ \xymatrix@=1pt{
     H_{n+1}[X]+\cdots+H_{n+1}[X]\ar[dddddd]_{\delta^{H_{n+1}[X]}_{n+1}} \ar@<-1,ex>[rddd]_{p_0+\cdots+p_0}\ar@<+1,ex>[rddd]^>>>>{p_1+\cdots+p_1} \ar@{->>}[rr]^{\theta_{H_{n+1}[X],\dots,H_{n+1}[X]}} && P_{H_{n+1}[X],\dots,H_{n+1}[X]} \ar@<-1,ex>[ddd]_{P_{p_0,\cdots,p_0}}\ar@<+1,ex>[ddd]^>>>>{P_{p_1,\cdots,p_1}}\\
      &&&&\\
      &&&&\\
    & X+\cdots+X  \ar@{->>}[r]^{\theta_{X,\cdots,X}} \ar[ddd]_{\delta^X_{n+1}} & P_{X,\dots,X} \ar[ddd]^{\mu_nX} \ar@{->>}[drr]^>>>>{P_{\epsilon_nX,\dots,\epsilon_nX}}  \\
      &&&&& P_{J^n(X),\dots,J^n(X)}\ar@{.>}[ddlll]  \\
      &&&&\\
  H_{n+1}[X]  \ar@<-1,ex>[r]_{p_0}\ar@<+1,ex>[r]^{p_1} &  X \ar@{->>}[r]_{\epsilon_nX} &  J^n(X)
                            }
  $$in which the existence of the dotted arrow has to be shown. By Lemma \ref{reflcoeq} $P_{J^n(X),\dots,J^n(X)}$ is the coequaliser of the reflexive pair $(P_{p_0,\dots,p_0},P_{p_1,\dots,p_1})$. It suffices thus to check that $\mu_nX$ coequalises the same pair. This follows by precomposition with the regular epimorphism $\theta_{H_{n+1}[X],\dots,H_{n+1}[X]}$ using the commutativity of the previous diagram.

  For the universal property of $\epsilon_nX:X\onto J^n(X)$ let us consider a morphism $f:X\to T$ with $n$-folded codomain $T$. By construction of $J^n(X)$, the following commutative diagram
  $$ \xymatrix@=10pt{
       X+\cdots+X  \ar@{->>}[rr]^{\theta_{X,\dots,X}} \ar[dd]_{\delta^X_{n+1}}\ar[dr]^{f+\cdots+f} && P_{X,\dots,X} \ar[dr]^{P_{f,\dots,f}}  \\
       & T+\cdots+T \ar@{->>}[rr]^{\theta_{T,\dots,T}} \ar[dr]^{\delta^T_{n+1}}  &&  P_{T,\dots,T} \ar[dl]^{m_T}\\
       X \ar[rr]_{f}& & T
                              }
    $$induces the desired factorisation.\endproof

\begin{lma}\label{reflcoeq}
In a $\sg$-pointed exact Mal'tsev category, the functor $P_{X_1,\dots,X_{n+1}}$ preserves reflexive coequalisers in each of its $n+1$ variables.
\end{lma}
\proof By exactness, it suffices to show that $P$ preserves regular epimorphisms in each variable, and that for a regular epimorphism $f_i:X_i\onto X'_i$ the induced map on kernel relations $P_{X_1,\dots,R[f_i],\dots,X_{n+1}}\to R[P_{X_1,\dots,f_i,\dots,X_{n+1}}]$ is a regular epimorphism as well. By symmetry, it is even sufficient to do so for the first variable.

We shall argue by induction on $n$ (since for $n=1$ there is nothing to prove) and consider the following downward-oriented pullback diagram

$$
\xymatrix@=20pt{
          P_{X_1,\dots,X_n,X_{n+1}}
          \ar@{->>}[rrrr]^{}\ar[d]_{}&& && P^{+X_{n+1}}_{X_1,\dots,X_n} \ar[d] \\
          X_1+\cdots+X_{n }\ar@{->>}[rrrr]_{\theta_{X_1,\dots,X_n}}\ar@<-1ex>[u] && && P_{X_1,\dots,X_n} \ar@<-1ex>[u]
                           }
$$which derives from square (n) of the beginning of this section. By induction hypothesis, the two lower corners and the upper right corner are functors preserving regular epimorphisms in the first variable. It follows then from the cogluing lemma (cf. the proof of Theorem \ref{folklore}a) and Corollary \ref{regpush} that the upper left corner also preserves regular epimorphisms in the first variable.
%\cite[Lemma 2.5.7]{BB}

It remains to be shown that for $f:X_1\onto X'_1$ we get an induced regular epimorphism on kernel relations. For this we denote by $F,G,H$ the functors induced on the lower left, lower right and upper right corners, and consider the following commutative diagram

$$ \xymatrix@=12pt{
  &    P(R[f]) \ar[rrr]^{}  \ar@{.>>}[dr]^{\rho_P} \ar[ddl]_<<<<{} &&& H(R[f]) \ar[ddl]_<<<<{t_{R[f]}}  \ar@{->>}[dr]^<<<{\rho_H}\\
 &&   R[P(f)] \ar[ddl]_<<<<{} \ar[rrr]^<<<<<<{} &&&  R[H(f)] \ar[ddl]_{R(t_X)} \\
   F(R[f]) \ar@{->>}[dr]_<<<{\rho_F} \ar[rrr]_>>>>>>{\theta_{R[f]}}\ar@<-1ex>[ruu] &&&  G(R[f]) \ar@{->>}[dr]_<<<{\rho_G}\ar@<-1ex>[ruu] \\
  &   R[F(f)] \ar[rrr]_{R(\theta_X)} \ar@<-1ex>[ruu] &&& R[G(f)] \ar@<-1ex>[ruu]
                  }
$$
in which the back vertical square is a downward-oriented pullback by definition of $P$. By commutation of limits the front vertical square is a downward oriented pullback as well. Again, according to the cogluing lemma and Corollary \ref{regpush}, the induced arrow $\rho_P$ is then a regular epimorphisms, since $\rho_F$, $\rho_G$ and $\rho_H$ are so by induction hypothesis.\endproof

\subsection{Higgins commutator relations and their normalisation}\label{Higginscommutator}--\vspace{1ex}

We shall now concentrate on the case $X=X_1=X_2=\cdots=X_{n+1}$. Accordingly, we abbreviate the $n+1$ ``contractions'' as follows:$$\pi_{\hat{X}_i}=\pi_i:\overbrace{X+\cdots+X}^{n+1}\to \overbrace{X+\cdots+X}^n,\quad i=1,\dots,n+1.$$Proposition \ref{intersection} reads then $R[\theta_{X,\dots,X}]=R[\pi_1]\cap\cdots\cap R[\pi_{n+1}]$.

We denote the direct image of the kernel relation $R[\theta_{X,\dots,X}]$ under the folding map $\delta^X_{n+1}:X+\dots+X\to X$ by a single bracket $[\nabla_X,\dots,\nabla_X]$ of length $n+1$ and call it the \emph{$(n+1)$-ary Higgins commutator relation} on $X$.

The proof of Theorem \ref{additivereflection} shows that in $\sg$-pointed exact Mal'tsev categories the universal $n$-folded quotient $J^n(X)$ of an object $X$ is obtained by quotienting out the $(n+1)$-ary Higgins commutator relation. The binary Higgins commutator relation coincides with the Smith commutator $[\nabla_X,\nabla_X]$ (cf. Section \ref{Smithcommutator}, Corollary \ref{stpushout} and Proposition \ref{abelianization}) which ensures consistency of our notation.

Recall that in a pointed category the \emph{normalisation} of an effective equivalence relation $R$ on $X$ is the kernel of its quotient map $X\onto X/R$. In $\sg$-pointed exact Mal'tsev categories normalisation commutes with direct image, cf. Corollary \ref{stpushout}. In particular, the normalisation of the Higgins commutator relation yields precisely the Higgins commutator of same length, cf. Remark \ref{nfoldedgroup}.

%The direct image of the kernel $K[\theta_{X,\dots,X}]$ under the $(n+1)$st folding map $\delta^X_{n+1}$ is commonly denoted by a single bracket $[X,\dots,X]$ of length $n+1$ and called the \emph{Higgins commutator of $X$ of length $n+1$}, see \cite{HL,HV,Hi,MM} for further details. In $\sg$-pointed exact Mal'tsev categories, these Higgins commutators are thus the normalisations of the Higgins commutator relations of same length.

\begin{prp}\label{case1}
In a $\sg$-pointed exact Mal'tsev category, the image of $R[\theta_{X,X,X}]$ under $\delta^X_2 + 1_X$ is the kernel relation of the pushout of $\theta_{X,X,X}$ along $\delta^X_2+1_X$
$$
\xymatrix@=20pt{
          X+X+X \ar@{->>}[rr]^{\theta_{X,X,X}}\ar@{->>}[d]_{\delta^X_2 + 1_X }&& P_{X,X,X} \ar@{->>}[d] \\
          X+X\ar@{->>}[rr]^{\zeta^X_{X,X}}   && J^X_{X,X}
                           }
$$which may be computed as an intersection: $R[\zeta^X_{X,X}]=[R[\pi_1],R[\pi_1]]\cap R[\pi_2]$.\vspace{1ex}

In particular, we get the inclusion $[\nabla_X,[\nabla_X,\nabla_X]]\subset[\nabla_X,\nabla_X,\nabla_X]$.
\end{prp}
\proof
By Corollary \ref{stpushout}, the pushout is regular so that the first assertion follows from Proposition \ref{regpushout1}. Consider the following diagram
$$ \xymatrix@=12pt{
  &    X+X+X \ar@{->>}[rrr]^{\delta^X_2+1_X}  \ar@{->>}[dr]^>>>>{\theta^{+X}_{X,X}} \ar[ddl]_{\pi_3} &&& X+X \ar[ddl]_>>>>{\pi_2}  \ar@{->>}[dr]^{\eta^1(\pi_1)}\\
  &&   (X+X)\times_X(X+X) \ar[ddl] \ar@{->>}[rrr] &&&  I^1(\pi_1) \ar[ddl]_{f'}  \\
  X+X \ar@<-1ex>@{->>}[dr]_<<<{\theta_{X,X}} \ar@{->>}[rrr]^>>>>>>>>>>>>>>>>>>>>>>>{\delta^X_2}\ar@<-1ex>[ruu] &&&  X \ar@<-1ex>@{->>}[dr]^{\eta^1(X)}\ar@<-1ex>[ruu] \\
  &   X\times X \ar@{->>}[rrr]_{} \ar@<-1ex>[ruu]  &&& I^1(X) \ar@<-1ex>[ruu]
                   }
$$in which top and bottom are regular pushouts by Corollary \ref{stpushout}.
The bottom square constructs the associated abelian object $I^1(X)$ of $X$, while the top square constructs the associated abelian object $I^1(\pi_1)$ of $\pi_1:X+X\rightleftarrows X$ in the fibre over $X$. The upward oriented back and front faces are pushouts of split monomorphisms. The left face is a specialisation of square (2) just before Proposition \ref{totalkernel}. We can therefore apply Corollary \ref{B-C3} and we get diagram
$$ \xymatrix@=10pt{
  X+X+X\ar@{->>}[rrr]^{\delta^X_2+1_X}  \ar@{->>}[ddrr]^{\theta_{X,X,X}} \ar@<-1ex>[ddddr]_{\pi_3} &&& X+X \ar@<-1ex>[ddddr]_>>>>>>>>>>{\pi_2}  \ar@{->>}[ddrr]^{\zeta_{X,X}^X}\\
  &&&&&\\
  &&  P_{X,X,X} \ar@<1ex>[ddl]_<<<<{} \ar@{->>}[rrr]^<<<<<<{}  &&& J^X_{X,X} \ar@<1ex>[ddl]_<<<<{}  \\
  &&&&&\\
  &  X+X \ar@{->>}[rrr]_{\delta^X_2} \ar@<-2ex>[ruu] \ar[uuuul] &&& X \ar@<-2ex>[ruu] \ar[uuuul]
                   }
$$
in which the kernel relation of the regular epimorphism $\zeta^X_{X,X}$ is given by $$R[\eta^1(\pi_1)]\cap R[\pi_2]=[R[\pi_1],R[\pi_1]]\cap R[\pi_2].$$

For the second assertion, observe first that the ternary folding map $\delta_3$ may be identified with the composition $\delta^X_2\circ(\delta^X_2+1_X)$. Therefore, the ternary Higgins commutator relation $[\nabla_X,\nabla_X,\nabla_X]$ is the direct image under $\delta^X_2:X+X\to X$ of the kernel relation of $\zeta^X_{X,X}$. Now we have the following chain of inclusions, where for shortness we write $R_1=R[\pi_1],R_2=R[\pi_2],R_{12}=R_1\cap R_2$:

$$[R_{12},[R_{12},R_{12}]]\subset R_{12}\cap[R_{12},R_{12}]\subset R_2\cap[R_1,R_1].$$

By exactness, the direct image of the leftmost relation is $[\nabla_X,[\nabla_X,\nabla_X]]$, while the direct image of the right most relation is $[\nabla_X,\nabla_X,\nabla_X]$.\endproof

\begin{prp}\label{case2}
In a $\sg$-pointed exact Mal'tsev category, the image of $R[\theta_{X,\dots,X}]$ under $\delta^X_{n-1} + 1_X$ is the kernel relation of the pushout of $\,\theta_{X,\dots,X}$ along $\delta^X_{n-1} +1_X$
$$
\xymatrix@=20pt{
          X+\cdots+X \ar@{->>}[rr]^{\theta_{X,\dots,X}}\ar@{->>}[d]_{\delta^X_{n-1} + 1_X }&& P_{X,\dots,X} \ar@{->>}[d] \\
          X+X\ar@{->>}[rr]^{\zeta^X_{X,\dots,X}}   && J^X_{X,\dots,X}
                           }
$$which may be computed as an intersection $R[\zeta^X_{X,\dots,X}]=\overbrace{[R[\pi_1],\dots,R[\pi_1]]}^{n-1}\cap\, R[\pi_2]$.\vspace{1ex}

In particular, we get the inclusion $[\nabla_X,\overbrace{[\nabla_X,\dots,\nabla_X]}^{n-1}]\subset\overbrace{[\nabla_X,\dots,\nabla_X]}^n$.
\end{prp}
\proof
The first assertion follows from Proposition \ref{regpushout1} and the following diagram
$$ \xymatrix@=10pt{
  &  X+\cdots+X\ar@{->>}[rrr]^{\delta^X_{n-1}+1_X}  \ar@{->>}[dr]^>>>>{\theta^{+X}_{X,\dots,X}} \ar[ddl]_{\pi_n} &&& X+X \ar[ddl]_>>>>{\pi_2}  \ar@{->>}[dr]^{q}\\
  &&   P_{X,\dots,X}^{+X} \ar[ddl]_<<<<{} \ar@{->>}[rrr]^<<<<<<{}  &&&  \frac{X+X}{[\nabla_{\pi_1},\dots,\nabla_{\pi_1}]} \ar[ddl]  \\
  X+\cdots+X \ar@<-1ex>@{->>}[dr]_<<<{\theta_{X,\dots,X}} \ar@{->>}[rrr]^>>>>>>>>>>>>{\delta^X_{n-1}}\ar@<-1ex>[ruu] &&&  X \ar@<-1ex>@{->>}[dr]\ar@<-1ex>[ruu] \\
  &   P_{X,\dots,X} \ar@{->>}[rrr]_{} \ar@<-1ex>[ruu]  &&& X/[\nabla_X,\dots,\nabla_X] \ar@<-1ex>[ruu]
                   }
$$in which top and bottom are regular pushouts by Corollary \ref{stpushout}.
The bottom square constructs the quotient of $X$ by the $(n-1)$-ary Higgins commutator relation  $[\nabla_X,\dots,\nabla_X]$. The top square constructs the quotient of $\pi_1:X+X\rightleftarrows X$ by the $(n-1)$-ary Higgins commutator relation  $[\nabla_{\pi_1},\dots,\nabla_{\pi_1}]$ in the fibre over $X$. The upward oriented back and front faces are pushouts of split monomorphisms. The left face is a specialisation of square (n) of the beginning of this section. We can therefore apply Corollary \ref{B-C3} and we get the following diagram
$$ \xymatrix@=10pt{
  X+\cdots+X\ar@{->>}[rrr]^{\delta^X_{n-1}+1_X}  \ar@{->>}[ddrr]^{\theta_{X,\dots,X}} \ar@<-1ex>[ddddr]_{\pi_n} &&& X+X \ar@<-1ex>[ddddr]_<<<<{\pi_2}  \ar@{->>}[ddrr]^{\zeta_{X,\dots,X}^X}\\
  &&&&&\\
  &&  P_{X,\dots,X} \ar@<1ex>[ddl]_<<<<{} \ar@{->>}[rrr]^<<<<<< {} &&& J^X_{X,\dots,X} \ar@<1ex>[ddl]_<<<<{}  \\
  &&&&&\\
  &  X+\cdots+X \ar@{->>}[rrr]_{\delta^X_{n-1}} \ar@<-2ex>[ruu] \ar[uuuul] &&& X \ar@<-2ex>[ruu] \ar[uuuul]
                   }
$$
in which the kernel relation of the regular epimorphism $\zeta^X_{X,\dots,X}$ is given by $$R[q]\cap R[\pi_2]=[\nabla_{\pi_1},\dots,\nabla_{\pi_1}]\cap R[\pi_2]=[R[\pi_1],\dots, R[\pi_1]]\cap R[\pi_2].$$

Since $\delta^X_n=\delta^X_2\circ(\delta^X_{n-1}+1_X)$, the proof of the second assertion is completely analogous to the proof of the corresponding part of Proposition \ref{case1}.\endproof

The semi-abelian part of the following corollary can also be derived from a direct analysis of ``iterated'' Higgins commutators, cf. \cite[Proposition 2.21(iv)]{HV}.

\begin{cor}\label{inclusion}
For any object $X$ of a $\sg$-pointed exact Mal'tsev category, the iterated Smith commutator $[\nabla_X,[\nabla_X,[\nabla_X,\cdots,[\nabla_X,\nabla_X]\cdots]]]$ is a subobject of the Higgins commutator relation $[\nabla_X,\dots,\nabla_X]$ of same length.

In a semi-abelian category, the iterated Huq commutator $[X,[X,\cdots,[X,X]\cdots]]$ is a subobject of the Higgins commutator $[X,\dots,X]$ of same length.\end{cor}

\proof The first statement follows inductively from Propositions \ref{case1} and \ref{case2}.

The second statement follows from the first and the fact that in a semi-abelian category, the iterated Huq commutator is the normalisation of the iterated Smith commutator by Remark \ref{Huq} and \cite[Proposition 2.2]{GL}.\endproof

\begin{prp}\label{degn>nnil}In a $\sg$-pointed exact Mal'tsev category $\EE$, each $n$-folded object is an $n$-nilpotent object, i.e. $\Fld^n(\EE)\subset\Nil^n(\EE)$. In particular, endofunctors of degree $\leq n$ take values in $n$-nilpotent objects.\end{prp}

\proof The second assertion follows from the first and from Proposition \ref{nadditive}. For any $n$-folded object $X$, the $(n+1)$-ary Higgins commutator relation of $X$ is discrete and hence, by Corollary \ref{inclusion}, the iterated Smith commutator of same length is discrete as well. By an iterated application of Theorem \ref{noBirkhoff} and Proposition \ref{univcentral}, this iterated Smith commutator is the kernel relation of $\eta^n_X:X\onto I^n(X)$, and hence $X\cong I^n(X)$, i.e. $X$ is $n$-nilpotent.\endproof

The following theorem generalises Theorem \ref{2discrepancy} to all positive integers.

\begin{thm}\label{boumboum}For a $\sg$-pointed exact Mal'tsev category $\DD$ such that the identity functor of $\Nil^{n-1}(\DD)$ is of degree $\leq n-1$, the following properties are equivalent:

\begin{itemize}\item[(a)]all objects are $n$-nilpotent;\item[(b)]for all objects $X_1,\dots,X_n$, the map $\theta_{X_1,\dots,X_n}$ is an affine extension;\item[(c)]for all objects $X_1,\dots,X_n$, the map $\theta_{X_1,\dots,X_n}$ is a central extension.\end{itemize}\end{thm}
\proof For an $n$-nilpotent category $\DD$, the Birkhoff reflection $I^{n-1}:\DD\to\Nil^{n-1}(\DD)$ is a central reflection. Since all limits involved in the construction of $P_{X_1,\dots,X_n}$ are preserved under $I^{n-1}$ by an iterative application of Proposition \ref{Diana}, we get $I^{n-1}(\theta_{X_1,\dots,X_n})=\theta_{I^{n-1}(X_1),\dots,I^{n-1}(X_n)}$. Since by assumption the identity functor of $\Nil^{n-1}(\DD)$ is of degree $\leq n-1$, the latter comparison map is invertible so that by Theorem \ref{pointwiseaffine}, the comparison map $\theta_{X_1,\dots,X_n}$ is an affine extension, i.e. (a) implies (b). By Proposition \ref{central}, (b) implies (c).

Specializing (c) to the case $X=X_1=X_2=\cdots=X_n$ we get the following commutative diagram
$$ \xymatrix@=25pt{
  R[\theta_{X,\dots,X}]  \ar@<-1,ex>[r] \ar@<+1,ex>[r]\ar@{->>}[d] &   X+\cdots+X \ar@{->>}[rr]^{\theta_{X,\dots,X}} \ar@{->>}[d]_{\delta^X_n} && P_{X,\dots,X} \ar@{->>}[d] \\
   [\nabla_X,\dots,\nabla_X]\ar@<-1,ex>[r] \ar@<+1,ex>[r] &   X  \ar@{->>}[rr] && X/[\nabla_X,\dots,\nabla_X]         }
  $$in which is the right square is a regular pushout by Corollary \ref{stpushout} so that the lower row represents the kernel relation of a central extension. We have already seen that the iterated Smith commutator $[\nabla_X,[\nabla_X,[\nabla_X,\cdots,[\nabla_X,\nabla_X]\cdots]]]$ of length $n$ is the kernel relation of the unit $\eta^{n-1}(X):X\onto I^{n-1}(X)$ of the $(n-1)$st Birkhoff reflection. Corollary \ref{inclusion} implies thus that this unit is a central extension as well so that $\DD$ is $n$-nilpotent, i.e. (c) implies (a).\endproof

\begin{dfn}\label{P_n}A $\sg$-pointed category $(\DD,\star_\DD)$ with pullbacks is said to satisfy condition $P_n$ if for all $X_1,\dots,X_n,Z$, pointed cobase-change $(\alpha_Z)_!:\DD\to\Pt_Z(\DD)$ takes the object $P_{X_1,\dots,X_n}$ to the object $P^{+Z}_{X_1,\dots,X_n}$.\end{dfn}

In particular, since $P_{X}=\star$ condition $P_1$ is void and just expresses that $(\alpha_Z)_!$ preserves the null-object. Since $P_{X,Y}=X\times Y$ condition $P_2$ expresses that $(\alpha_Z)_!$ preserves binary products. Therefore, the following result extends Corollary \ref{abelian} ($n=1$) and Theorem \ref{quad} ($n=2$) to all positive integers.

\begin{prp}\label{ndeg}The identity functor of a $\sg$-pointed exact Mal'tsev category $\,\DD$ is of degree $\leq n$ if and only if all objects are $n$-nilpotent and the Birkhoff subcategories $\Nil^k(\DD)$ satisfy condition $P_k$ for $1\leq  k\leq n$.\end{prp}

\proof Since the statement is true for $n=1$ by Corollary \ref{abelian} we can argue by induction on $n$ and assume that the statement is true up to level $n-1$. In particular, we can assume that $\Nil^{n-1}(\DD)$ has an identity functor of degree $\leq n-1$. Let us then consider the following substitute of the main diagram \ref{maindiagram}:

$$ \xymatrix@=12pt{
     {(X_1+\cdots+X_n)\diamond Z\;} \ar@{>->}[r] \ar[d]_{\theta_{X_1,\dots,X_n}\diamond Z}  & (X_1+\cdots+X_n)+Z \ar@{}[rrd]|*+[o][F-]{a} \ar@{->>}[rr]^{\theta_{X_1+\cdots+X_n,Z}} \ar[d]^{\theta_{X_1,\dots,X_n}+ Z} && (X_1+\cdots+X_n)\times Z\ar[d]^{\theta_{X_1,\dots X_n}\times Z} \\
     {P_{X_1,\dots,X_n}\diamond Z\;} \ar@{>->}[r] \ar[d]_{\varphi^Z_{X_1,\dots,X_n}} & P_{X_1,\dots,X_n}+Z \ar@{}[rrd]|*+[o][F-]{b} \ar@{->>}[rr]^{\theta_{P_{X_1,\dots,X_n},Z}} \ar[d]^{\phi_{X_1,\dots,X_n}^Z} && P_{X_1,\dots,X_n} \times Z \ar[d]^{\mu_{X_1,\dots,X_n}^Z} \\
     {P^{\diamond Z}_{X_1,\dots,X_n}\;} \ar@{>->}[r] & P^{+Z}_{X_1,\dots,X_n} \ar@{->>}[rr] && P^{\times Z}_{X_1,\dots,X_n}
                       }
    $$
in which the composite vertical morphisms from left to right are respectively $$\theta_{X_1,\dots,X_n}^{\diamond Z}\text{ and }\theta_{X_1,\dots,X_n}^{+Z}\text{ and }\theta_{X_1,\dots,X_n}^{\times Z},$$ and the morphism $\mu_{X_1,\dots,X_n}^Z$ is the canonical isomorphism. Exactly as in the proof of Proposition \ref{quad1} it follows that the identity functor of $\DD$ is of degree $\leq n$ if and only if squares (a) and (b) are pullback squares.

Square (b) is a pullback if and only if $\phi_{X_1,\dots,X_n}^Z$ is invertible which is the case precisely when condition $P_n$ holds. By Proposition \ref{caraffine} and Theorem \ref{boumboum}, square (a) is a pullback if and only if $\theta_{X_1,\dots,X_n}$ is an affine (resp. central) extension, which is the case for all objects $X_1,\dots,X_n$ precisely when $\DD$ is $n$-nilpotent.\endproof

\begin{cor}\label{crosseffects}A semi-abelian category has an identity functor of degree $\leq n$ if and only if either of the following three equivalent conditions is satisfied:\begin{itemize}\item[(a)]all objects are $n$-nilpotent, and the comparison maps$$\varphi^Z_{X_1,\dots,X_n}:P_{X_1,\dots,X_n}\diamond Z\to P^{\diamond Z}_{X_1,\dots,X_n}$$ are invertible for all objects $X_1,\dots,X_n,Z$;\item[(b)]the $(n+1)$st cross-effects of the identity functor $$cr_{n+1}(X_1,\dots,X_n,Z)=K[\theta_{X_1,\dots,X_n,Z}]$$ vanish for all objects $X_1,\dots,X_n,Z$;\item[(c)]the co-smash product is of degree $\leq n-1$, i.e. the comparison maps$$\theta^{\diamond Z}_{X_1,\dots,X_n}:(X_1+\cdots+X_n)\diamond Z\to P^{\diamond Z}_{X_1,\dots,X_n}$$are invertible for all objects $X_1,\dots,X_n,Z.$\end{itemize}\end{cor}

\proof Condition (a) expresses that squares (a) and (b) of the main diagram are pullbacks. By Proposition \ref{ndeg} this amounts to an identity functor of degree $\leq n$.

For (b) note that by protomodularity the cross-effect $cr_{n+1}(X_1,\dots,X_n,Z)$ is trivial if and only if the regular epimorphism $\theta_{X_1,\dots,X_n,Z}$ is invertible.

The equivalence of conditions (b) and (c) follows from the isomorphism of kernels $K[\theta_{X_1,\dots,X_n,Z}]\cong K[\theta^{\diamond Z}_{X_1,\dots,X_n}]$. The latter is a consequence of the $3\times 3$-lemma which, applied to main diagram \ref{ndeg} and to square (n), yields a chain of isomorphisms:$$K[\theta_{X_1,\dots,X_n}^{\diamond Z}]\cong K[K[\theta_{X_1,\dots,X_n}^{+Z}]\onto K[\theta_{X_1,\dots,X_n}^{\times Z}]]\cong K[\theta_{X_1,\dots,X_n,Z}].$$\endproof

\subsection{Higher duality and multilinear cross-effects}--\vspace{1ex}

In Section \ref{linearquadratic} we obtained a precise criterion for when $2$-nilpotency implies quadra-ticity, namely algebraic distributivity, cf. Corollary \ref{parexemple}. We now look for a similar criterion for when $n$-nilpotency implies an identity functor of degree $\leq n$. Proposition \ref{ndeg} gives us an explicit exactness condition in terms of certain limit-preservation properties (called $P_n$) of pointed cobase-change along initial maps. In order to exploit the latter we first need to dualise condition $P_n$ into a colimit-preservation property, extending Proposition \ref{tau=iota}. Surprisingly, this dualisation process yields the simple condition that in each variable, the functor $P_{X_1,\dots,-,\dots,X_n}$ takes binary sums to binary sums in the fibre over $P_{X_1,\dots,\star,\dots,X_n}$. For $n$-nilpotent semi-abelian categories, this in turn amounts to the condition that the $n$-th cross-effect of the identity functor $cr_n(X_1,\dots,X_n)=K[\theta_{X_1,\dots,X_n}]$ is \emph{multilinear}.

Such a characterisation of degree $n$ functors in terms of the multilinearity of their $n$-th cross-effect is already present in the original treatment of Eilenberg-Mac Lane \cite{EM} for functors between abelian categories. It plays also an important role in Goodwillie's \cite{G2} homotopical context (where however linearity has a slightly different meaning). The following lemma is known in contexts close to our's.

\begin{lma}\label{diagonal}Let $\DD$ be a $\sg$-pointed category and let $\,\EE$ be an abelian category. Any multilinear functor $F:\DD^n\to\EE$ has a diagonal $G:\DD\overset{\Delta_\DD^n}{\lto}\DD^n\overset{F}{\lto}\EE$ of degree $\leq n$.\end{lma}

\proof This is a consequence of the decomposition formula of Eilenberg-Mac Lane for functors taking values in abelian categories, cf. Remark \ref{history}. Indeed, an induction on $k$ shows that the $k$-th cross-effect of the diagonal $cr_k^G(X_1,\dots,X_k)$ is the direct sum of all terms $F(X_{j_1},\dots,X_{j_n})$ such that the sequence $(j_1,\dots,j_n)$ contains only integers $1,2,\dots,k$, but each of them at least once. In particular, $$cr_n^G(X_1,\dots,X_n)\cong\bigoplus_{\sigma\in\Sigma_n}F(X_{\sigma(1)},\dots,X_{\sigma(n)})$$ and the cross-effects of $G$ of order $>n$ vanish, whence $G$ is of degree $\leq n$.\endproof

\begin{lma}[cf. Proposition 2.9 in \cite{HV}]\label{crosslemma}For each $n\geq 1$, the $n$-th cross-effect of the identity functor of a semi-abelian category preserves regular epimorphisms in each variable.\end{lma}

\proof The first cross-effect is the identity functor and the second cross-effect is the co-smash product. Proposition \ref{regpushout2} and Lemma \ref{kerpush} imply that the co-smash product preserves regular epimorphisms in both variables.

The general case $n+1$ follows from the already treated case $n=1$. By symmetry it suffices to establish the preservation property for the last variable which we shall denote $Z$. We have the following formula:$$cr_{n+1}(X_1,\dots,X_n,Z)=cr_n^{\diamond Z}(X_1,\dots,X_n)\quad\quad(n\geq 1)$$where $cr_n^{\diamond Z}(X_1,\dots,X_n)=K[\theta_{X_1,\dots,X_n}^{\diamond Z}]$ denotes the $n$-th cross-effect of the functor $-\diamond Z$. Indeed, this kernel has already been identified with $K[\theta_{X_1,\dots,X_n,Z}]$ in the proofs of Corollaries \ref{bilinear} and \ref{crosseffects}. It is now straightforward to deduce preservation of regular epimorphisms in $Z$ using that $(-)\diamond(-)$ preserves regular epimorphisms in both variables.\endproof

\begin{cor}[cf. Proposition 2.21 in \cite{HV}]\label{Higgins}In a semi-abelian category, the image of a Higgins commutator $[X,\dots,X]$ of $X$ under a regular epimorphism $f:X\onto Y$ is the corresponding Higgins commutator $[Y,\dots,Y]$ of $\,Y$.\end{cor}

\proof The Higgins commutator of length $n$ is the image of the diagonal $n$-th cross-effect $K[\theta_{X,\dots,X}]$ under the folding map $\delta^X_n:X+\cdots+X\to X$, cf. Section \ref{Higginscommutator}. By Lemma \ref{crosslemma}, any regular epimorphism $f:X\onto Y$ induces a regular epimorphism $K[\theta_{X,\dots,X}]\onto K[\theta_{Y,\dots,Y}]$ on diagonal cross-effects, whence the result.\endproof

Note that the commutative square (n) of the beginning of this section induces the following pullback square
$$
\xymatrix@=20pt{
        &  P_{X_1,\dots,X_{n+1}}\ar@{->}[rrrr]^{\chi_{X_1,\dots,X_{n+1}}}\ar[d]_{P_{X_1,\dots,X_n,\omega_{X_{n+1}}}}&& && P^{+X_{n+1}}_{X_1,\dots,X_n} \ar[d]_{P^{+\omega_{X_{n+1}}}_{X_1,\dots,X_n}} \\
     P_{X_1,\dots,X_n,\star} \ar@{=}[r]  &  X_1+\cdots+X_{n }\ar@{->}[rrrr]_{\theta_{X_1,\dots,X_n}}\ar@<-1ex>[u]_{P_{X_1,\dots,X_n,\alpha_{X_{n+1}}}} && && P_{X_1,\dots,X_n} \ar@<-1ex>[u]_{P^{+\alpha_{X_{n+1}}}_{X_1,\dots,X_n}}
                                }
$$in which the identification $P_{X_1,\dots,X_{n-1},\star}=X_1+\cdots+X_{n-1}$ is exploited to give the left vertical morphisms names. Recall that $\alpha_X:\star\to X$ and $\omega_X:X\to\star$ denote the initial and terminal maps.

\begin{prp}\label{higherduality}For any objects $X_1,\dots,X_{n-1},Y,Z$ of a $\sigma$-pointed category $(\DD,\star_\DD)$ with pullbacks consider the following diagram

$$ \xymatrix@=10pt{
       P_{X_1,\dots,X_{n-1},Y}+Z \ar@{->>}[dd]_{P_{X_1,\dots,X_{n-1},\omega_{Y}}+Z} \ar[rr]^{\rho_{X_{1},\dots,X_{n-1},Y,Z}}&& P_{X_1,\dots,X_{n-1},Y+Z} \ar@{->>}[dd]_{P_{X_{1},\dots,X_{n-1},\pi_{Z}^{Y}}}\\
        && \\
       X_1+\cdots+X_{n-1}+Z\ar[rr]_{\theta_{X_{1},\dots,X_{n-1},Z}} \ar@<-1ex>@{->}[uu]_{} && P_{X_1,\dots,X_{n-1},Z} \ar@<-1ex>@{->}[uu]_{P_{X_{1},\dots,X_{n-1},\iota_{Z}^{Y}}}
                           }
     $$
     in which the horizontal map $\rho_{X_1,\dots,X_{n-1},Y,Z}$ is induced by the pair $P_{X_{1},\dots,X_{n-1},\iota_{Y}^{Z}}:P_{X_1,\dots,X_{n-1},Y} \to P_{X_1,\dots,X_{n-1},Y+Z}$ and $P_{\alpha_{X_1},\dots,\alpha_{X_{n-1}},\iota_{Y}^{Z}}:Z\to P_{X_1,\dots,X_{n-1},Y+Z}$;
\begin{itemize}\item[(1)]The functor $P_{X_1,\dots,X_{n-1},-}:\DD\to\Pt_{X_1+\cdots+X_{n-1}}(\DD)$ preserves binary sums if and only if the upward-oriented square is a pushout for all objects $Y,Z$;
     \item[(2)]the category $\DD$ satisfies condition $P_n$ (cf. Definition \ref{P_n}) if and only if the downward-oriented square is a pullback for all objects $X_1,\dots,X_{n-1},Y,Z$.\end{itemize}

In particular, (1) and (2) hold simultaneously whenever $\theta_{X_1,\dots,X_{n-1},Z}$ is an affine extension for all objects $X_1,\dots,X_{n-1},Z$.\end{prp}

     \proof The second assertion follows from the discussion in Section \ref{pointfibration}. For (1), observe that the left upward-oriented square of the following diagram
     $$ \xymatrix@=10pt{
         P_{X_1,\dots,X_{n-1},Y}  \ar[rr]_{} \ar@<2ex>[rrrr]^{P_{X_{1},\dots,X_{n-1},\iota_{Y}^{Z}}} \ar@{->>}[dd]_{} &&  P_{X_1,\dots,X_{n-1},Y}+Z \ar@{->>}[dd]_{} \ar[rr]_<<{\rho_{X_{1},\dots,X_{n-1},Y,Z}} && P_{X_1,\dots,X_{n-1},Y+Z} \ar@{->>}[dd]_{P_{X_{1},\dots,X_{n-1},\pi_{Z}^Y}}\\
            &&  && \\
          X_1+\cdots+X_{n-1}\ar@<-1ex>[uu]_{P_{X_{1},\dots,X_{n-1},\alpha_Y}} \ar[rr]^{} \ar@<-2ex>[rrrr]_{P_{X_{1},\dots,X_{n-1},\alpha_Z}} &&  X_1+\cdots+X_{n-1}+Z\ar[rr]^{\theta_{X_1,\dots,Z}} \ar@<-1ex>[uu]_{} && P_{X_1,\dots,X_{n-1},Z} \ar@<-1ex>[uu]_{P_{X_{1},\dots, X_{n-1},\iota^{Y}_Z}}
                                }
          $$
         is a pushout so that the whole upward-oriented rectangle is a pushout if and only if the right upward-oriented square is a pushout, establishing (1).

\noindent For (2) observe that the right downward-oriented square of the following diagram
          $$ \xymatrix@=10pt{
                     P_{X_1,\dots X_{n-1},Y}+Z \ar@{->>}[dd]_{} \ar[rr]^{\rho_{X_{1},\dots,X_{n-1},Y,Z}}  && P_{X_1,\dots,X_{n-1},Y+Z} \ar@{->>}[dd]_{}
                             \ar[rr]^{\chi_{X_{1},\dots,X_{n-1},Y+Z}} && P^{+Y+Z}_{X_1,\dots,X_{n-1}} \ar@{->>}[dd]_{P^{+\pi^Y_{Z}}_{X_1,\dots, X_{n-1}}}\\
                            &&  && \\
                           X_1+\cdots+X_{n-1}+Z\ar[rr]^{\theta_{X_{1},\dots,X_{n-1},Z}} \ar@<-1ex>[uu]_{} \ar@<-2ex>[rrrr]_{\theta^{+Z}_{X_1,\cdots,X_{n-1}}} && P_{X_1,\dots,X_{n-1},Z} \ar@<-1ex>[uu]{}
                            \ar[rr]^{\chi_{X_{1},\dots,X_{n-1},Z}} && P^{+Z}_{X_1,\dots,X_{n-1}} \ar@<-1ex>[uu]_{P^{+\iota_Z^{Y}}_{X_1,\dots,X_{n-1}}}
                                                }
                    $$
          is a pullback (see below) so that the whole downward-oriented rectangle is a pullback if and only if the left downward-oriented square is a pullback. The whole downward-oriented rectangle is a pullback if and only if the comparaison map $$P_{X_1,\dots,X_{n-1},Y}+Z\to P^{+Z}_{X_1,\dots,X_{n-1},Y}$$ is invertible (i.e. if and only if condition $P_n$ holds) since the following square is by definition a pullback in the fibre $\Pt_Z(\DD)$:
          $$
          \xymatrix@=20pt{
P^{+Z}_{X_1,\dots,X_{n-1},Y}\ar@{->}[rrrr]^{\chi^{+Z}_{X_1,\dots,X_{n-1},Y}}\ar[d]_{P^{+Z}_{X_1,\dots,X_{n-1},\pi^Y_{Z}}}&& && P^{+Y+Z}_{X_1,\dots,X_{n-1}} \ar[d]_{P^{+\pi^Y_{Z}}_{X_1,\dots, X_{n-1}}} \\
               X_1+\cdots+X_{n-1}+Z\ar@{->}[rrrr]_{\theta^{+Z}_{X_1,\dots,X_{n-1}}}\ar@<-1ex>[u] && && P^{+Z}_{X_1,\dots,X_{n-1}} \ar@<-1ex>[u]
                                          }
          $$
Thus (2) is established. The pullback property of the right square above follows finally from considering the following diagram
         $$ \xymatrix@=10pt{
                    P_{X_1,\dots,X_{n-1},Y+Z} \ar@{->>}[dd]_{} \ar[rr]^{\chi_{X_{1},\dots,X_{n-1},Y+Z}} && P^{+Y+Z}_{X_1,\dots,X_{n-1}} \ar@{->>}[dd]_{P^{+\pi^Y_{Z}}_{X_{1},\dots,X_{n-1}}}\\
                          && \\
                         P_{X_1,\dots,X_{n-1},Z} \ar@<-1ex>[uu]^{}\ar@{->>}[dd]_{}
                           \ar[rr]_{\chi_{X_{1},\dots,X_{n-1},Z}} && P^{+Z}_{X_1,\dots,X_{n-1}} \ar@<-1ex>[uu]_{P^{+\iota_Z^{Y}}_{X_{1},\dots, X_{n-1}}}\ar@{->>}[dd]_{}\\
                           && \\
                        X_1+\cdots+X_{n-1} \ar[rr]_{\theta_{X_1,\dots,X_{n-1}}} \ar@<-1ex>[uu] && P_{X_1,\dots,X_{n-1}} \ar@<-1ex>[uu]
                                               }
                   $$in which whole rectangle and lower square are downward-oriented pullbacks.\endproof

\begin{thm}\label{multilinear}Let $\,\DD$ be an $n$-nilpotent $\sg$-pointed exact Mal'tsev category such that the identity functor of $\,\Nil^{n-1}(\DD)$ is of degree $\leq n-1$.\vspace{1ex}

\noindent The following properties are equivalent:
\begin{itemize}\item[(a)]the identity functor of $\,\DD$ is of degree $\leq n$;
\item[(b)]the category $\,\DD$ satisfies condition $P_n$ (cf. Definition \ref{P_n});
\item[(c)]the functor $P_{X_1,\dots,X_{n-1},-}:\DD\to\Pt_{X_1+\cdots+X_{n-1}}(\DD)$ preserves binary sums for all objects $X_1,\dots,X_{n-1}$.\end{itemize}

\noindent If $\,\DD$ is semi-abelian then the former properties are also equivalent to:
\begin{itemize}\item[(d)]the $n$-th cross-effect of the identity is coherent in each variable;
\item[(e)]the $n$-th cross-effect of the identity is linear in each variable;
\item[(f)]the diagonal $n$-th cross-effect of the identity is a functor of degree $\leq n$.
\end{itemize}\end{thm}

\proof It follows from Proposition \ref{ndeg} that properties (a) and (b) are equivalent, while properties (b) and (c) are equivalent by Theorem \ref{boumboum} and Proposition \ref{higherduality}. For the equivalence between (c) and (d), note first that the $n$-th cross-effect preserves regular epimorphisms in each variable by Lemma \ref{crosslemma} so that coherence (in the last variable) amounts to the property that the canonical map$$cr_n(X_1,\dots,X_{n-1},Y)+cr_n(X_1,\dots,X_{n-1},Z)\to cr_n(X_1,\dots,X_{n-1},Y+Z)$$ is a regular epimorphism. Since by Theorem \ref{boumboum} for $W=Y,Z,Y+Z$ the regular epimorphism $X_1+\cdots+X_{n-1}+W\onto P_{X_1,\dots,X_{n-1},W}$ is an affine extension, Proposition \ref{affinequotient} establishes the equivalence between (c) and (d). Finally, consider the following commutative diagram in $\Nil^1(\DD)=\Ab(\DD)$

$$ \xymatrix@=15pt{
          I^1(cr_n(X_{1\dots n-1},Y)+cr_n(X_{1\dots n-1},Z)) \ar[rr]^{\cong}\ar[d] && cr_n(X_{1\dots n-1},Y)\times cr_n(X_{1\dots n-1},Z) \\
          cr_n(X_{1\dots n-1},Y+Z) \ar[rr]_{cr_n(X_1,\dots,X_{n-1},\theta_{Y,Z})}   && cr_n(X_{1\dots n-1},Y\times Z)\ar[u]
                           }
     $$in which the upper horizontal map is invertible because the $n$-th cross-effect takes values in abelian group objects. It follows that the left vertical map is a section so that property (d) is equivalent to the invertibility of this left vertical map. Therefore, (d) is equivalent to the invertibility of the diagonal map $$cr_n(X_1,\dots,X_{n-1},Y+Z)\to cr_n(X_1,\dots,X_{n-1},Y)\times cr_n(X_1,\dots,X_{n-1},Z)$$ which expresses linearity in the last variable, i.e. property (e).

Property (e) implies property (f) by Lemma \ref{diagonal}. It suffices now to prove that (f) implies (a). The Higgins commutator $[X,\dots,X]$ of length $n$ is the image of diagonal $n$-th cross-effect $cr_n(X,\dots,X)$ under the $n$-th folding map $\delta_n^X:X+\cdots+X\to X$. The Higgins commutator of length $n$ is thus a quotient-functor of the diagonal $n$-th cross-effect and as such a functor of degree $\leq n$ by Theorem \ref{folklore}a. Corollary \ref{inclusion} and Remark \ref{Huq} imply that the kernel $K[\eta^{n-1}_X:X\onto I^{n-1}(X)]$ (considered as a functor in $X$) is a subfunctor of the Higgins commutator of length $n$ and hence, again by Theorem \ref{folklore}a, a functor of degree $\leq n$. It follows then from the short exact sequence of endofunctors$$\star\lto K[\eta^{n-1}]\lto id_\DD\lto I^{n-1}\lto\star$$ (by a third application of Theorem \ref{folklore}a) that the identity functor of $\DD$ is also of degree $\leq n$, whence (f) implies (a).\endproof

\subsection{Homogeneous nilpotency towers}--\vspace{1ex}

One of the starting points of this article was the existence of a functorial nilpotency tower for any $\sg$-pointed exact Mal'tsev category $\EE$, cf. Section \ref{Birkhoffcentral}. It is not surprising that for a semi-abelian category $\EE$ the successive kernels of the nilpotency tower capture the essence of the whole tower.\vspace{1ex}

To make this more precise, we denote by $$L_\EE(X)=\bigoplus_{n\geq 1} L_\EE^n(X)=\bigoplus_{n\geq 1}K[I^n(X)\onto I^{n-1}(X)]\in\Ab(\EE)$$ the graded abelian group object defined by the successive kernels. This construction is a functor in $X$. The \emph{nilpotency tower} of $\EE$ is said to be \emph{homogeneous} if for each $n$, the $n$-th kernel functor $L^n_\EE:\EE\to\Ab(\EE)$ is a functor of degree $\leq n$.

The degree of a functor does not change under composition with conservative left exact functors. We can therefore consider $L^n_\EE$ as an endofunctor of $\EE$. Observe also that the binary sum in $\Nil^n(\EE)$ is obtained as the reflection of the binary sum in $\EE$. This implies that the degree of $L^n_\EE$ is the same as the degree of $L^n_{\Nil^n(\EE)}$. We get the following short exact sequence of endofunctors of $\Nil^n(\EE)$ $$\star\lto L_{\Nil^n(\EE)}^n\lto id_{\Nil^n(\EE)}\lto I_\EE^{n,n-1}\lto\star$$where the last term is the relative Birkhoff reflection $I_\EE^{n,n-1}:\Nil^n(\EE)\to\Nil^{n-1}(\EE)$.

A more familiar way to express the successive kernels $L_\EE^n(X)$ of the nilpotency tower of $X$ is to realise them as subquotients of the lower central series of $X$. Indeed, the $3\times 3$-lemma implies that there is a short exact sequence$$\star\lto L_\EE^n(X)=\gamma_{n}(X)/\gamma_{n+1}(X)\lto X/\gamma_{n+1}(X)\lto X/\gamma_{n}(X)\lto\star$$where $\gamma_{n+1}(X)$ denotes the iterated Huq commutator of $X$ of length $n+1$, i.e. the kernel of the $n$-th Birkhoff reflection $\eta^n_X:X\onto I^n(X)$, cf. Remark \ref{Huq}.

The conclusion of the following theorem is folklore among those who are familiar with Goodwillie calculus in homotopy theory (cf. \cite{BD,G2}). Ideally, we would have liked to establish Theorem \ref{folklore}c by checking inductively one of the conditions of Theorem \ref{multilinear} without using any computation involving elements.

%\footnote{Statement (a) holds more generally for short exact sequences of functors taking values in a homological category with binary sums.}

\begin{thm}\label{folklore}Let $\EE$ be a semi-abelian category.\begin{itemize}\item[(a)]For any short exact sequence $\star\lto F_1\lto F\lto F_2\lto \star$ of endofunctors of $\,\EE$, $F$ is of degree $\leq n$ if and only if $F_1$ and $F_2$ are both of degree $\leq n$;\item[(b)]The nilpotency tower of $\;\EE$ is homogeneous if and only if the identity functors of $\,\Nil^n(\EE)$ are of degree $\leq n$ for all $n$;\item[(c)]The category of groups and the category of Lie algebras have homogeneous nilpotency towers.\end{itemize}\end{thm}

\proof For (a) we need the following \emph{cogluing lemma} for regular epimorphisms in regular categories: for any quotient-map of cospans
$$ \xymatrix@=25pt{ X \ar@{->>}[d]_f\ar@{->>}[r]&
        Z \ar@{->>}[d]_h\ar@{<-}[r] & Y \ar@{->>}[d]^g\\
    X'\ar@{->>}[r] &     Z' \ar@{<-}[r]\ & Y'
                       }
    $$in which the left naturality square is a \emph{regular pushout} (cf. Section \ref{regpushout0}), the induced map on pullbacks $f\times_hg:X\times_ZY\to X'\times_{Z'}Y'$ is again a regular epimorphism. Indeed, a diagram chase shows that the following square
    $$ \xymatrix@=20pt{
      X  \ar@{<-}[rr] \ar@{->>}[d] && X\times_ZY\ar@{.>>}[d] \\
   X'\times_{Z'}Z \ar@{<-}[rr]  && (X'\times_{Z'}Y')\times_{Y'}Y
                       }
    $$is a pullback. The left vertical map is a regular epimorphism by assumption so that the right vertical map is a regular epimorphism as well. Since $g$ is a regular epimorphism, the projection $(X'\times_{Z'}Y')\times_{Y'}Y\to X'\times_{Z'}Y'$ is again a regular epimorphism so that $f\times_hg$ is the composite of two regular epimorphisms.

The limit construction $P^F_{X_1,\dots,X_{n+1}}$ is an iterated pullback along \emph{split} epimorphisms. Therefore, Corollary \ref{regpush} and the cogluing lemma show inductively that the morphism $P^F_{X_1,\dots,X_{n+1}}\to P^{F_2}_{X_1,\dots,X_{n+1}}$ induced by the quotient-map $F\onto F_2$ is a regular epimorphism. The $3\times 3$-lemma yields then the exact $3\times 3$-square
$$
\xymatrix@=15pt{&\star\ar@{->}[d]&\star\ar@{->}[d]&\star\ar@{->}[d]&\\\star\;\ar@{->}[r]&cr_{n+1}^{F_1}(X_1,\dots,X_{n+1})\;\ar@{->}[r]\ar@{->}[d]&
cr_{n+1}^F(X_1,\dots,X_{n+1})\ar@{->}[r]\ar@{->}[d] & cr_{n+1}^{F_2}(X_1,\dots,X_{n+1}) \ar@{->}[r]\;\ar@{->}[d]&\star\\
\star\;\ar@{->}[r]&F_1(X_1+\cdots+X_{n+1})\ar@{->}[d]\;\ar@{->}[r]&F(X_1+\cdots+X_{n+1})\ar@{->}[r]\ar@{->}[d]&F_2(X_1+\cdots+X_{n+1})\ar@{->}[r]\ar@{->}[d]&\star\\
\star\ar@{->}[r]&P^{F_1}_{X_1,\dots,X_{n+1}}\ar@{->}[r]\ar@{->}[d]&P^{F}_{X_1,\dots,X_{n+1}}\ar@{->}[r]\ar@{->}[d]&
P^{F_2}_{X_1,\dots,X_{n+1}}\ar@{->}[r]\ar@{->}[d]&\star\\&\star&\star&\star&
              }
$$from which (a) immediately follows, cf. Corollary \ref{crosseffects}b.

For (b) we can assume inductively that $\Nil^{n-1}(\EE)$ has an identity functor of degree $\leq n-1$ so that the Birkhoff reflection $I^{n,n-1}_\EE:\Nil^n(\EE)\to\Nil^{n-1}(\EE)$ is of degree $\leq n-1$, and finally $I^{n,n-1}_\EE$ is also of degree $\leq n-1$ when considered as an endofunctor of $\Nil^n(\EE)$. Statement (b) follows then from (a) by induction on $n$.

For (c) we treat the group case, the Lie algebra case being very similar. In the category of groups, the graded object $L_\Gr(X)$ is a \emph{graded Lie ring} with Lie bracket $[-,-]:L^m_\Gr\otimes L^n_\Gr(X)\to L^{m+n}_\Gr(X)$ induced by the \emph{commutator map} $(x,y)\mapsto xyx^{-1}y^{-1}$ in $X$. This graded Lie ring is generated by its elements of degree $1$, cf. Lazard \cite[Section I.2]{La}. In particular, there is a regular epimorphism of abelian groups $L^1_\Gr(X)^{\otimes n}\to L^n_\Gr(X)$ which is natural in $X$. The functor which assigns to $X$ the tensor power $L^1_\Gr(X)^{\otimes n}$ is the diagonal of a multilinear abelian-group-valued functor in $n$ variables, and hence a functor of degree $\leq n$ by Lemma \ref{diagonal}. It follows from (a) that its quotient-functor $L^n_{\Gr}$ is of degree $\leq n$ as well, whence the homogeneity of the nilpotency tower in the category of groups.\endproof

\begin{thm}\label{stableHuq=Higgins}Let $\EE$ be a semi-abelian category.\vspace{1ex}

The following conditions are equivalent:
\begin{itemize}\item[(a)] The nilpotency tower of $\,\EE$ is homogeneous;
\item[(b)]For each $n$, the $n$-th Birkhoff reflection $I^n:\EE\to\Nil^n(\EE)$ is of degree $\leq n$;
\item[(c)]For each $n$, an object of $\,\EE$ is $n$-nilpotent if and only if it is $n$-folded;
\item[(d)]For each object $X$ of $\,\EE$, iterated Huq commutator $[X,[X,\cdots,X]\cdots]]$ and Higgins commutator $[X,X,\dots,X]$ of same length coincide.\end{itemize}

If $\,\EE$ satisfies one and hence all of these conditions, then so does any reflective Birkhoff subcategory of $\,\EE$. If $\,\EE$ is algebraically extensive and satisfies one and hence all of these conditions then so does the fibre $\Pt_X(\EE)$ over any object $X$ of $\,\EE$.\end{thm}

\proof We have already seen that (b) is equivalent to condition (b) of Theorem \ref{folklore}, which implies the equivalence between (a) and (b). Propositions \ref{nadditive} and \ref{degn>nnil} show that (b) implies (c), while Theorem \ref{additivereflection} shows that (c) implies (b). The equivalence between (c) and (d) is proved in exactly the same way as Corollary \ref{Huq=Higgins}.

Let $\DD$ be a reflective Birkhoff subcategory of $\EE$. We shall show that $\DD$ inherits (c) from $\EE$. By Proposition \ref{degn>nnil}, it suffices to show that in $\DD$ each $n$-nilpotent object $X$ is $n$-folded. Since the inclusion $\DD\inc\EE$ is left exact, it preserves $n$-nilpotent objects so that $X$ is $n$-nilpotent in $\EE$, and hence by assumption $n$-folded in $\EE$. The Birkhoff reflection $\EE\to\DD$ preserves sums and the limit construction $P_{X_1,\dots,X_{n+1}}$ by an iterated application of Proposition \ref{Diana}. Therefore, $X$ is indeed $n$-folded in $\DD$.

By Lemma \ref{semiabelian} algebraic extensivity implies that all pointed base-change functors are exact. By Lemma \ref{exact} this implies that the following square of functors
$$ \xymatrix@=10pt{
     \Pt_X(\EE)  \ar@{->}[rr]^{I^n_{\Pt_X(\EE)}} \ar[dd]_{\omega_X^*} && \Nil^n(\Pt_X(\EE)) \ar[dd]^{\omega_X^*}  \\
      && \\
     \EE\ar@{->}[rr]_{I^n_\EE}& & \Nil^n(\EE)
                            }
  $$
commutes up to isomorphism. The vertical functors are exact and conservative. Therefore, if $I^n_\EE$ is of degree $\leq n$ then $I^n_{\Pt_X(\EE)}$ is of degree $\leq n$ as well.\endproof

%It follows from Theorem \ref{folklore}b and Theorem \ref{additivereflection} that the nilpotency tower of a semi-abelian category is homogeneous if and only if $n$-nilpotency and $n$-foldedness are equivalent properties, or what amounts to the same (cf. the proof of Corollary \ref{Huq=Higgins}) if and only if iterated Huq commutator and Higgins commutator of same length coincide for each object. It would be highly desirable to have a categorical property responsible for this. It is certainly significative that both, the category of groups and the category of Lie algebras, are algebraically extensive.

\subsection{On Moufang loops and triality groups}\label{Moufang}--\vspace{1ex}

We end this article by giving an example of a semi-abelian category in which $2$-foldedness is not equivalent to $2$-nilpotency, namely the semi-abelian variety of Moufang loops. In particular, the semi-abelian subvariety of $2$-nilpotent Moufang loops is neither quadratic (cf. Proposition \ref{nadditive}) nor algebraically distributive (cf. Corollaries \ref{parexemple} and \ref{Huq=Higgins}). The nilpotency tower of the semi-abelian category of Moufang loops is thus inhomogeneous (cf. Theorem \ref{stableHuq=Higgins}). Nevertheless, the category of Moufang loops fully embeds into the category of \emph{triality groups} \cite{D,G,Ha} which, as we will see, is a semi-abelian category with homogeneous nilpotency tower.

Recall \cite{Br} that a \emph{loop} is a unital magma $(L,\cdot,1)$ such that left and right translation by any element $z\in L$ are bijective. A \emph{Moufang loop} \cite{Mou} is a loop $L$ such that $(z(xy))z=(zx)(yz)=z((xy)z)$ for all $x,y,z\in L$. Moufang loops form a semi-abelian variety which contains the variety of groups as a reflective Birkhoff subvariety. Moufang loops share many properties of groups, but the lack of a full associative law complicates the situation. The main example of a non-associative Moufang loop is the set of invertible elements of a non-associative \emph{alternative} algebra (i.e. in characteristic $\not=2$ a unital algebra in which the difference $(xy)z-x(yz)$ alternates in sign whenever two variables are permuted). In particular, the set $\OO^*$ of non-zero \emph{octonions} forms a Moufang loop. Taking the standard real basis of the octonions together with their additive inverses yields a Moufang subloop$$\OO_{16}=\{\pm 1,\pm e_1,\dots,\pm e_7\}$$with sixteen elements. We will see that $\OO_{16}$ is $2$-nilpotent, but not $2$-folded.

By Moufang's theorem \cite{Mou}, any Moufang loop, which can be generated by two elements, is associative and hence a group. In particular, for any element of a Moufang loop, left and right inverse coincide. The kernel of the reflection of a Moofang loop $L$ into the category of groups is the so-called \emph{associator subloop} $[L,L,L]_{ass}$ of $L$. For a Moufang loop $L$, the associator subloop is generated by the elements of the form $[x,y,z]=((xy)z)(x(yz))^{-1}$. Such an ``associator'' satisfies $[1,y,z]=[x,1,z]=[x,y,1]=1$ and is thus \emph{$3$-reducible}, cf. Remark \ref{nfoldedgroup}. This implies that for a Moufang loop $L$, the associator subloop $[L,L,L]_{ass}$ is contained in the ternary Higgins commutator $[L,L,L]$, cf. Proposition \ref{intersection} and Section \ref{Higginscommutator}. In conclusion, \emph{any $2$-folded Moufang loop has a trivial associator subloop and is therefore a $2$-folded group}. In particular, $\OO_{16}$ cannot be $2$-folded since $\OO_{16}$ is not a group. One can actually show that $[\OO_{16},\OO_{16},\OO_{16}]=\{\pm 1\}$. On the other hand, the centre of $\OO_{16}$ is also $\{\pm 1\}$, and the quotient by the centre $\OO_{16}/\{\pm 1\}$ is isomorphic to $(\ZZ/2\ZZ)^3$. This implies that $\OO_{16}$ is $2$-nilpotent, i.e. $[\OO_{16},[\OO_{16},\OO_{16}]]=\{1\}.$

The variety of Moufang loops is interesting with respect to the existence of centralisers. Since algebraic distributivity fails, such centralisers do not exist for general subloops, cf. \cite{BGr}. Nevertheless, each Moufang loop $L$ has a \emph{centre} $Z(L)$ in the sense of Section \ref{centralizer}, i.e. a centraliser $Z(1_L)$ for its identity $1_L:L\to L$. This centre $Z(L)$ is a \emph{normal} subloop of $L$, and is the intersection $Z(L)=M(L)\cap N(L)$ of the \emph{Moufang centre} $M(L)=\{z\in L\,|\, zx=xz\,\forall x\in L\}$ with the so-called \emph{nucleus} $N(L)=\{z\in L\,|\,[z,x,y]=[x,z,y]=[x,y,z]=1\,\forall x,y\in L\}$, cf. Bruck \cite{Br}.

Groups with triality have been introduced in the context of Moufang loops by Glauberman \cite{G} and Doro \cite{D}. A \emph{triality} on a \emph{group} $G_0$ is an action (by automorphisms) of the symmetric group $\Sg_3$ on three letters such that for all $g\in G_0$ and $\sg\in\Sg_3$ (resp. $\rho\in\Sg_3$) of order $2$ (resp. $3$), the identity$$[\sg,g](\rho.[\sg,g])(\rho^2.[\sg,g])=1$$holds where $[\sg,g]=(\sg.g)g^{-1}$. We denote the split epimorphism associated to the group action by $p:G_0\rtimes\Sg_3\lrto\Sg_3:i$ and call it the associated \emph{triality group}. The defining relations for a group with triality are equivalent to the following condition on the associated triality group $p:G\lrto\Sg_3:i$ (cf. Liebeck \cite{Li} and Hall \cite{Ha}):\begin{center}\emph{for any two special elements $g,h\in G$ such that $p(g)\not=p(h)$ one has $(gh)^3=1$}\end{center}where $g\in G$ is called \emph{special} if $g$ is conjugate in $G$ to some element of order $2$ in $i(\Sg_3)$. For the obvious notion of morphism, the category $\TG$ of triality groups is a full subcategory of the fibre $\Pt_{\Sg_3}(\Gr)$ over the symmetric group $\Sg_3$.

The category $\TG$ is closed under taking subobjects, products and quotients in $\Pt_{\Sg_3}(\Gr)$. Moreover, quotienting out the normal subgroup generated by the products $(gh)^3$ for all pairs of special elements $(g,h)$ such that $p(g)\not=p(h)$ defines a reflection $\Pt_{\Sg_3}(\Gr)\to\TG$. Therefore, $\TG$ is a reflective Birkhoff subcategory of $\Pt_{\Sg_3}(\Gr)$. Since the category of groups is an algebraically extensive semi-abelian category (cf. Section \ref{extensive}) with homogeneous nilpotency tower (cf. Theorem \ref{folklore}), so is its fibre $\Pt_{\Sg_3}(\Gr)$ by Lemma \ref{fibreextensive} and Theorem \ref{stableHuq=Higgins}.  The reflective Birkhoff subcategory $\TG$ formed by the triality groups is thus also a semi-abelian category with homogeneous nilpotency tower, again by Theorem \ref{stableHuq=Higgins}.

This result is remarkable because the category of triality groups contains the category of Moufang loops as a full \emph{coreflective} subcategory, and the latter has an inhomogeneous nilpotency tower. The embedding of Moufang loops and its right adjoint have been described by Doro \cite{D} for groups with triality, and by Hall \cite{Ha} for the associated triality groups, see also Grishkov-Zavarnitsine \cite{GZ}. Moufang loops can thus up to equivalence of categories be identified with triality groups for which the counit of the adjunction is invertible. Considering them inside the category of triality groups permits the construction of a homogeneous nilpotency tower.

\section*{Acknowledgements}We are grateful to Georg Biedermann, Rosona Eldred, Marino Gran, James Gray, Jonathan Hall, George Janelidze, Daniel Tanr\'e and Tim Van der Linden for helpful discussions. We are also grateful to the referee for his careful reading of our manuscript. Special thanks are due to Manfred Hartl whose seminar talk in September 2013 in Nice was the starting point for this work. The first author acknowledges financial support of the French ANR grant HOGT.


\begin{thebibliography}{99}

%\bibitem{BGV}M. Barr, P. A. Grillet and D. H. van Osdol -- \emph{Exact categories and categories of sheaves}, Lect. Notes Math. \textbf{236}, 1971.

\bibitem{BP}H.-J. Baues and T. Pirashvili -- \emph{Quadratic endofunctors of the category of groups}, Adv. Math. \textbf{141} (1999), 167--206.

%\bibitem{Be}C. Berger -- \emph{Algebraic and homotopical nilpotency}, talk at CT2015, available at http://math.unice.fr/$\sim$cberger/CT2015.pdf.

\bibitem{BG}I. Berstein and and T. Ganea -- \emph{Homotopical nilpotency}, Illinois J. Math. \textbf{5}  (1961), 99--130.

\bibitem{BD}G. Biedermann and B. Dwyer -- \emph{Homotopy nilpotent groups}, Algebr. Geom. Topol. \textbf{10} (2010), 33--61.

\bibitem{BB}F. Borceux and D. Bourn -- \emph{Mal'cev, protomodular, homological and semi-abelian categories}, Math. Appl. \textbf{566}, Kluwer Acad. Publ., 2004.

\bibitem{Bourn}D. Bourn -- \emph{Normalization equivalence, kernel equivalence and affine categories}, Lect. Notes Math. \textbf{1488}, Springer Verlag 1991, 43--62.

\bibitem{Bourn0}D. Bourn -- \emph{Mal'tsev categories and fibration of pointed objects}, Appl. Categ. Struct. \textbf{4} (1996), 307--327.

\bibitem{Bourn3}D. Bourn -- \emph{The denormalized $3\times 3$ lemma}, J. Pure Appl. Algebra \textbf{177} (2003), 113--129.

\bibitem{Bourn4}D. Bourn -- \emph{Commutator theory in regular Mal'tsev categories}, AMS Fields Inst. Commun. \textbf{43} (2004), 61--75.

\bibitem{Bourn5}D. Bourn -- \emph{Commutator theory in strongly protomodular categories}, Theory Appl. Categ. \textbf{13} (2004), 27--40.

\bibitem{Bourn6}D. Bourn -- \emph{On the monad of internal groupoids}, Theory Appl.Categ. \textbf{28} (2013), 150--165.

\bibitem{BG2}D. Bourn and M. Gran -- \emph{Central extensions in semi-abelian categories}, J. Pure Appl. Algebra \textbf{175} (2002), 31--44.

\bibitem{BG3}D. Bourn and M. Gran -- \emph{Centrality and connectors in Maltsev categories}, Algebra Universalis \textbf{48} (2002), 309--331.

\bibitem{BGr}D. Bourn and J.R.A. Gray -- \emph{Aspects of algebraic exponentiation}, Bull. Belg. Math. Soc. \textbf{19} (2012), 823--846.

\bibitem{BR}D. Bourn and D. Rodelo -- \emph{Comprehensive factorization and $I$-central extensions}, J. Pure Appl. Algebra \textbf{216} (2012), 598--617.

\bibitem{Br}R. H. Bruck -- \emph{A survey of binary systems}, Ergebnisse der Mathematik und ihrer Grenzgebiete \textbf{20}, Springer Verlag 1958.

\bibitem{CJ}A. Carboni and G. Janelidze -- \emph{Smash product of pointed objects in lextensive categories}, J. Pure Appl. Algebra \textbf{183} (2003), 27--43.

\bibitem{CKP}A. Carboni, G. M. Kelly and M. C. Pedicchio -- \emph{Some remarks on Mal'tsev and Goursat categories}, Appl. Categ. Struct. \textbf{1} (1993), 385--421.

\bibitem{CLW}A. Carboni, S. Lack and R. F. C. Walters -- \emph{Introduction to extensive and distributive categories}, J. Pure Appl. Algebra \textbf{84} (1993), 145--158.

\bibitem{CLP}A. Carboni, J. Lambek, M. C. Pedicchio, -- \emph{Diagram chasing in Malcev categories}, J. Pure Appl. Algebra \textbf{69} (1991), 271--284.

\bibitem{CGV}A. Cigoli, J. R. A. Gray, T. Van der Linden -- \emph{Algebraically coherent categories}, Theory Appl. Categ. \textbf{30} (2015), 1864--1905.

\bibitem{CSV}C. Costoya, J. Scherer, A. Viruel -- \emph{A torus theorem for homotopy nilpotent groups}, arXiv:1504.06100.

\bibitem{D} S. Doro -- \emph{Simple Moufang loops}, Math. Proc. Cambridge Philos. Soc. \textbf{83} (1978), 377--392.

\bibitem{EM}S. Eilenberg and S. Mac Lane -- \emph{On the groups H($\pi$,n). II. Methods of computation}, Ann. of Math. (2) \textbf{60} (1954), 49--139.

\bibitem{E}R. Eldred -- \emph{Goodwillie calculus via adjunction and LS cocategory}, arXiv:1209.2384.

\bibitem{EV}T. Everaert and T. Van der Linden -- \emph{Baer invariants in semi-abelian categories I: general theory}, Theory Appl. Categ. \textbf{12} (2004), 1--33.

\bibitem{EV2}T. Everaert and T. Van der Linden -- \emph{A note on double central extensions in exact Mal'tsev categories}, Cah. Topol. G\'eom. Differ. Cat\'eg. \textbf{51} (2010), 143--153.

\bibitem{FM}R. S. Freese and R. N. McKenzie -- \emph{Commutator theory for congruence modular varieties}, London Math. Soc. Lect. Note Series \textbf{125}, Cambridge Univ. Press, Cambridge, 1987.

\bibitem{G}G. Glauberman -- \emph{On loops of odd order II}, J. Algebra \textbf{8} (1968), 383--414.

\bibitem{G2}T. G. Goodwillie -- \emph{Calculus III. Taylor series}, Geom. Topol. \textbf{7} (2003), 645--711.

\bibitem{Gran1}M. Gran -- \emph{Central extensions and internal groupoids in Maltsev categories}, J. Pure Appl. Alg. \textbf{155} (2001), 139--166.

\bibitem{Gran2}M. Gran -- \emph{Applications of categorical Galois theory in universal algebra}, AMS Fields Inst. Commun. \textbf{43} (2004), 243--280.

%\bibitem{GKV}M. Gran, G. Kadjo and J. Vercruysse -- \emph{A torsion theory in the category of cocommutative Hopf algebras}, arXiv:1502.03130.

\bibitem{GR}M. Gran and D. Rodelo -- \emph{Beck-Chevalley condition and Goursat categories}, arXiv:1512.04066.

\bibitem{GL}M. Gran and T. Van der Linden -- \emph{On the second cohomology group in semi-abelian categories}, J. Pure Appl. Algebra \textbf{212} (2008), 636--651.

\bibitem{GZ}A. N. Grishkov and A. V. Zavarnitsine -- \emph{Groups with triality}, J. Algebra Appl. \textbf{5} (2006), 441--463.

\bibitem{Gr}J. R. A. Gray -- \emph{Algebraic exponentiation in general categories}, Appl. Categ. Struct. \textbf{20} (2012), 543--567.

\bibitem{Gr1}J. R. A. Gray -- \emph{Algebraic exponentiation for categories of Lie algebras}, J. Pure Appl. Algebra \textbf{216} (2012), 1964--1967.

\bibitem{Ha}J. I. Hall -- \emph{Central automorphisms, $Z^*$-theorems, and loop structures}, Quasigroups and Related Systems \textbf{19} (2011), 69--108.

\bibitem{HL}M. Hartl and B. Loiseau -- \emph{On actions and strict actions in homological categories}, Theory Appl. Categ. \textbf{27} (2013), 347--392.

\bibitem{HV}M. Hartl and T. Van der Linden -- \emph{The ternary commutator obstruction for internal crossed modules}, Adv. Math. \textbf{232} (2013), 571--607.

\bibitem{HVe}M. Hartl and C. Vespa -- \emph{Quadratic functors on pointed categories}, Adv. Math. \textbf{226} (2011), 3927--4010.

\bibitem{Hi}P. J. Higgins -- \emph{Groups with multiple operators}, Proc. London Math. Soc. \textbf{6} (1956), 366--416.

\bibitem{Ho}M. Hovey -- \emph{Lusternik-Schnirelmann cocategory}, Illinois J. Math. \textbf{37} (1993), 224--239.

\bibitem{Huq}S. A. Huq -- \emph{Commutator, nilpotency and solvability in categories}, Quart. J. Math. Oxford \textbf{19} (1968), 363--389.

\bibitem{JK}G. Janelidze and G. M. Kelly -- \emph{Galois theory and a general notion of central extension}, J. Pure Appl. Algebra \textbf{97} (1994), 135--161.

\bibitem{JK2}G. Janelidze and G. M. Kelly -- \emph{Central extensions in universal algebra: a unification of three notions}, Algebra Universalis \textbf{44} (2000), 123--128.

\bibitem{JK3}G. Janelidze and G. M. Kelly -- \emph{Central extensions in Mal'tsev varieties}, Theory Appl. Categ. \textbf{7} (2000), 219--226.

\bibitem{JMT}G. Janelidze, L. M\'arki and W. Tholen -- \emph{Semi-abelian categories}, J. Pure Appl. Algebra \textbf{168} (2002), 367--386.

\bibitem{JST}G. Janelidze, M. Sobral and W. Tholen -- \emph{Beyond Barr exactness: effective descent morphisms},  Cambridge Univ. Press, Encycl. Math. Appl. \textbf{97} (2004), 359--405.

\bibitem{JP}M. Jibladze and T. Pirashvili -- \emph{Linear extensions and nilpotence of Maltsev theories}, Contributions to Algebra and Geometry \textbf{46} (2005), 71--102.

\bibitem{JMc}B. Johnson and R. McCarthy -- \emph{A classification of degree $n$ functors I/II}, Cah. Topol. G\'eom. Diff\'er. Cat\'eg. \textbf{44} (2003), 2--38, 163--216.

%\bibitem{Lack}S. Lack -- \emph{The 3-by-3 lemma for regular Goursat categories}, Homology, Homotopy Appl. \textbf{6} (2004), 1--3.

\bibitem{La}M. Lazard -- \emph{Sur les groupes nilpotents et les anneaux de Lie}, Ann. Sci. E. N. S. \textbf{71} (1954), 101--190.

\bibitem{Li}M. W. Liebeck -- \emph{The classification of finite simple Moufang loops}, Math. Proc. Cambridge Philos. Soc. \textbf{102} (1987), 33--47.

\bibitem{Ma}A. I. Mal'cev -- \emph{On the general theory of algebraic systems}, Mat. Sbornik N. S. \textbf{35} (1954), 3--20.

\bibitem{MM}S. Mantovani and G. Metere -- \emph{Normalities and commutators}, J. of Algebra \textbf{324} (2010), 2568--2588.

\bibitem{Mos}J. Mostovoy -- \emph{Nilpotency and dimension series for loops}, Comm. Algebra \textbf{36} (2008), 1565--1579.

\bibitem{Mou}R. Moufang -- \emph{Zur Struktur von Alternativk\"orpern}, Math. Ann. \textbf{110} (1935), 416--430.

\bibitem{Pe}M. C. Pedicchio -- \emph{A categorical approach to commutator theory}, J. of Algebra \textbf{177} (1995), 647--657.

\bibitem{Pen}J. Penon -- \emph{Sur les quasi-topos}, Cah. Topol. G\'eom. Diff. \textbf{18} (1977), 181--218.

\bibitem{Q}D. Quillen -- \emph{Homotopical algebra},  Lect. Notes Math. \textbf{43}, Springer Verlag 1967.

\bibitem{Sm}J. D. H. Smith -- \emph{Mal'cev varieties}, Lect. Notes Math. \textbf{554}, Springer Verlag 1976.

\bibitem{SV}D. Stanovsky and P. Vojt\v{e}chovsk\`y, \emph{Commutator theory for loops}, J. of Algebra \textbf{399} (2014), 290--322.

\bibitem{Tim}T. Van der Linden -- \emph{Simplicial homotopy in semi-abelian categories}, J. K-theory \textbf{4} (2009), 379--390.

\end{thebibliography}
\end{document}